\documentclass[11pt,a4paper]{article}
\usepackage[utf8]{inputenc}
\usepackage[T1]{fontenc}
\usepackage{lmodern}
\usepackage{amsmath,amssymb,amsthm}
\usepackage[margin=2.6cm]{geometry}
\usepackage{enumitem}
\usepackage{xcolor}

\usepackage{tikz}
\usetikzlibrary{calc}
\usepackage[colorlinks=true,linkcolor=blue!50!black,
            citecolor=blue!50!black,urlcolor=blue!50!black]{hyperref}

\newcommand{\R}{\mathbb{R}}
\newcommand{\F}{\mathcal{F}}
\newcommand{\M}{\mathcal{M}}
\newcommand{\aff}{\operatorname{aff}}
\newcommand{\conv}{\operatorname{conv}}
\newcommand{\cone}{\operatorname{cone}}
\newcommand{\lin}{\operatorname{lin}}
\newcommand{\relint}{\operatorname{relint}}
\newcommand{\inte}{\operatorname{int}}
\newcommand{\Int}{\operatorname{Int}}
\newcommand{\rk}{\operatorname{rk}}
\newcommand{\hup}[1]{#1^{\uparrow}}
\newcommand{\hdn}[1]{#1^{\downarrow}}

\theoremstyle{plain}
\newtheorem{theorem}{Theorem}[section]
\newtheorem{lemma}[theorem]{Lemma}
\newtheorem{proposition}[theorem]{Proposition}
\newtheorem{corollary}[theorem]{Corollary}
\newtheorem{conjecture}[theorem]{Conjecture}

\theoremstyle{definition}
\newtheorem{definition}[theorem]{Definition}
\newtheorem{example}[theorem]{Example}
\newtheorem{remark}[theorem]{Remark}

\begin{document}

\begin{center}
\LARGE{\bf Face lattices of polytopes and bipolar polytopes}\\
\vspace{0.5 cm}
\normalsize{Jan Pavl\'ik} \\
\vspace{0.5 cm}
\small Faculty of Mechanical Engineering, Brno University of Technology,\\
Technick\'a 2896/2, 616 69 Brno, Czech Republic\\
{\tt pavlik@fme.vutbr.cz}\\
\end{center}

\begin{abstract}
The intervals of the face lattice of a convex polytope $Q$, ordered by inclusion, form an
abstract polytope, as follows dually from a theorem of Gleason and Hubard. Lindstr\"om asked in 1971 whether it is always the face lattice of a
convex polytope, an \emph{interval polytope} of $Q$, whose dual is an antiprism; Dobbins
showed that the answer is negative in dimension four. We study these objects axiomatically
as \emph{bipolar polytopes}: polytopes with two distinguished vertices, the poles, such
that every facet contains exactly one of them and an intersection property holds. The
vertices of every bipolar polytope carry a lattice order whose intervals are the faces,
and this lattice is the face lattice of the vertex figure at either pole; conversely, a
finite lattice is the vertex lattice of a bipolar polytope if and only if it is
atomistic, coatomistic and an abstract polytope. We relate the interval condition, under
which the polar of the Cayley polytope of $Q$ and $Q^\circ$ realizes the interval
polytope, to perfectly centered polytopes, recovering results of Broadie, of Fukuda and
Weibel, and of Dobbins. Finally, every combinatorial type of $3$-polytope has a
realization that decomposes into convex bipolar polytopes, one for each vertex.
\end{abstract}


\section{Introduction}\label{sec:intro}

The face lattice of a polytope provides a complete combinatorial description of its
structure, encoding the incidence relations among faces of all dimensions. While this
representation is fundamental in polytope theory, it is natural to ask whether the
same information can be expressed through the vertices and faces of another polytope.
In particular, we investigate representations in which the elements of a face lattice
correspond to vertices of a higher-dimensional polytope equipped with two
distinguished vertices, called poles. The poles provide a geometric interpretation of
the least and greatest elements of the lattice, while the remaining vertices represent
its intermediate elements. Such a representation establishes a connection between
lattice-theoretic and geometric structures: the order relation can be recovered from
the incidence of faces containing either pole, and the original face lattice is
closely related to the vertex figure at a pole. This approach offers a geometric
interpretation of abstract order relations and provides a natural setting for studying
duality, intervals, and incidence properties of face lattices.

A further motivation arises from the problem of decomposing convex polytopes into
smaller cells with prescribed combinatorial properties. In particular, cubical
subdivisions \cite{BuchstaberPanov,Eppstein} of three-dimensional polytopes suggest associating a cell with each
vertex of the original polytope, using an interior point together with suitably chosen
points on its edges and faces. For a vertex of order $3$, the resulting cell is
combinatorially a cube; in a more general case, it has the structure of a
trapezohedron whose two opposite vertices correspond to the original vertex and the
common interior point. These cells naturally lead to bipolar polytopes, and their
combinatorial structure reflects the vertex figure of the original polytope. This
raises the question of when such decompositions can be realized with convex cells and
how far the construction extends beyond three dimensions. Besides their intrinsic
geometric interest, these questions are related to polyhedral subdivisions, structured
mesh generation, and the construction of cell complexes with controlled local
combinatorics.

More broadly, the study of bipolar polytopes is motivated by the interaction between
the combinatorics of face lattices and the geometry of convex realizations. Every
interval of the face lattice of a polytope is itself the face lattice of a
polytope, suggesting that a suitable
geometric representation should preserve this hereditary structure. Bipolar polytopes
provide a framework in which lattice intervals can be interpreted as faces, with their
endpoints serving as local poles. The symmetry between the two poles also reflects
order duality and connects the construction to the duality theory of polytopes. An
important question is therefore whether arbitrary abstractly polytopal lattices
(Definition~\ref{def:apl}) admit bipolar representations and, more strongly, whether these
representations can always be realized by convex polytopes; the latter question goes
back to Lindstr\"om and was answered negatively by Dobbins (see the related work below). Understanding this distinction between combinatorial
existence and geometric realizability may provide new tools for constructing polytopes
with prescribed incidence structures and for investigating the constraints imposed by
convexity.

\medskip
\medskip
\noindent\textbf{Main results.}
\begin{description}[leftmargin=1em,style=nextline]
\item[Theorem A (structure, Section~\ref{sec:bipolar}).]
The vertices of a bipolar polytope carry a lattice order whose intervals are exactly
its nonempty faces; this lattice is the face lattice of the vertex figure at the south
pole.

\item[Theorem B (characterization, Section~\ref{sec:char}).]
A finite lattice is abstractly polytopal, i.e.\ atomistic, coatomistic and an abstract
polytope (Definition~\ref{def:apl}), if and only if it is the vertex lattice of a bipolar
polytope. In particular, the faces of every convex polytope $Q$ are the vertices
of a bipolar polytope $W(\F(Q))$, whose faces are the intervals of the face lattice
$\F(Q)$.

\item[Theorem C (decomposition in $\R^3$, Section~\ref{sec:geom}).]
Every $3$-polytope with a regular configuration is tiled by convex bipolar polytopes,
one for each vertex, with a common pole. Every combinatorial type of $3$-polytope has
such a realization.
\end{description}

\noindent In more detail: in Theorem~A the order is $x\le y$ iff $y$ lies in the
smallest north face containing $x$, the facets are the intervals $[a,\nu]$ and
$[\sigma,c]$ for the atoms $a$ and coatoms $c$, and every face is again bipolar.
Theorem~B is constructive: the bipolar polytope $W(L)$ has the elements of $L$ as
vertices and the principal filters of the atoms and principal ideals of the coatoms as
facets; it is the poset of intervals of $L$ considered by Lindstr\"om \cite{Lindstrom},
and its dual is the antiprism of $L$ in the sense of Gleason and Hubard
\cite{GleasonHubard}. For $W(L)$ itself, the description of its faces as intervals and of
the vertex figures at the poles appears, in dual form, in \cite[Section~2]{GleasonHubard}.
That $W(\F(Q))$ has a convex realization whenever $Q$ is perfectly
centered, namely the polar of the Cayley polytope $\mathrm{Cay}(Q,Q^\circ)$, is, in polar
form, a theorem of Broadie \cite{Broadie} (see also \cite[Lemma~2.1]{Lawrence}); in our terms it is
Theorem~\ref{thm:cayley}. Every combinatorial type of dimension
at most three has a perfectly centered realization \cite[Theorem~2]{GrunbaumAnti},
\cite[Theorem~3.4]{Dobbins}, but not every $4$-polytope \cite[Theorem~3.9]{Dobbins}.
Perfect centering is equivalent to the interval condition (I) of
Definition~\ref{def:intcond}, which concerns the normal fan and the face fan of $Q$ about
an interior point (Theorem~\ref{thm:cayley}). In the companion paper \cite{Thales} this realization is given a metric form: the cutout
of the apex of a pyramid over $Q$ by orthogonal projections, an orthogonally bipolar
polytope inscribed in a sphere \cite[Theorem~4.10]{Thales}; for cubes, cross-polytopes and
regular simplices it yields Hjelmslev cuboids that are not boxes. The regular
configurations of Theorem~C (Definition~\ref{def:conf}) exist in
particular for every $3$-polytope with a midsphere centred at an interior point.

\medskip
\noindent The orthogonal case is studied in the companion paper \cite{Thales}: there
the cutouts of the orthogonal projections of a point $B$ onto the faces of a polytope
are shown to be orthogonally bipolar and to tile the polytope exactly when $B$ lies in
the \emph{orthogonal kernel}, the open convex set of points whose projections onto all
faces lie in their relative interiors, i.e.\ about which the polytope is perfectly
centered \cite[Theorem~6.5]{Thales}.

\medskip
\noindent\textbf{Related work.} Abstract polytopes were introduced by McMullen and
Schulte \cite{McMullenSchulte}; the interval poset construction of
Lemma~\ref{lem:AP} is in the same spirit. Poset-theoretic aspects of face lattices go
back to Bj\"orner \cite{Bjorner} and Stanley's Eulerian posets \cite{Stanley}. In
dimension $3$, the canonical realizations used in Section~\ref{sec:geom} rest on the
Koebe--Andreev--Thurston circle-packing theorem \cite{Koebe,Andreev,Thurston},
\cite[Theorems~4.12 and~4.13]{Ziegler}.
Limitations of edge-tangency in higher dimensions are discussed by Schulte
\cite{Schulte} and Eppstein--Kuperberg--Ziegler \cite{EKZ}.

\emph{Interval polytopes and antiprisms.} The poset of intervals of the face lattice of
a polytope was considered by Lindstr\"om \cite{Lindstrom}, who asked whether it is
always the face lattice of a convex polytope. A polytope with this face lattice is an
\emph{interval polytope}, and its dual is an \emph{antiprism}: a polytope with two
disjoint facets, realizations of $Q$ and of its dual, that together contain all vertices
\cite{Dobbins}. Broadie \cite{Broadie} showed that $Q$ and $Q^\circ$ are the two bases of an
antiprism if and only if $Q$ is \emph{perfectly centered}, i.e.\ for every nonempty face
$F$ the foot of the perpendicular from the origin to $\aff F$ lies in the relative
interiors of $F$ and of its normal cone \cite{FukudaWeibel}; in the terminology of
Section~\ref{sec:convex}, $Q$ is orthogonally reciprocal. Every $3$-polytope has a
perfectly centered, namely a midscribed, realization \cite[Theorem~2]{GrunbaumAnti},
\cite[Theorem~3.4]{Dobbins}, whereas Dobbins \cite[Theorem~3.9]{Dobbins} constructed a
combinatorial $4$-polytope without any antiprism, answering Lindstr\"om's question in
the negative. Gleason and Hubard \cite{GleasonHubard} defined the antiprism of every abstract
polytope, as the poset of its intervals ordered by reverse inclusion with a new maximum,
and proved that it is again an abstract polytope \cite[Theorem~2.1]{GleasonHubard}; they
also described its facets and co-faces and showed that a polytope and its dual have
isomorphic antiprisms \cite[Proposition~2.2]{GleasonHubard}. Lawrence \cite{Lawrence}
showed that the face lattice of an antiprism is dual to the lattice of intervals of the
face lattice of its base \cite[Lemma~2.1]{Lawrence}, and that the construction
$P\mapsto\conv\bigl((P\times\{-1\})\cup(P^\circ\times\{1\})\bigr)$, i.e.\ the Cayley polytope
of $P$ and $P^\circ$, preserves perfect centering \cite[Theorem~2.2]{Lawrence}; so it can
be iterated, and starting from a segment it yields polytopes realizing all the iterated
interval lattices, whose face numbers are counted by binary partitions. In this language $W(\F(Q))$
is the interval polytope of $Q$, Theorem~\ref{thm:cayley} contains Broadie's theorem in polar
form, and the
``if'' part of Theorem~B is, dually, the theorem of Gleason and Hubard for lattices. A
bipolar polytope is dual to a prismoid, a polytope with two disjoint facets that
together contain all vertices \cite{Dobbins}; Theorem~A can be read dually as a characterization, by the
intersection property, of the prismoids that are antiprisms of their bases. What the
present paper adds is this axiomatic and lattice-theoretic description
(Theorems~A and~B) and the
decompositions of Section~\ref{sec:geom} (Theorem~C).
Combinatorially, the cutout of a vertex $v$ corresponds to the closed
star of $v$ in the barycentric subdivision of $P$ (with $B$ as the
barycentre of $P$), coarsened so that its cells correspond to the
intervals $[G,H]$ of faces containing $v$; for a simplex this is the
classical decomposition into cubes. For a simple polytope these cells
are combinatorial cubes, and the decomposition is the cubical subdivision
$\mathcal C(P)$ of Buchstaber and Panov \cite[Construction~2.10.7]{BuchstaberPanovCh2},

whose maximal cells $C_v$ are indexed by the vertices $v$ of $P$; for a
tetrahedron it is the splitting into four hexahedra used in hexahedral
mesh generation \cite{Eppstein}. For cubical complexes the analogous
\emph{cubical barycentric subdivision} has a face poset isomorphic to the
poset of closed intervals of the face poset with a least element adjoined
\cite{Savvidou}, the counterpart of the poset $\widetilde I(L)$ of
Lemma~\ref{lem:AP}. What is new here is the convexity and
the metrical realization of these cells.

\section{Preliminaries}\label{sec:prelim}

This section fixes notation and recalls the notions used throughout the paper. Most of
them are standard, and we refer to the literature for their definitions; the few
notions introduced here are elementary conventions. The new notions of the paper are
defined where they are first needed.

\subsection{Polytopes and their face lattices}\label{sec:poly}

\paragraph{Polytopes.} A \emph{polytope} is the convex hull of finitely many points
of $\R^n$. The \emph{face lattice} $\F(Q)$ is the set of faces of $Q$ including
$\emptyset$ and $Q$, ordered by inclusion; its atoms are the vertices, forming a set
$V(Q)$, and its coatoms are the facets. For a vertex $v$ the faces containing $v$ form
the face lattice of the \emph{vertex figure} $Q/v$; equivalently, they correspond to
the faces of the tangent cone $T_v(Q)=v+\cone\{u-v:u\in Q\}$. Every \emph{section}
$[F,G]=\{H:F\subseteq H\subseteq G\}$ of $\F(Q)$ is again the face lattice of a
polytope \cite[Theorem~2.7(ii)]{Ziegler}.

\paragraph{Abstract polytope.} A poset $L$ with a strictly monotone rank
function $\rk\colon L\to\{-1,\dots,n\}$, a least element $F_{-1}$ of rank
$-1$ and a greatest element $F_{n}$ of rank $n$, is an \emph{abstract
polytope of rank $n$} (in the sense of \cite{McMullenSchulte}) if
\begin{enumerate}[label=(P\arabic*)]
\item every flag (a maximal chain) contains exactly one element of each
      rank;
\item (diamond) every section $G/F=[F,G]$ with $\rk G-\rk F=2$ has exactly
      four elements;
\item (strong flag-connectedness) every section is flag-connected, i.e.\ any
      two of its flags can be joined by a sequence of flags in which
      consecutive ones differ in exactly one element.
\end{enumerate}
\label{def:AP}

\subsection{Notation and conventions}\label{sec:notation}

\paragraph{Points and vectors.} We work in $\R^n$ with the standard inner product
$\langle\cdot,\cdot\rangle$ and the norm $\|x\|=\langle x,x\rangle^{1/2}$. The origin is denoted by $O$ or $O_n$ to emphasize the dimension of the space. Whenever an
origin is fixed (this is the case for polarity and in Section~\ref{sec:char}), a point
$x\in\R^n$ is identified with its position vector $x-0$. Thus, for points $x,y$, the
expression $\langle x,y\rangle$ is the inner product of their position vectors,
$\|x\|$ is the distance of $x$ from the origin, and $x-y$ is the vector from $y$ to
$x$.

\paragraph{Hulls and dimension.} For $S\subseteq\R^n$ let $\aff S$ be the
\emph{affine hull}, the smallest affine subspace containing $S$, i.e.\ the set of all
affine combinations $\sum\lambda_ix_i$ with $x_i\in S$, $\lambda_i\in\R$,
$\sum\lambda_i=1$ (negative coefficients are allowed, in contrast to $\conv S$);
$\conv S$ is the convex hull, $\cone S$ the set of non-negative combinations of points
of $S$, and $\lin S$ the linear span. $\relint S$ and $\inte S$ denote the relative
interior (the interior within $\aff S$) and the interior. The dimension of a convex set
is $\dim S=\dim\aff S$. For a face $F$ of a polytope, $\aff F$ is the plane (line,
hyperplane, \dots) of $F$; the vertex $\{v\}$ has $\aff\{v\}=\{v\}$.

\paragraph{Perpendicularity.} A \emph{vertex vector} of a polytope $W$ is a difference
$x-y$ of two of its vertices. Two faces $F,G$ of $W$ are \emph{perpendicular}, written
$F\perp G$, if every vertex vector of $F$ is orthogonal to every vertex vector of $G$;
equivalently, the direction spaces of $\aff F$ and $\aff G$ are orthogonal.

\paragraph{Direction space.} For a face $F$ of a polytope we denote by $U_F$ the
direction space of $\aff F$, i.e.\ $U_F=\{x-y:x,y\in\aff F\}$; equivalently, $U_F$ is
the linear span of the directions of the edges of $F$ at any vertex of $F$.

\paragraph{Normal cones.} For a face $F$ of a polytope $Q\subset\R^d$ let
\[
N_Q(F)=\{c\in\R^d:\ \langle c,y\rangle=\max_{u\in Q}\langle c,u\rangle\ \text{ for all }y\in F\},
\]
the (outer) \emph{normal cone} of $F$. A vector $c$ lies in its relative interior iff
$F$ is \emph{exactly} the set of maximizers of $\langle c,\cdot\rangle$ on $Q$; we say
that $c$ \emph{exposes} $F$. The relative interiors of the cones $N_Q(F)$ partition
$\R^d$ (the normal fan); $\dim N_Q(F)=d-\dim F$.

\paragraph{Polar polytope and dual faces.} Let $Q\subset\R^d$ be a polytope with
$O\in\inte Q$. Its \emph{polar} is
\[
Q^{\circ}=\{z\in\R^d:\ \langle z,y\rangle\le1\ \text{ for all }y\in Q\}.
\]
It is again a  polytope with $O\in\inte Q^\circ$, and $(Q^\circ)^\circ=Q$ \cite[Theorem~2.11(v)]{Ziegler}.
The vertices of $Q^\circ$ correspond to the facets of $Q$: the facet
$\{y\in Q:\langle n,y\rangle=1\}$ corresponds to the vertex $n$. More generally, for a
face $F$ of $Q$ the \emph{dual face} is
\[
F^{*}=\{z\in Q^\circ:\ \langle z,y\rangle=1\ \text{ for all }y\in F\},
\]
a face of $Q^\circ$ of dimension $d-1-\dim F$ (with $\emptyset^{*}=Q^\circ$ and
$Q^{*}=\emptyset$). The map $F\mapsto F^{*}$ is an inclusion-reversing bijection
$\F(Q)\to\F(Q^\circ)$ with $(F^{*})^{*}=F$; hence
$\F(Q^\circ)\cong\F(Q)^{\mathrm{op}}$. Moreover $N_Q(F)=\cone(F^{*})$ for every
nonempty proper face $F$. Examples: the polar of the cube $[-1,1]^d$ is the
cross-polytope $\conv\{\pm e_i\}$; the polar of a regular $k$-gon centred at $O$ is a
regular $k$-gon (rotated by $\pi/k$). We say that a property of $Q$ is \emph{invariant
under polarity} if it holds for $Q$ iff it holds for $Q^\circ$ (with respect to the
same origin).

\begin{figure}[htbp]
\centering
\begin{tikzpicture}[scale=0.85,line join=round,line cap=round]
\filldraw[fill=blue!14,draw=blue!70!black,line width=0.9pt] (-1.000,-0.600)--(1.000,-0.200)--(0.100,1.000)--cycle;
\filldraw[fill=green!18,draw=green!50!black,line width=0.45pt,double=white,double distance=0.8pt] (1.500,-0.500)--(2.500,-0.500)--(2.500,0.500)--(1.500,0.500)--cycle;
\fill[black!6] (3.100,-0.100)--(3.100,-1.100)--(4.100,-1.100)--(6.100,-0.700)--(6.100,0.300)--(5.200,1.500)--(4.200,1.500)--cycle;
\draw[green!50!black,line width=0.55pt,double=white,double distance=1.1pt,line join=miter] (3.100,-0.100)--(3.100,-1.100)--(4.100,-1.100);
\draw[green!50!black,line width=0.55pt,double=white,double distance=1.1pt] (6.100,-0.700)--(6.100,0.300);
\draw[green!50!black,line width=0.55pt,double=white,double distance=1.1pt] (5.200,1.500)--(4.200,1.500);
\draw[blue!75!black,line width=2.1pt,line cap=round] (4.100,-1.100)--(6.100,-0.700);
\draw[blue!75!black,line width=2.1pt,line cap=round] (6.100,0.300)--(5.200,1.500);
\draw[blue!75!black,line width=2.1pt,line cap=round] (4.200,1.500)--(3.100,-0.100);
\node at (1.15,0.1) {$+$}; \node at (3.0,0.1) {$=$};
\draw[black,line width=0.5pt] (0.000,0.000) ++(-2.2pt,-2.2pt)--++(4.4pt,4.4pt) (0.000,0.000) ++(-2.2pt,2.2pt)--++(4.4pt,-4.4pt);
\node[font=\scriptsize,anchor=north west,inner sep=1.5pt] at (0.000,0.000) {$O$};
\draw[black,line width=0.5pt] (2.000,0.000) ++(-2.2pt,-2.2pt)--++(4.4pt,4.4pt) (2.000,0.000) ++(-2.2pt,2.2pt)--++(4.4pt,-4.4pt);
\node[font=\scriptsize,anchor=north west,inner sep=1.5pt] at (2.000,0.000) {$O$};
\draw[black,line width=0.5pt] (4.600,0.000) ++(-2.2pt,-2.2pt)--++(4.4pt,4.4pt) (4.600,0.000) ++(-2.2pt,2.2pt)--++(4.4pt,-4.4pt);
\node[font=\scriptsize,anchor=north west,inner sep=1.5pt] at (4.600,0.000) {$O$};
\node[font=\small] at (0,-1.6) {$A$}; \node[font=\small] at (2,-1.6) {$B$}; \node[font=\small] at (4.6,-1.6) {$A+B$};
\end{tikzpicture}\hspace{0.8cm}
\begin{tikzpicture}[scale=1.45,line join=round,line cap=round]
\fill[blue!8] (0.346,2.085)--(-1.125,-0.096)--(0.909,-0.146)--cycle;
\fill[green!10] (-0.617,2.170)--(-1.125,-0.096)--(0.346,2.085)--cycle;
\fill[green!10] (0.346,2.085)--(0.909,-0.146)--(0.617,2.387)--cycle;
\fill[green!18] (-0.346,2.472)--(-0.617,2.170)--(0.346,2.085)--(0.617,2.387)--cycle;
\draw[blue!70!black,line width=0.9pt] (-1.125,-0.096)--(0.909,-0.146);
\draw[blue!70!black,dashed,line width=0.5pt,opacity=0.6] (0.909,-0.146)--(0.367,0.294);
\draw[black,line width=0.9pt] (-1.125,-0.096)--(0.346,2.085);
\draw[green!50!black,line width=0.45pt,double=white,double distance=0.8pt] (-0.617,2.170)--(0.346,2.085);
\draw[black,line width=0.9pt] (0.909,-0.146)--(0.617,2.387);
\draw[black,line width=0.9pt] (-1.125,-0.096)--(-0.617,2.170);
\draw[black,line width=0.9pt] (0.909,-0.146)--(0.346,2.085);
\draw[black,dashed,line width=0.5pt,opacity=0.6] (-1.125,-0.096)--(-0.346,2.472);
\draw[blue!70!black,dashed,line width=0.5pt,opacity=0.6] (-1.125,-0.096)--(0.367,0.294);
\draw[green!50!black,line width=0.45pt,double=white,double distance=0.8pt] (0.346,2.085)--(0.617,2.387);
\draw[black,dashed,line width=0.5pt,opacity=0.6] (0.367,0.294)--(-0.346,2.472);
\draw[green!50!black,line width=0.45pt,double=white,double distance=0.8pt] (0.617,2.387)--(-0.346,2.472);
\draw[green!50!black,line width=0.45pt,double=white,double distance=0.8pt] (-0.617,2.170)--(-0.346,2.472);
\draw[black,dashed,line width=0.5pt,opacity=0.6] (0.367,0.294)--(0.617,2.387);
\draw[green!50!black,line width=0.45pt,double=white,double distance=0.8pt,densely dashed] (-0.736,1.188)--(-0.871,1.037);
\draw[green!50!black,line width=0.55pt,double=white,double distance=1.1pt] (-0.871,1.037)--(-0.390,0.994);
\draw[green!50!black,line width=0.55pt,double=white,double distance=1.1pt] (0.627,0.970)--(0.763,1.121);
\draw[green!50!black,line width=0.45pt,double=white,double distance=0.8pt,densely dashed] (0.492,1.341)--(0.011,1.383);
\draw[blue!75!black,line width=2.1pt,line cap=round] (-0.390,0.994)--(0.627,0.970);
\draw[blue!75!black,line width=0.8pt,densely dashed,line cap=round] (0.763,1.121)--(0.492,1.341);
\draw[blue!75!black,line width=0.8pt,densely dashed,line cap=round] (0.011,1.383)--(-0.736,1.188);
\fill[black] (-1.125,-0.096) circle (1.1pt);
\fill[black] (0.909,-0.146) circle (1.1pt);
\fill[black] (0.367,0.294) circle (1.1pt);
\fill[black] (-0.617,2.170) circle (1.1pt);
\fill[black] (0.346,2.085) circle (1.1pt);
\fill[black] (0.617,2.387) circle (1.1pt);
\fill[black] (-0.346,2.472) circle (1.1pt);
\node[blue!70!black,font=\small,anchor=north] at (-0.051,-0.289) {$A\times\{0\}$};
\node[green!45!black,font=\small,anchor=south] at (0.162,2.460) {$B\times\{1\}$};
\node[font=\small,anchor=west] at (1.313,1.148) {$\tfrac12(A+B)\times\{\tfrac12\}$};
\draw[->,black!60] (1.313,1.148)--(0.813,1.148);
\end{tikzpicture}
\caption{The Cayley trick for a triangle $A$ and a square $B$. Left: the Minkowski sum
$A+B$; its edges are the edges of $A$ (blue, single lines) and of $B$ (green,
double lines), translated. Right:
the Cayley polytope $\mathrm{Cay}(A,B)$ with $A$ at height $0$ and $B$ at height $1$
(heights stretched). Its middle slice is $\tfrac12(A+B)$. Each of its $7$ mixed
facets is $\conv(F\times\{0\}\cup G\times\{1\})$ for faces $F$ of $A$ and $G$ of $B$
with a common outer normal --- an edge of $A$ with a vertex of $B$ (single line) or a
vertex of $A$ with an edge of $B$ (double line) --- and meets the slice in the corresponding edge
$\tfrac12(F+G)$ of $\tfrac12(A+B)$. The origin $O$ is marked in each polygon; it
affects $A+B$ and $\mathrm{Cay}(A,B)$ only by a translation.}
\label{fig:cayley}
\end{figure}
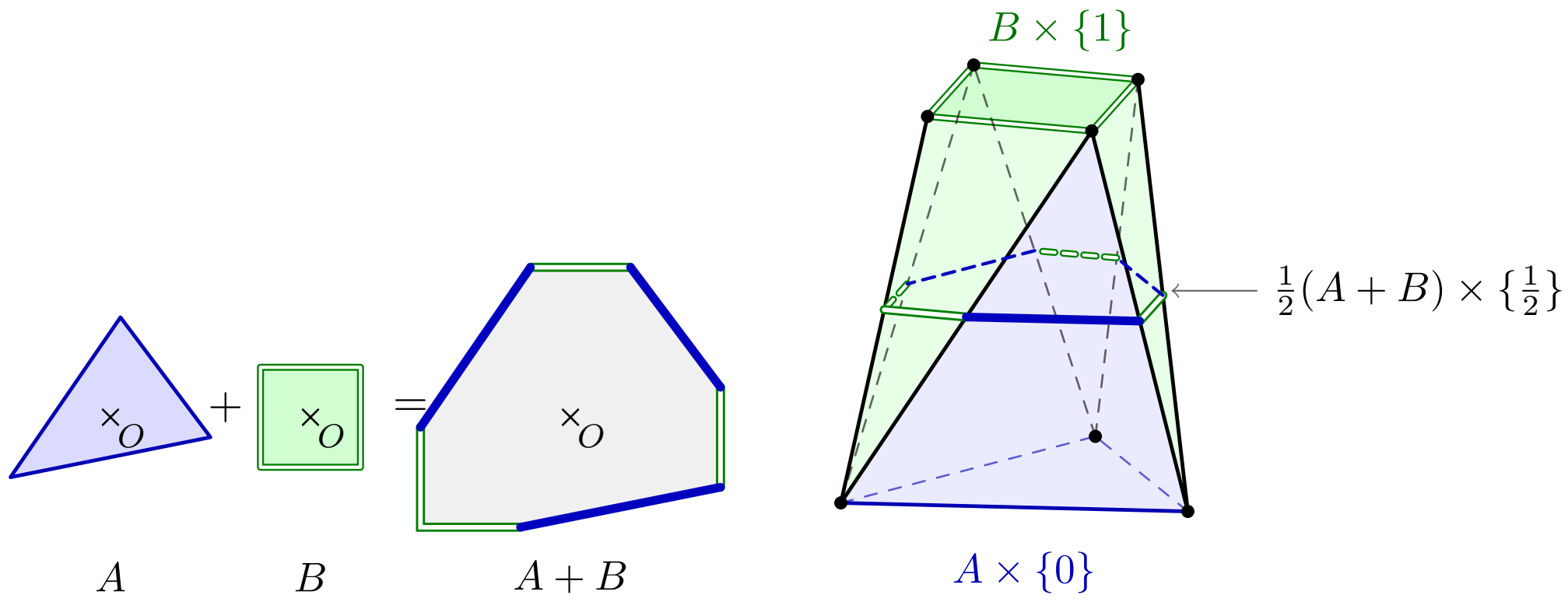

\paragraph{Minkowski sums and Cayley polytopes.} The \emph{Minkowski sum} of
polytopes $A,B\subset\R^d$ is $A+B=\{a+b:a\in A,\ b\in B\}$. For $c\in\R^d$ let
$A^c$ denote the face of $A$ on which $\langle c,\cdot\rangle$ attains its maximum
$h_A(c)=\max_{a\in A}\langle c,a\rangle$; the function $h_A$ is the \emph{support
function} of $A$. Then $(A+B)^c=A^c+B^c$, so the proper faces of $A+B$ correspond to the pairs
$(F,G)$ of faces of $A$ and $B$ with a common normal in the relative interiors,
$\relint N_A(F)\cap\relint N_B(G)\ne\emptyset$; these are the cells of the common
refinement of the two normal fans. The \emph{Cayley polytope}
\[
\mathrm{Cay}(A,B)=\conv\bigl(A\times\{0\}\cup B\times\{1\}\bigr)\subset\R^{d+1}
\]
encodes the same information in a single polytope (the \emph{Cayley trick}
\cite{HRS}). A linear functional $(c,\lambda)$ on $\R^d\times\R$ attains on it the
maximum $\max\bigl(h_A(c),\,h_B(c)+\lambda\bigr)$. Hence the nonempty faces of
$\mathrm{Cay}(A,B)$ are the polytope itself, the faces of the two bases
$A\times\{0\}$ and $B\times\{1\}$, and, for $c\ne O_d$ and $\lambda=h_A(c)-h_B(c)$, the
\emph{mixed faces} $\conv\bigl(A^c\times\{0\}\cup B^c\times\{1\}\bigr)$. The mixed
faces therefore correspond bijectively to the proper faces of $A+B$; indeed the slice
of $\mathrm{Cay}(A,B)$ at height $t\in(0,1)$ is $(1-t)A+tB$
(Figure~\ref{fig:cayley}). Every vertex of $A$ or $B$ gives a vertex of
$\mathrm{Cay}(A,B)$, and since normal fans do not change under positive scaling, the
face lattice of $\mathrm{Cay}(\lambda A,\mu B)$ does not depend on $\lambda,\mu>0$.

\paragraph{Realizations.} Two polytopes are \emph{combinatorially
equivalent} if their face lattices are isomorphic; the class of all polytopes
combinatorially equivalent to $Q$ is the \emph{combinatorial type} of $Q$, and each of
its members is a \emph{realization} of the type.

\paragraph{Projective transformations.} A \emph{projective transformation} of $\R^n$ is a
map $\varphi(x)=(Ax+b)/(\langle c,x\rangle+\delta)$, where the $(n+1)\times(n+1)$ matrix
with blocks $A,b,c^{\mathsf T},\delta$ is invertible; it is defined outside the hyperplane
$H_\varphi=\{x:\langle c,x\rangle+\delta=0\}$, which it sends to infinity, and affine maps are
the case $c=0$. It is \emph{admissible} for a polytope $P$ if $P\cap H_\varphi=\emptyset$.
Then $\varphi(P)$ is a polytope combinatorially equivalent to $P$: $\varphi$ maps the faces of
$P$ onto the faces of $\varphi(P)$, relative interiors onto relative interiors, collinear
and coplanar points to collinear and coplanar points, and $\varphi(\conv S)=\conv\varphi(S)$
for every $S\subseteq P$ \cite[Section~2.6]{Ziegler}. Two polytopes are \emph{projectively
equivalent} if one is the image of the other under an admissible projective
transformation; this is much stronger than combinatorial equivalence.

\paragraph{Tilings.} A finite family of $n$-polytopes $C_1,\dots,C_k\subseteq\R^n$
\emph{tiles} an $n$-polytope $P$, or forms a \emph{tiling} of $P$, if
$C_1\cup\dots\cup C_k=P$ and the interiors of the $C_i$ are pairwise disjoint. The
tiling is \emph{face-to-face} if moreover any two of the $C_i$ meet in a common face of
both (possibly empty).

\paragraph{Voronoi diagrams, power diagrams and Delaunay subdivisions.} Let
$S\subset\R^n$ be a finite set of points, the \emph{sites}. The \emph{Voronoi cell} of
$s\in S$ is $\{p:\ \|p-s\|\le\|p-t\|\text{ for all }t\in S\}$; the cells are convex
polyhedra that cover $\R^n$ and have pairwise disjoint interiors, and they form the
\emph{Voronoi diagram} of $S$. Two cells are separated by the perpendicular bisector of
their sites. Dually, the cells of the \emph{Delaunay subdivision} of $S$ are the convex
hulls of the sets of sites whose Voronoi cells share a vertex of the diagram; for sites
in general position it is a triangulation, the \emph{Delaunay triangulation},
characterized by the property that the circumsphere of each simplex has no site in its
interior. More generally, let each site carry a weight $w_s\in\R$, and define the
\emph{power} of a point $p$ with respect to $s$ by
\[
\operatorname{pow}_s(p)=\|p-s\|^2-w_s ;
\]
for $w_s=r^2>0$ this is the power of $p$ with respect to the sphere with centre $s$ and
radius $r$, which for $p$ outside the sphere is the squared length of a tangent
segment from $p$. Replacing the distance by the power gives the \emph{power diagram}
(or Laguerre diagram); two of its cells are separated by the \emph{radical hyperplane}
of the two spheres, which is perpendicular to the segment joining the sites but in
general does not pass through its midpoint, and some cells may be empty. The Voronoi
diagram is the power diagram with equal weights. Since
$\operatorname{pow}_s(p)=\|p\|^2-\bigl(2\langle p,s\rangle-\|s\|^2+w_s\bigr)$, the
power cell of $s$ is the set of points where the affine function
$2\langle p,s\rangle-\|s\|^2+w_s$ is largest among these functions; conversely, every
family of affine functions $\langle p,a_s\rangle+b_s$ arises in this way, with the
sites $a_s/2$ and the weights $b_s+\|a_s\|^2/4$. So power diagrams are exactly the
subdivisions of space into the regions where one of finitely many affine functions is
largest \cite[Observation~4 and Theorem~4]{Aurenhammer87}, \cite{Aurenhammer91}. Power
diagrams of spheres through a common
point appear in \cite[Corollary~5.3 and Remark~6.4]{Thales}.

\subsection{Discrete polytopes}\label{sec:discrete}

\begin{definition}[Discrete polytope]\label{def:discrete}
A \emph{discrete polytope} is a pair $(P,\M)$ where $P$ is a finite set
(\emph{vertices}) and $\M$ is a family of nonempty subsets of $P$ (\emph{facets}; by Lemma~\ref{lem:sep}(b) they are exactly the
coatoms of $\F(P,\M)$, so the name agrees with the usual one)
such that
\begin{enumerate}[label=(D\arabic*)]
\item $\bigcup\M=P$ and no member of $\M$ contains another;
\item the family $\F(P,\M)=\{\bigcap\mathcal N:\mathcal N\subseteq\M\}\cup\{\emptyset\}$
      (with $\bigcap\emptyset=P$), ordered by inclusion, is the face poset of an
      abstract polytope \cite{McMullenSchulte}, and every singleton $\{x\}$, $x\in P$,
      belongs to $\F(P,\M)$.
\end{enumerate}
\end{definition}

\begin{remark}\label{rem:D2}
Let $W$ be a convex $n$-polytope. Then $\F(W)$ is an abstract polytope of rank $n$ in the
sense of \S\ref{sec:poly}; that is, $(V(W),\{\text{vertex sets of facets}\})$
satisfies (D2). We record the two parts that are immediate and refer for the third.

\emph{(P1) Grading.} Put $\rk F=\dim F$; thus $\rk\emptyset=-1$, vertices have rank $0$
and $\rk W=n$. Then
$G$ covers $F$ iff $\dim G=\dim F+1$, and every chain can be refined to a complete flag
$\emptyset\subset F_0\subset\dots\subset F_n=W$; hence every flag has exactly one element
of each rank \cite[Theorem~2.7(i)]{Ziegler} (where the rank is $\dim F+1$).

\emph{(P2) Diamond.} For faces $F\subseteq G$ the section $[F,G]$ is isomorphic to
$\F(G/F)$, where $G/F$ is a polytope of dimension $\dim G-\dim F-1$
\cite[Theorem~2.7(ii)]{Ziegler}. If $\dim G-\dim F=2$, then $G/F$ is a segment and $[F,G]$ has four
elements; this is the diamond property \cite[Theorem~2.7(iii)]{Ziegler},
\cite[\S3.1]{Grunbaum}.

\emph{(P3) Strong flag-connectedness.} This is the only non-elementary part. Its content
is that flags of a polytope are connected by \emph{flips} --- replacing one member of a
flag by the other member of the corresponding diamond --- and it is a classical property of
face lattices of convex polytopes; see \cite[\S 2A]{McMullenSchulte}, where the face
lattice of a convex polytope is treated as the motivating example of an abstract
polytope.

Nothing in Sections~\ref{sec:bipolar}--\ref{sec:char} depends on (P3) beyond the
inheritance ``sections of abstract polytopes are abstract polytopes''; we never use
flag-connectedness in a geometric argument. Finally, every vertex is a face of dimension
$0$, so all singletons belong to $\F(W)$, which completes (D2).
\end{remark}

Elements of $\F(P,\M)$ are \emph{faces}; $P$ itself is the top face. If $Q$ is a
convex polytope, $(V(Q),\{\text{vertex sets of facets}\})$ is a discrete polytope:
(D1) is immediate, and (D2) is the standard statement that the face lattice of a
convex polytope is an abstract polytope of rank $n=\dim Q$, with the rank function
$\rk F={}$$\dim F$ (so $\rk\emptyset=-1$); see \cite[Theorem~2.7(i)--(iii)]{Ziegler} for the
grading, the intervals and the diamond property (also \cite[\S3.1]{Grunbaum}), and
\cite[{\S}2A]{McMullenSchulte}, and
Remark~\ref{rem:D2} for the short argument. The
poset $\F(P,\M)$ is a lattice (the join of two faces is the intersection of all
facets containing both) which is atomistic and coatomistic. 

\begin{definition}[abstractly polytopal lattice]\label{def:apl}
A finite lattice $L$ is \emph{abstractly polytopal} if
\begin{enumerate}[label=(L\arabic*)]
\item $L$ is atomistic, i.e.\ every element is the join of the atoms below it;
\item $L$ is coatomistic, i.e.\ every element is the meet of the coatoms above it;
\item $L$, as a poset, is an abstract polytope in the sense of
      \S\ref{sec:poly}, i.e.\ it satisfies (P1)--(P3).
\end{enumerate}
The \emph{rank} of $L$ is its rank as an abstract polytope.
\end{definition}

\begin{figure}[htbp]
\centering
\begin{tikzpicture}[scale=1.25,line join=round]
\fill[blue!10] (0.000,0.000)--(1.000,0.000)--(0.500,0.866)--cycle;
\fill[orange!14] (1.000,0.000)--(1.500,0.866)--(0.500,0.866)--cycle;
\fill[blue!18] (1.000,0.000)--(2.000,0.000)--(1.500,0.866)--cycle;
\fill[orange!24] (2.000,0.000)--(2.500,0.866)--(1.500,0.866)--cycle;
\fill[blue!10] (2.000,0.000)--(3.000,0.000)--(2.500,0.866)--cycle;
\fill[orange!14] (3.000,0.000)--(3.500,0.866)--(2.500,0.866)--cycle;
\fill[blue!18] (3.000,0.000)--(4.000,0.000)--(3.500,0.866)--cycle;
\fill[orange!24] (4.000,0.000)--(4.500,0.866)--(3.500,0.866)--cycle;
\fill[blue!10] (4.000,0.000)--(5.000,0.000)--(4.500,0.866)--cycle;
\fill[orange!14] (5.000,0.000)--(5.500,0.866)--(4.500,0.866)--cycle;
\fill[blue!18] (5.000,0.000)--(6.000,0.000)--(5.500,0.866)--cycle;
\fill[orange!24] (6.000,0.000)--(6.500,0.866)--(5.500,0.866)--cycle;
\fill[blue!10] (6.000,0.000)--(7.000,0.000)--(6.500,0.866)--cycle;
\fill[orange!14] (7.000,0.000)--(7.500,0.866)--(6.500,0.866)--cycle;
\draw[line width=0.5pt] (0.000,0.000)--(0.500,0.866);
\draw[line width=0.5pt] (1.000,0.000)--(1.500,0.866);
\draw[line width=0.5pt] (2.000,0.000)--(2.500,0.866);
\draw[line width=0.5pt] (3.000,0.000)--(3.500,0.866);
\draw[line width=0.5pt] (4.000,0.000)--(4.500,0.866);
\draw[line width=0.5pt] (5.000,0.000)--(5.500,0.866);
\draw[line width=0.5pt] (6.000,0.000)--(6.500,0.866);
\draw[line width=0.5pt] (7.000,0.000)--(7.500,0.866);
\draw[line width=0.5pt] (1.000,0.000)--(0.500,0.866);
\draw[line width=0.5pt] (2.000,0.000)--(1.500,0.866);
\draw[line width=0.5pt] (3.000,0.000)--(2.500,0.866);
\draw[line width=0.5pt] (4.000,0.000)--(3.500,0.866);
\draw[line width=0.5pt] (5.000,0.000)--(4.500,0.866);
\draw[line width=0.5pt] (6.000,0.000)--(5.500,0.866);
\draw[line width=0.5pt] (7.000,0.000)--(6.500,0.866);
\draw[line width=0.9pt,red!70!black] (0.000,0.000)--(7.000,0.000);
\draw[line width=0.9pt,red!70!black] (0.500,0.866)--(7.500,0.866);
\draw[line width=0.9pt,green!50!black] (0.000,0.000)--(0.500,0.866) (7.000,0.000)--(7.500,0.866);
\fill[white] (0.000,0.000) circle (5.2pt); \draw[line width=0.4pt] (0.000,0.000) circle (5.2pt); \node[font=\scriptsize] at (0.000,0.000) {0};
\fill[white] (0.500,0.866) circle (5.2pt); \draw[line width=0.4pt] (0.500,0.866) circle (5.2pt); \node[font=\scriptsize] at (0.500,0.866) {3};
\fill[white] (1.000,0.000) circle (5.2pt); \draw[line width=0.4pt] (1.000,0.000) circle (5.2pt); \node[font=\scriptsize] at (1.000,0.000) {1};
\fill[white] (1.500,0.866) circle (5.2pt); \draw[line width=0.4pt] (1.500,0.866) circle (5.2pt); \node[font=\scriptsize] at (1.500,0.866) {4};
\fill[white] (2.000,0.000) circle (5.2pt); \draw[line width=0.4pt] (2.000,0.000) circle (5.2pt); \node[font=\scriptsize] at (2.000,0.000) {2};
\fill[white] (2.500,0.866) circle (5.2pt); \draw[line width=0.4pt] (2.500,0.866) circle (5.2pt); \node[font=\scriptsize] at (2.500,0.866) {5};
\fill[white] (3.000,0.000) circle (5.2pt); \draw[line width=0.4pt] (3.000,0.000) circle (5.2pt); \node[font=\scriptsize] at (3.000,0.000) {3};
\fill[white] (3.500,0.866) circle (5.2pt); \draw[line width=0.4pt] (3.500,0.866) circle (5.2pt); \node[font=\scriptsize] at (3.500,0.866) {6};
\fill[white] (4.000,0.000) circle (5.2pt); \draw[line width=0.4pt] (4.000,0.000) circle (5.2pt); \node[font=\scriptsize] at (4.000,0.000) {4};
\fill[white] (4.500,0.866) circle (5.2pt); \draw[line width=0.4pt] (4.500,0.866) circle (5.2pt); \node[font=\scriptsize] at (4.500,0.866) {0};
\fill[white] (5.000,0.000) circle (5.2pt); \draw[line width=0.4pt] (5.000,0.000) circle (5.2pt); \node[font=\scriptsize] at (5.000,0.000) {5};
\fill[white] (5.500,0.866) circle (5.2pt); \draw[line width=0.4pt] (5.500,0.866) circle (5.2pt); \node[font=\scriptsize] at (5.500,0.866) {1};
\fill[white] (6.000,0.000) circle (5.2pt); \draw[line width=0.4pt] (6.000,0.000) circle (5.2pt); \node[font=\scriptsize] at (6.000,0.000) {6};
\fill[white] (6.500,0.866) circle (5.2pt); \draw[line width=0.4pt] (6.500,0.866) circle (5.2pt); \node[font=\scriptsize] at (6.500,0.866) {2};
\fill[white] (7.000,0.000) circle (5.2pt); \draw[line width=0.4pt] (7.000,0.000) circle (5.2pt); \node[font=\scriptsize] at (7.000,0.000) {0};
\fill[white] (7.500,0.866) circle (5.2pt); \draw[line width=0.4pt] (7.500,0.866) circle (5.2pt); \node[font=\scriptsize] at (7.500,0.866) {3};
\end{tikzpicture}
\caption{The $7$-vertex triangulation of the torus \cite{Csaszar}: $14$ triangles
$\{i,i+1,i+3\}$ and $\{i,i+2,i+3\}$, $i\in\mathbb{Z}/7$. The strip is glued along its two
slanted ends and along its top and bottom, where the vertices with the same
labels are identified; the bottom row is identified with the top row shifted by three
steps. Every two of the $7$ vertices span an edge. The face lattice is abstractly
polytopal of rank~$3$ but is not the face lattice of a convex polytope
(Remark~\ref{rem:apl}).}
\label{fig:torus}
\end{figure}
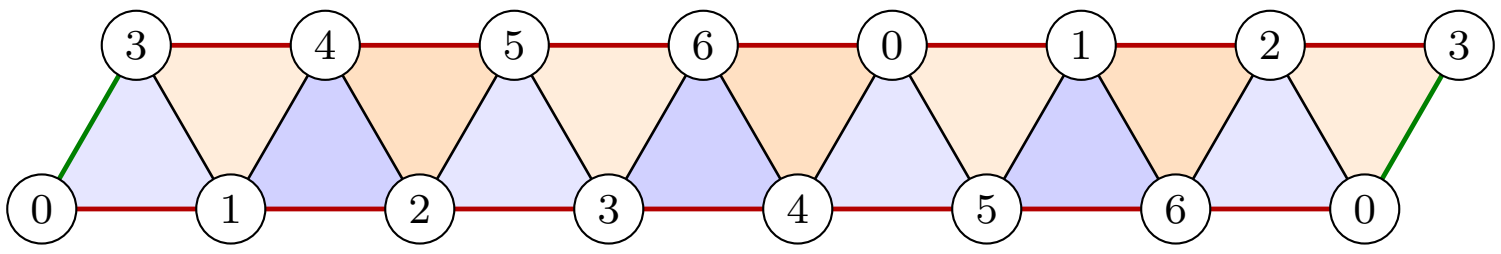

\begin{remark}\label{rem:apl}
(a) The face lattice $\F(Q)$ of a convex $n$-polytope $Q$ is abstractly
polytopal of rank $n$: (L3) is Remark~\ref{rem:D2}, (L1) holds because
every face is the convex hull of its vertices, and (L2) because every
proper face is the intersection of the facets containing it. Likewise
$\F(P,\M)$ is abstractly polytopal for every discrete polytope $(P,\M)$,
by (D2) and the preceding paragraph.

(b) We avoid the shorter term \emph{polytopal lattice}, which in the
literature usually means the face lattice of a \emph{convex} polytope;
which lattices are of this kind is the Steinitz problem
\cite{Steinitz}, \cite[Lecture~4]{Ziegler}. Every such lattice is abstractly polytopal, but
not conversely: the face lattices of non-polytopal simplicial spheres,
such as the Barnette sphere \cite{Barnette}, are abstractly polytopal.
Non-spherical examples occur already in rank~$3$: the face lattice of the
$7$-vertex triangulation of the torus (Figure~\ref{fig:torus}) satisfies
(L1)--(L3), as is checked directly, but it is not the face lattice of a convex
$3$-polytope, since the boundary of a convex $3$-polytope is a $2$-sphere. The
torus does have a non-convex realization in $\R^3$ with flat triangles, the Cs\'asz\'ar
polyhedron \cite{Csaszar}; so the obstruction is convexity, not realizability as such.
For such a lattice $L$, the discrete bipolar polytope $W(L)$ of Section~\ref{sec:char}
exists and its vertex lattice is $L$, but it has no convex realization, since the vertex
figure at a pole would be a convex polytope with face lattice $L$.
\end{remark}

\begin{lemma}\label{lem:sep}
Let $(P,\M)$ be a discrete polytope.
\begin{enumerate}[label=(\alph*)]
\item Every nonempty face $F\ne P$ is the intersection of the
      facets containing it, and for distinct $x,y\in P$ there is $M\in\M$ with $x\in M$,
      $y\notin M$.
\item If $|P|\ge2$, the members of $\M$ are exactly the maximal proper
      faces, i.e.\ the coatoms of $\F(P,\M)$; if $\F(P,\M)$ has rank $n$,
      they are its faces of rank $n-1$.
\end{enumerate}
\end{lemma}
\begin{proof}
(a) If $F=\bigcap\mathcal N$ with $\mathcal N\subseteq\M$, then $F$ is contained in every
member of $\mathcal N$ and so $F=\bigcap\{M\in\M:F\subseteq M\}$. Apply this to the
face $\{x\}$.

(b) A proper face $F\ne P$ is $\bigcap\mathcal N$ with $\mathcal N\ne\emptyset$,
so it is contained in a member of $\M$. No $M\in\M$ equals $P$: otherwise
$\M=\{P\}$ by (D1), so $\F(P,\M)=\{\emptyset,P\}$ contains no singleton,
contradicting (D2) as $|P|\ge2$. If $M\in\M$ were contained in a proper
face $F\supsetneq M$, then $F\subseteq M'$ for some $M'\in\M$ and
$M\subsetneq M'$, contradicting (D1). Hence every member of $\M$ is a
maximal proper face, and every maximal proper face, lying in some
$M\in\M$, is one of them. By (P1) every coatom lies in a flag directly
below the top element of rank $n$, so it has rank $n-1$.
\end{proof}

\section{Bipolar polytopes}\label{sec:bipolar}

\subsection{Discrete and geometric approach}\label{sec:bipdef}

\paragraph{Poles.} Let $W$ be a convex polytope with two distinct vertices
$\sigma\ne\nu$, the \emph{poles} (\emph{south} $\sigma$, \emph{north} $\nu$). A face
is a \emph{south} face if it contains $\sigma$ and a \emph{north} face if it contains
$\nu$; for a vertex $x$ we write $[\sigma,x]$ and $[x,\nu]$ for the corresponding
faces.

\begin{definition}[Bipolar polytope]\label{def:bipolar}
A \emph{bipolar polytope} is a discrete polytope $(P,\M)$ with two distinct vertices
$\sigma,\nu$ (the \emph{south} and the \emph{north pole}) such that
\begin{enumerate}[label=(B\arabic*)]
\item every $M\in\M$ contains exactly one of $\sigma,\nu$;
\item (\emph{intersection property}) with
      $\F_\sigma=\{F\in\F(P,\M):\sigma\in F\}$ and
      $\F_\nu=\{F\in\F(P,\M):\nu\in F\}$ (south and north faces; both contain $P$):
      for every $H\in\F_\sigma$ the family
      $\mathcal N(H)=\{G\in\F_\nu:G\cap H\neq\emptyset\}$ is closed under
      intersection, and for every $G\in\F_\nu$ the family
      $\mathcal S(G)=\{H\in\F_\sigma:H\cap G\neq\emptyset\}$ is closed under
      intersection.
\end{enumerate}
\end{definition}

We write $\M_\sigma,\M_\nu$ for the facets containing $\sigma$, resp.\ $\nu$. For
$x\in P$ put
\[
\hup{x}=\bigcap\{N\in\M_\nu:x\in N\},\qquad
\hdn{x}=\bigcap\{D\in\M_\sigma:x\in D\},
\]
the smallest north (south) face containing $x$.

\begin{definition}[Geometric bipolar polytope]\label{def:geombip}
A convex polytope $W\subset\R^n$ with two vertices $\sigma,\nu$ is a \emph{geometric
bipolar polytope} if $(V(W),\{\text{vertex sets of facets}\})$ with the poles
$\sigma,\nu$ satisfies (B1) and (B2). Axioms (D1) and (D2) hold automatically, by the paragraph following Definition~\ref{def:discrete}.
\end{definition}

A convex polytope $W$ with poles $\sigma,\nu$ is \emph{orthogonally bipolar} if it
is geometric bipolar and $[\sigma,a]\perp[a,\nu]$ for every atom $a$ of its vertex
lattice; an orthogonally bipolar polytope whose vertex lattice is the Boolean lattice $B_n$ is a \emph{Hjelmslev
$n$-cuboid}. These notions belong to the companion paper \cite{Thales}, where $W(\F(Q))$ is realized
as an orthogonally bipolar polytope \cite[Theorem~4.10]{Thales}; here we use them only in
Remark~\ref{rem:canon-orth}.

\subsection{Properties of discrete bipolar polytopes}

\begin{lemma}\label{lem:hull}
\begin{enumerate}[label=(\alph*)]
\item If $F\in\F_\nu$, then $F=\bigcap\{N\in\M_\nu:F\subseteq N\}$; dually for $\F_\sigma$.
\item If $x\in F\in\F_\nu$, then $\hup{x}\subseteq F$; dually for $\F_\sigma$.
\item If $y\in\hup{x}$ and $y\neq x$, there is $D\in\M_\sigma$ with $x\in D$, $y\notin D$.
      Dually, if $y\in\hdn{x}$, $y\ne x$, there is $N\in\M_\nu$ with $x\in N$, $y\notin N$.
\end{enumerate}
\end{lemma}
\begin{proof}
(a) By Lemma~\ref{lem:sep}, $F$ is the intersection of the facets containing it;
they contain $\nu$, hence lie in $\M_\nu$ by (B1). (b) By (a) $F$ is an intersection of
members of $\M_\nu$, all containing $x$. (c) Lemma~\ref{lem:sep} gives $M\in\M$ with
$x\in M\not\ni y$. If $M\in\M_\nu$, then $y\in\hup{x}\subseteq M$, a contradiction; so
$M\in\M_\sigma$.
\end{proof}

\begin{lemma}\label{lem:least}
For $G\in\F_\nu$ the family $\mathcal S(G)$ has a least element $k(G)$, and for
$H'\in\F_\sigma$ we have $H'\cap G\ne\emptyset\iff k(G)\subseteq H'$. Dually,
$\mathcal N(H)$ has a least element $l(H)$ for $H\in\F_\sigma$.
\end{lemma}
\begin{proof}
$\mathcal S(G)$ is finite, contains $P$ and is closed under intersection by (B2); its
intersection $k(G)$ lies in it. The equivalence is clear.
\end{proof}

\begin{definition}
Given $G\in\F_\nu$ and $H\in\F_\sigma$, the least elements of $\mathcal S(G)$ and
$\mathcal N(H)$ are denoted by $k(G)$ and $l(H)$.
\end{definition}

\begin{lemma}[singleton lemma]\label{lem:single}
For $G\in\F_\nu$ the set $G\cap k(G)$ consists of exactly one vertex $i_G$. Dually
$H\cap l(H)=\{j_H\}$ for $H\in\F_\sigma$.
\end{lemma}
\begin{proof}
Put $H=k(G)$; $S=G\cap H\ne\emptyset$ by Lemma~\ref{lem:least}. Suppose $x\neq y$ in
$S$.

\emph{Step 1: no $D\in\M_\sigma$ contains exactly one of $x,y$.} If $x\in D\not\ni y$,
then $D\cap G\ni x$, so $D\supseteq H\ni y$ by definition of $k$, a contradiction; the
same with $x,y$ exchanged.

\emph{Step 2.} By Lemma~\ref{lem:sep} there are $M\ni x$, $M\not\ni y$ and $M'\ni y$,
$M'\not\ni x$ in $\M$; by Step~1 both lie in $\M_\nu$. Hence $y\notin\hup{x}$ and
$x\notin\hup{y}$. Both $\hup{x}$ and $\hup{y}$ meet $H$ (in $x$, resp.\ $y$), so by
(B2) $\hup{x}\cap\hup{y}$ meets $H$: pick $z\in\hup{x}\cap\hup{y}\cap H$.

\emph{Step 3.} If $z\ne x$, Lemma~\ref{lem:hull}(c) gives $D\in\M_\sigma$ with $x\in
D\not\ni z$. But $D\cap G\ni x$ implies $D\supseteq H\ni z$, a contradiction. So
$z=x$; symmetrically $z=y$, contradicting $x\ne y$.
\end{proof}

\begin{lemma}\label{lem:gen}
Let $G\in\F_\nu$. Then $G=\hup{i_G}$, where $i_G$ is the vertex of
Lemma~\ref{lem:single}. Dually $H=\hdn{j_H}$ for $H\in\F_\sigma$.
\end{lemma}
\begin{proof}
Let $i=i_G$ and $H=k(G)$. By Lemma~\ref{lem:hull}(b), $\hup{i}\subseteq G$. Conversely
let $y\in G$; we show $y\in\hup{i}$, i.e.\ $y\in N$ for every $N\in\M_\nu$ with $i\in
N$. As $y\in G$ we have $y\in y^{\uparrow}$ and, by Lemma~\ref{lem:hull}(b),
$y^{\uparrow}\subseteq G$. Let $D=k(y^{\uparrow})$; every south face meeting
$y^{\uparrow}$ meets $G$, hence $D\supseteq H\ni i$. Now $N\cap D\ni i$ and
$y^{\uparrow}\cap N\ne\emptyset$, so by (B2) there is $z\in N\cap y^{\uparrow}\cap D$.
If $z\ne y$, Lemma~\ref{lem:hull}(c) gives $E\in\M_\sigma$ with $y\in E\not\ni z$; but
$E$ meets $y^{\uparrow}$, so $E\supseteq D\ni z$, a contradiction. Thus $z=y$, and
$y\in N$.
\end{proof}

\begin{proposition}\label{prop:xup}
For every $x\in P$: $i_{\hup{x}}=x$ and $j_{\hdn{x}}=x$.
\end{proposition}
\begin{proof}
Let $G=\hup{x}$ and $i=i_G$. By Lemma~\ref{lem:gen}, $G=\hup{i}$; so $i\in\hup{x}$. If
$i\ne x$, Lemma~\ref{lem:hull}(c) gives $D\in\M_\sigma$ with $x\in D\not\ni i$. But
$D\cap G\ni x$ implies $D\supseteq k(G)\ni i$. Hence $i=x$.
\end{proof}

\begin{proposition}\label{prop:kl}
The maps $k:\F_\nu\to\F_\sigma$ and $l:\F_\sigma\to\F_\nu$ are mutually inverse and
order reversing. Hence $\F_\nu\cong\F_\sigma^{\mathrm{op}}$.
\end{proposition}
\begin{proof}
Both maps reverse inclusion. Let $H=k(G)$, $K=l(H)$. As $G$ meets $H$, $K\subseteq G$.
Also $K\cap H\ne\emptyset$ and $K\cap H\subseteq G\cap H=\{i_G\}$, so $i_G\in K$ and
$K\supseteq\hup{i_G}=G$ by Lemmas~\ref{lem:hull}(b) and~\ref{lem:gen}. So $l(k(G))=G$;
dually $k(l(H))=H$.
\end{proof}

\begin{lemma}\label{lem:sym}
$y\in\hup{x}\iff x\in\hdn{y}$.
\end{lemma}
\begin{proof}
($\Rightarrow$) Let $G=\hup{x}\ni y$. The south face $\hdn{y}$ contains $y$, so it
meets $G$ and contains $k(G)\ni i_G=x$ (Proposition~\ref{prop:xup}). ($\Leftarrow$)
Dual.
\end{proof}

\subsection{Structural theorems}

\begin{theorem}[structure of bipolar polytopes; Theorem A]\label{thm:lattice}
Let $(P,\M)$ be bipolar. Define $x\le y:\iff y\in\hup{x}$ $(\iff x\in\hdn{y})$.
\begin{enumerate}[label=(\alph*)]
\item $\le$ is a partial order on $P$ with least element $\sigma$ and greatest element
      $\nu$; $\hup{x}=\{z:z\ge x\}$ and $\hdn{y}=\{z:z\le y\}$.
\item $(P,\le)$ is a lattice.
\item The nonempty faces of $(P,\M)$ are exactly the intervals $[x,y]$, $x\le y$.
\item $\M_\nu=\{\hup{a}:a\text{ an atom of }(P,\le)\}$ and
      $\M_\sigma=\{\hdn{c}:c\text{ a coatom}\}$.
\item $(P,\le)$ is atomistic and coatomistic.
\end{enumerate}
\end{theorem}
\begin{proof}
(a) Reflexivity: $x\in\hup{x}$. Transitivity: if $y\in\hup{x}$ then
$\hup{x}\supseteq\hup{y}$ by Lemma~\ref{lem:hull}(b), so $z\in\hup{y}$ implies
$z\in\hup{x}$. Antisymmetry: $x\le y\le x$ gives $\hup{x}=\hup{y}$, hence
$x=i_{\hup{x}}=i_{\hup{y}}=y$ by Proposition~\ref{prop:xup}. No north facet contains
$\sigma$, so $\hup{\sigma}=P$; every $\hup{x}$ is an intersection of north facets, so
it contains $\nu$. The descriptions follow from the definition and
Lemma~\ref{lem:sym}.

(b) $\hup{x}\cap\hup{y}$ is a north face, so by Lemma~\ref{lem:gen} it equals
$\hup{z}$ for $z=i_{\hup{x}\cap\hup{y}}$; $z$ is the least upper bound of $x,y$.
Dually for meets.

(c) $[x,y]=\hup{x}\cap\hdn{y}$ is a face. Conversely, for a nonempty face $F$,
Lemma~\ref{lem:sep} gives $F=F^\sigma\cap F^\nu$ with
$F^\nu=\bigcap\{N\in\M_\nu:F\subseteq N\}\in\F_\nu$ and $F^\sigma\in\F_\sigma$ defined
dually. By Lemma~\ref{lem:gen}, $F^\nu=\hup{a}$, $F^\sigma=\hdn{b}$, so $F=[a,b]$.

(d) The facets are the maximal proper intervals. A proper interval $[x,y]$ has
$x\ne\sigma$ or $y\ne\nu$; in the first case it lies in $[a,\nu]$ for an atom $a\le x$,
in the second in $[\sigma,c]$ for a coatom $c\ge y$. By Lemma~\ref{lem:sep}(b) the facets are precisely
the maximal proper faces.

(e) $\hup{x}=\bigcap\{\hup{a}:a\text{ atom},a\le x\}=\hup{\bigvee\{a\}}$, whence
$x=\bigvee\{a\text{ atom}:a\le x\}$. Dually for coatoms.
\end{proof}

\begin{theorem}[faces of bipolar polytopes]\label{thm:faces}
Let $x<y$ in $(P,\le)$. Then the face $[x,y]$ with its faces inherited from $P$ is a
bipolar polytope with poles $x,y$, and its order is the restriction of $\le$.
\end{theorem}
\begin{proof}
Faces of the face $[x,y]$ are the faces of $P$ contained in it, i.e.\ a section of an
abstract polytope, so (D1), (D2) hold. Its facets are $[x',y]$ with $x'$ covering $x$,
and $[x,y']$ with $y'$ covered by $y$; they contain $y$ (resp.\ $x$) but not the other
pole $\Rightarrow$ (B1). North faces of $[x,y]$ are $[u,y]$, south faces $[x,v]$, and
$[u,y]\cap[x,v]\neq\emptyset\iff u\le v$. If $u,u'\le v$ then $u\vee u'\le v$, so
$[u,y]\cap[u',y]=[u\vee u',y]$ still meets $[x,v]$. Together with the dual
argumentation we get (B2).
\end{proof}

\begin{proposition}[vertex figures]\label{prop:vf}
The map $y\mapsto\hdn{y}=[\sigma,y]$ is an isomorphism of $(P,\le)$ onto $\F_\sigma$
(the faces containing $\sigma$, ordered by inclusion), i.e.\ onto the face lattice of
the vertex figure of $(P,\M)$ at $\sigma$. Dually
$\F_\nu\cong(P,\le)^{\mathrm{op}}$.
\end{proposition}
\begin{proof}
By Theorem~\ref{thm:lattice}, $\F_\sigma=\{\hdn{y}\}$ and
$\hdn{y}\subseteq\hdn{y'}\iff y\le y'$.
\end{proof}

\subsection{Examples}

\begin{example}\label{ex:trap}
Let $m\ge3$ and let $Z_m$ be the $m$-gonal trapezohedron with poles $s,n$ and
equatorial vertices $h_1,\dots,h_m,t_1,\dots,t_m$ (indices mod $m$), with south faces
$s\,h_i\,t_i\,h_{i+1}$ and north faces $n\,t_{i-1}\,h_i\,t_i$. Then $Z_m$ is a bipolar
polytope. Indeed (B1) is clear, and (B2) is a short check: the north faces meeting
the south edge $sh_i$ are $G_i=nt_{i-1}h_it_i$ and $P$; those meeting the
south quadrilateral $s h_it_ih_{i+1}$ are $G_i,G_{i+1}$, the edge $nt_i$ and $P$, and
$G_i\cap G_{i+1}=nt_i$ is in the family; the remaining cases are trivial, and the dual
statements are symmetric.
\end{example}

\begin{example}\label{ex:cube}
An $n$-dimensional cube is a bipolar polytope if $\sigma$, $\nu$ are opposite vertices, both (B1) and (B2) follow from the Boolean structure of the vertex
lattice. If it is a rectangular box (for instance the unit cube), it is even an
orthogonally bipolar polytope, since each edge at $\sigma$ is perpendicular to the facet through its
other endpoint not containing $\sigma$. A general combinatorial cube need
not be an orthogonally bipolar polytope, and a non-box orthogonally bipolar polytope of cube type exists
\cite[Figure~3]{Thales}.
\end{example}

\section{The combinatorial model and Theorem B}\label{sec:char}

The construction is suggested by the cutout of the apex of a pyramid
(\cite[Section~4.4]{Thales}): the faces of a pyramid over $Q$ containing its apex are the
cones over the faces of $Q$, so a cutout has one vertex for each face of $Q$, and its
facets are the sets $\{F:F\ni a\}={\uparrow}a$ for the vertices $a$ of $Q$ and
$\{F:F\subseteq\Phi\}={\downarrow}\Phi$ for the facets $\Phi$. For lattices this can be
written down directly.

\begin{definition}[$W(L)$]\label{def:W}
For a finite lattice $L$ with $0\neq1$ put
\[
W(L)=\bigl(L,\ \{{\uparrow}a: a\text{ an atom of }L\}\cup\{{\downarrow}c: c\text{ a coatom of }L\}\bigr),
\]
where ${\uparrow}a=\{z\ge a\}$, ${\downarrow}c=\{z\le c\}$. The poles are $0$ and $1$.
\end{definition}

\begin{lemma}\label{lem:Int}
Let $L$ be atomistic and coatomistic. The intersections of members of $W(L)$ are exactly
the intervals $[x,y]={\uparrow}x\cap{\downarrow}y$ ($x\le y$), the empty set and $L$
itself.
\end{lemma}
\begin{proof}
Atomistic: ${\uparrow}x=\bigcap\{{\uparrow}a:a\le x\}$ since
${\uparrow}a\cap{\uparrow}a'={\uparrow}(a\vee a')$. Dually
${\downarrow}y=\bigcap\{{\downarrow}c:c\ge y\}$.
\end{proof}

\begin{lemma}\label{lem:AP}
Let $L$ be an abstract polytope of rank $n$ (with rank function $\rk$ and extremal
elements $0_L,1_L$). Then $\widetilde I(L)=\Int(L)\cup\{\emptyset\}$, the set of
intervals of $L$ ordered by inclusion with a new least element, is an abstract polytope
of rank $n+1$ with the rank function $\rho(\emptyset)=-1$, $\rho([x,y])=\rk y-\rk x$.
\end{lemma}
\noindent
\begin{proof}
Ordered by reverse inclusion and with a new greatest element instead of $\emptyset$, the
intervals of $L$ form the antiprism of $L$ in the sense of Gleason and Hubard, which is an
abstract polytope of rank $n+1$ \cite[Theorem~2.1]{GleasonHubard}. Since the dual of an
abstract polytope is again an abstract polytope \cite{McMullenSchulte}, so is
$\widetilde I(L)$, and $\rho$ is its rank function: the singletons $[x,x]$ have rank $0$,
reversing inclusion raises $\rho$ by one in each covering step, and
$\rho([0_L,1_L])=n+1$.
\end{proof}

\begin{theorem}\label{thm:W}
If $L$ is an abstractly polytopal lattice, then $W(L)$ is a bipolar polytope with poles $0,1$, and
the lattice order from Theorem~\ref{thm:lattice} is the order of $L$.
\end{theorem}
\begin{proof}
(D1): the members $\uparrow a$, $\downarrow c$ are pairwise incomparable and cover $L$.
(D2): by Lemma~\ref{lem:Int}, $\F(W(L))=\widetilde I(L)$, which is an abstract polytope
by Lemma~\ref{lem:AP}; the singletons are $[x,x]$. (B1): $\uparrow a\ni1$, $\not\ni0$;
$\downarrow c\ni0$, $\not\ni1$. (B2): north faces are the filters ${\uparrow}x$ and south
faces the ideals ${\downarrow}y$; they meet iff $x\le y$. If $x,x'\le y$ then
${\uparrow}x\cap{\uparrow}x'={\uparrow}(x\vee x')$ and $x\vee x'\le y$; dually. Finally
$\hup{x}={\uparrow}x$ in $W(L)$.
\end{proof}

\begin{theorem}[Theorem B]\label{thm:B}
A finite lattice $L$ is abstractly polytopal if and only if $L\cong(P,\le)$ for a
bipolar polytope $(P,\M)$. If $L$ is abstractly polytopal, one may take
$(P,\M)=W(L)$. In particular, for every convex polytope $Q$, $W(\F(Q))$ is a bipolar
polytope whose vertices are the faces of $Q$, whose vertex lattice is $\F(Q)$, and whose
faces are the intervals of $\F(Q)$.
\end{theorem}
\begin{proof}
If $L\cong(P,\le)$, then by Proposition~\ref{prop:vf} $L$ is the lattice of faces
containing $\sigma$, i.e.\ a section of the face poset of the discrete polytope
$(P,\M)$. Sections of abstract polytopes are abstract polytopes, and $(P,\le)$ is
atomistic and coatomistic by Theorem~\ref{thm:lattice}(e); so $L$ is abstractly
polytopal. The converse is Theorem~\ref{thm:W}. The last statement follows, as
$\F(Q)$ is abstractly polytopal (Remark~\ref{rem:apl}(a)) and the faces of $W(L)$ are the
intervals of $L$ (Theorem~\ref{thm:lattice}).
\end{proof}

\begin{remark}\label{rem:Wdual}
Reversing the order exchanges the poles: $W(L^{\mathrm{op}})$ is obtained from $W(L)$
by exchanging the two families of facets; dually, this is
\cite[Proposition~2.2]{GleasonHubard}. A \emph{geometric realization} of $W(L)$
is a convex polytope $P$ with $\F(P)\cong\F(W(L))=\widetilde I(L)$
(Lemma~\ref{lem:Int}). Since the vertices and the facets of a polytope are the atoms
and the coatoms of its face lattice (Lemma~\ref{lem:sep}(b)), such an isomorphism
identifies the discrete polytope $(V(P),\{\text{vertex sets of facets}\})$ with
$W(L)$; so $P$ is a geometric bipolar polytope with vertex lattice $L$ by
Theorem~\ref{thm:W}. For $L=\F(Q)$ the vertices $w_x$ of $P$ are indexed by the faces
$x$ of $Q$ (including $\emptyset$ and $Q$), and its facets have the vertex sets
$\{w_x:a\in x\}$, $a$ a vertex of $Q$, and $\{w_x:x\subseteq\Phi\}$, $\Phi$ a facet of
$Q$. Geometric realizations of $W(\F(Q))$ are the subject of
Section~\ref{sec:convex}.
\end{remark}

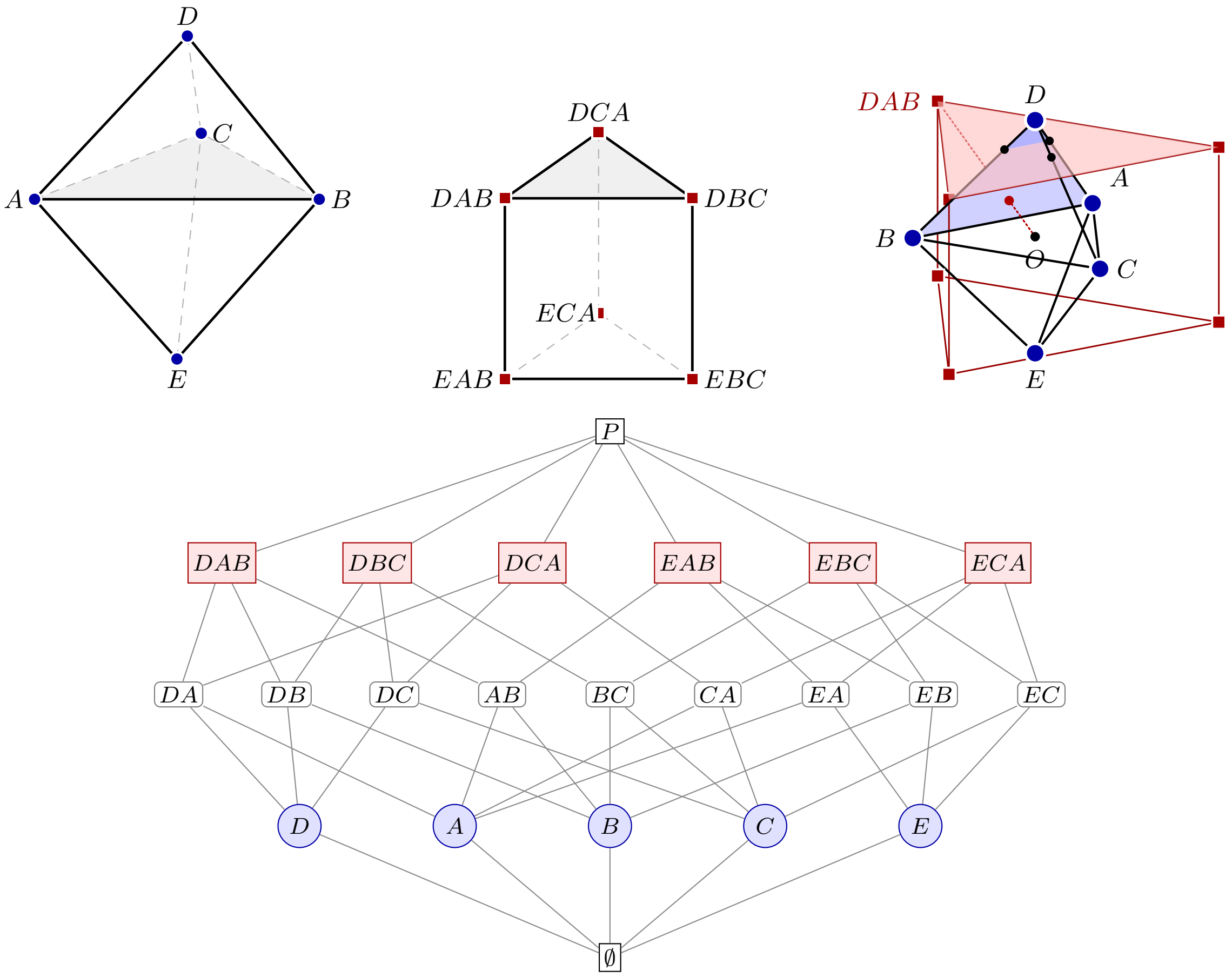
\begin{figure}[htbp]
\centering
\begin{tikzpicture}[scale=0.85,line join=round,line cap=round,
                    every node/.style={font=\footnotesize}]
  \coordinate (D) at (0.15,2.45);
  \coordinate (A) at (-2.05,0.10);
  \coordinate (B) at (2.05,0.10);
  \coordinate (C) at (0.35,1.05);
  \coordinate (E) at (0.00,-2.20);
  \fill[gray!12] (A)--(B)--(C)--cycle;
  \draw[gray!60,dashed] (D)--(C) (E)--(C) (A)--(C) (B)--(C);
  \draw[line width=0.9pt] (D)--(A) (D)--(B) (E)--(A) (E)--(B) (A)--(B);
  \foreach \p in {D,A,B,C,E}{%
    \fill[white] (\p) circle (3.2pt);
    \fill[blue!65!black] (\p) circle (2.3pt);}
  \node[above] at (D) {$D$};
  \node[left]  at (A) {$A$};
  \node[right] at (B) {$B$};
  \node[right] at (C) {$C$};
  \node[below] at (E) {$E$};
\end{tikzpicture}\hspace{0.7cm}%
\begin{tikzpicture}[scale=0.85,line join=round,line cap=round,
                    every node/.style={font=\footnotesize}]
  \coordinate (t1) at (0,2.8);
  \coordinate (t2) at (2.7,2.8);
  \coordinate (t3) at (1.35,3.75);
  \coordinate (b1) at (0,0.2);
  \coordinate (b2) at (2.7,0.2);
  \coordinate (b3) at (1.35,1.15);
  \fill[gray!12] (t1)--(t2)--(t3)--cycle;
  \draw[gray!60,dashed] (t3)--(b3) (b1)--(b3) (b2)--(b3);
  \draw[line width=0.9pt] (t1)--(t2)--(t3)--cycle (t1)--(b1) (t2)--(b2) (b1)--(b2);
  \foreach \p in {t1,t2,t3,b1,b2,b3}{%
    \begin{scope}[shift={(\p)}]
      \fill[white] (-3.0pt,-3.0pt) rectangle (3.0pt,3.0pt);
      \fill[red!65!black] (-2.1pt,-2.1pt) rectangle (2.1pt,2.1pt);
    \end{scope}}
  \node[left]  at (t1) {$DAB$};
  \node[right] at (t2) {$DBC$};
  \node[above] at (t3) {$DCA$};
  \node[left]  at (b1) {$EAB$};
  \node[right] at (b2) {$EBC$};
  \node[left,fill=white,inner sep=0.5pt] at (b3) {$ECA$};
\end{tikzpicture}
\hspace{0.7cm}%
\begin{tikzpicture}[scale=1.5,line join=round,line cap=round]
\draw[red!60!black,line width=0.55pt] (-0.796,1.106)--(-0.703,0.304);
\draw[red!60!black,line width=0.55pt] (-0.703,0.304)--(-0.676,0.309);
\draw[red!60!black,line width=0.55pt] (0.315,0.501)--(1.499,0.730);
\draw[red!60!black,line width=0.55pt] (1.499,0.730)--(-0.796,1.106);
\draw[red!60!black,line width=0.55pt] (-0.796,-0.320)--(-0.703,-1.123);
\draw[red!60!black,line width=0.55pt] (-0.703,-1.123)--(1.499,-0.697);
\draw[red!60!black,line width=0.55pt] (1.499,-0.697)--(-0.796,-0.320);
\draw[red!60!black,line width=0.55pt] (-0.796,1.106)--(-0.796,-0.320);
\draw[red!60!black,line width=0.55pt] (-0.703,0.304)--(-0.703,0.268);
\draw[red!60!black,line width=0.55pt] (-0.703,0.054)--(-0.703,-1.123);
\draw[red!60!black,line width=0.55pt] (1.499,0.730)--(1.499,-0.697);
\draw[black,line width=0.8pt] (0.312,0.238)--(0.530,-0.262);
\draw[black,line width=0.8pt] (0.000,-0.951)--(0.469,0.273);
\draw[black,line width=0.8pt] (0.000,-0.951)--(-0.999,-0.011);
\draw[black,line width=0.8pt] (0.000,-0.951)--(0.530,-0.262);
\draw[black,line width=0.8pt] (-0.999,-0.011)--(0.530,-0.262);
\draw[black,line width=0.8pt] (0.530,-0.262)--(0.469,0.273);
\draw[red!70!black,densely dotted,line width=0.6pt] (0.000,0.000)--(-0.119,0.166);
\draw[red!70!black,densely dotted,line width=0.6pt] (-0.209,0.290)--(-0.796,1.106);
\fill[white] (0.530,-0.262) circle (2.6pt); \fill[blue!65!black] (0.530,-0.262) circle (1.9pt);
\fill[white] (0.000,-0.951) circle (2.6pt); \fill[blue!65!black] (0.000,-0.951) circle (1.9pt);
\fill[black] (0.000,0.000) circle (1.1pt);
\node[draw=white,fill=red!65!black,line width=0.7pt,minimum size=4.6pt,inner sep=0pt] at (-0.796,1.106) {};
\node[draw=white,fill=red!65!black,line width=0.7pt,minimum size=4.6pt,inner sep=0pt] at (-0.703,0.304) {};
\node[draw=white,fill=red!65!black,line width=0.7pt,minimum size=4.6pt,inner sep=0pt] at (1.499,0.730) {};
\node[draw=white,fill=red!65!black,line width=0.7pt,minimum size=4.6pt,inner sep=0pt] at (-0.796,-0.320) {};
\node[draw=white,fill=red!65!black,line width=0.7pt,minimum size=4.6pt,inner sep=0pt] at (-0.703,-1.123) {};
\node[draw=white,fill=red!65!black,line width=0.7pt,minimum size=4.6pt,inner sep=0pt] at (1.499,-0.697) {};
\fill[blue!18] (0.469,0.273)--(-0.999,-0.011)--(-0.250,0.711)--(0.117,0.782)--cycle;
\draw[red!60!black,line width=0.55pt] (-0.676,0.309)--(0.315,0.501);
\draw[red!60!black,line width=0.55pt] (-0.703,0.268)--(-0.703,0.054);
\draw[black,line width=0.8pt] (0.117,0.782)--(0.469,0.273);
\draw[black,line width=0.8pt] (-0.250,0.711)--(-0.999,-0.011);
\draw[black,line width=0.8pt] (0.133,0.648)--(0.312,0.238);
\draw[black,line width=0.8pt] (0.469,0.273)--(-0.999,-0.011);
\draw[red!70!black,densely dotted,line width=0.6pt] (-0.119,0.166)--(-0.209,0.290);
\fill[white] (0.469,0.273) circle (2.6pt); \fill[blue!65!black] (0.469,0.273) circle (1.9pt);
\fill[white] (-0.999,-0.011) circle (2.6pt); \fill[blue!65!black] (-0.999,-0.011) circle (1.9pt);
\fill[red!70!black] (-0.212,0.295) circle (1.1pt);
\fill[red!25,opacity=0.6] (-0.796,1.106)--(-0.703,0.304)--(1.499,0.730)--cycle;
\draw[black,line width=0.8pt] (0.080,0.769)--(0.133,0.648);
\fill[black] (0.133,0.648) circle (1.0pt);
\fill[blue!30] (0.000,0.951)--(0.117,0.782)--(-0.250,0.711)--cycle;
\draw[black,line width=0.8pt] (0.000,0.951)--(0.117,0.782);
\draw[black,line width=0.8pt] (0.000,0.951)--(-0.250,0.711);
\draw[black,line width=0.8pt] (0.000,0.951)--(0.080,0.769);
\fill[white] (0.000,0.951) circle (2.6pt); \fill[blue!65!black] (0.000,0.951) circle (1.9pt);
\fill[black] (0.117,0.782) circle (1.0pt);
\fill[black] (-0.250,0.711) circle (1.0pt);
\node[font=\footnotesize,above right=2pt] at (0.469,0.273) {$A$};
\node[font=\footnotesize,left=2pt] at (-0.999,-0.011) {$B$};
\node[font=\footnotesize,right=2pt] at (0.530,-0.262) {$C$};
\node[font=\footnotesize,above=2pt] at (0.000,0.951) {$D$};
\node[font=\footnotesize,below=2pt] at (0.000,-0.951) {$E$};
\node[font=\footnotesize,red!60!black,left=2pt] at (-0.796,1.106) {$DAB$};
\node[font=\footnotesize,below=1pt] at (0.000,0.000) {$O$};
\end{tikzpicture}
\par\medskip
\begin{tikzpicture}[xscale=1.0,yscale=1.15,
   every node/.style={font=\scriptsize,inner sep=1.6pt},
   ver/.style={circle,draw=blue!65!black,fill=blue!12,minimum size=15pt},
   fac/.style={rectangle,draw=red!65!black,fill=red!10,minimum size=14pt},
   edg/.style={rectangle,rounded corners=2pt,draw=black!50,fill=white},
   ext/.style={rectangle,draw=black,fill=white},
   hl/.style={black!45,line width=0.4pt}]
  \node[ext] (bot) at (0,0) {$\emptyset$};
  \node[ver] (vD) at (-3.80,1.4) {$D$};
  \node[ver] (vA) at (-1.90,1.4) {$A$};
  \node[ver] (vB) at (0.00,1.4) {$B$};
  \node[ver] (vC) at (1.90,1.4) {$C$};
  \node[ver] (vE) at (3.80,1.4) {$E$};
  \node[edg] (eDA) at (-5.28,2.8) {$DA$};
  \node[edg] (eDB) at (-3.96,2.8) {$DB$};
  \node[edg] (eDC) at (-2.64,2.8) {$DC$};
  \node[edg] (eAB) at (-1.32,2.8) {$AB$};
  \node[edg] (eBC) at (0.00,2.8) {$BC$};
  \node[edg] (eCA) at (1.32,2.8) {$CA$};
  \node[edg] (eEA) at (2.64,2.8) {$EA$};
  \node[edg] (eEB) at (3.96,2.8) {$EB$};
  \node[edg] (eEC) at (5.28,2.8) {$EC$};
  \node[fac] (fDAB) at (-4.75,4.2) {$DAB$};
  \node[fac] (fDBC) at (-2.85,4.2) {$DBC$};
  \node[fac] (fDCA) at (-0.95,4.2) {$DCA$};
  \node[fac] (fEAB) at (0.95,4.2) {$EAB$};
  \node[fac] (fEBC) at (2.85,4.2) {$EBC$};
  \node[fac] (fECA) at (4.75,4.2) {$ECA$};
  \node[ext] (top) at (0,5.6) {$P$};
  \draw[hl] (bot)--(vD) (bot)--(vA) (bot)--(vB) (bot)--(vC) (bot)--(vE) (vD)--(eDA);
  \draw[hl] (vA)--(eDA) (vD)--(eDB) (vB)--(eDB) (vD)--(eDC) (vC)--(eDC) (vA)--(eAB);
  \draw[hl] (vB)--(eAB) (vB)--(eBC) (vC)--(eBC) (vC)--(eCA) (vA)--(eCA) (vE)--(eEA);
  \draw[hl] (vA)--(eEA) (vE)--(eEB) (vB)--(eEB) (vE)--(eEC) (vC)--(eEC) (eDA)--(fDAB);
  \draw[hl] (eDB)--(fDAB) (eAB)--(fDAB) (eDB)--(fDBC) (eDC)--(fDBC) (eBC)--(fDBC) (eDA)--(fDCA);
  \draw[hl] (eDC)--(fDCA) (eCA)--(fDCA) (eAB)--(fEAB) (eEA)--(fEAB) (eEB)--(fEAB) (eBC)--(fEBC);
  \draw[hl] (eEB)--(fEBC) (eEC)--(fEBC) (eCA)--(fECA) (eEA)--(fECA) (eEC)--(fECA) (fDAB)--(top);
  \draw[hl] (fDBC)--(top) (fDCA)--(top) (fEAB)--(top) (fEBC)--(top) (fECA)--(top);
\end{tikzpicture}
\caption{Theorem~B for a three-dimensional base. (a) The triangular bipyramid
$P=\conv\{A,B,C,D,E\}$: $5$ vertices, $9$ edges, $6$ triangular facets; the triangle
$ABC$ is not a face. (b) The dual prism $P^{*}$, whose six vertices are the facets of
$P$ (labelled by the vertices they span) and whose five facets are dual to the five
vertices of $P$. (c) The duality: $P$ (blue round vertices) together with a
homothetic copy of $P^{*}$ (red square vertices), the polar of $P$ with respect to a sphere centred at $O$. The
vertex of $P^{*}$ labelled $DAB$ lies on the normal through $O$ (dotted) of the facet
$DAB$ of $P$ (shaded blue, darker), and the facet of $P^{*}$ dual to the vertex $D$ is the
triangle shaded red (lighter); it lies in front of the lower part of $P$ and is pierced by the
three edges from $D$ (dots), so that only the small cap at $D$ rises above it. (d) The Hasse diagram of $\F(P)$ ($22$ elements), with
the atoms (vertices) marked as blue circles and the coatoms (facets) marked as red
squares. By Theorem~\ref{thm:B} this lattice is the vertex lattice of
the discrete bipolar polytope $W(\F(P))$ with poles $\emptyset$ and $P$;
its facets are the $5$ filters ${\uparrow}a$ and the $6$ ideals
${\downarrow}\Phi$.}
\label{fig:thmB}
\end{figure}

\begin{example}\label{ex:thmB}
Let $P$ be the triangular bipyramid of Figure~\ref{fig:thmB}. Its face lattice has
$22$ elements, $5$ atoms (the vertices) and $6$ coatoms (the facets); its
Hasse diagram is Figure~\ref{fig:thmB}(d). By
Theorem~\ref{thm:B} it is the vertex lattice of $W(\F(P))$.
Explicitly, $W(\F(P))$ is the discrete bipolar polytope whose vertex set
consists of the $22$ faces of $P$, with the poles $\sigma=\emptyset$ and
$\nu=P$, and whose $11$ facets are
\[
{\uparrow}a=\{F\in\F(P):a\in F\}\quad(a\in\{A,B,C,D,E\}),\qquad
{\downarrow}\Phi=\{F\in\F(P):F\subseteq\Phi\}\quad(\Phi\text{ a facet}).
\]
For instance ${\uparrow}D$ consists of $D$, the edges $DA,DB,DC$, the
facets $DAB,DBC,DCA$ and $P$ ($8$ elements), and ${\downarrow}DAB$
consists of $\emptyset$, the vertices $D,A,B$, the edges $DA,DB,AB$ and
$DAB$ ($8$ elements). Every facet contains exactly one pole, and
by Lemma~\ref{lem:Int} the faces of $W(\F(P))$ are the intervals of
$\F(P)$; it is an abstract polytope of rank $4$ by Lemma~\ref{lem:AP}.
The order of Theorem~\ref{thm:lattice} on its vertices is inclusion of
faces of $P$, so its vertex lattice is the lattice drawn in
Figure~\ref{fig:thmB}(d). A convex $4$-dimensional realization of
$W(\F(P))$ is given in Example~\ref{ex:cayley-bip}.
\end{example}

\section{Convex realizations}\label{sec:convex}

Throughout this section $Q\subset\R^d$ is a $d$-polytope with $O\in\inte Q$, $Q^\circ$
is its polar, and $F\mapsto F^*$ is the duality of Section~\ref{sec:notation}. Recall
that $N_Q(H)=\cone(H^*)$ and, dually, $N_{Q^\circ}(G^*)=\cone(G)$: the cones
$N_Q(H)$ form the normal fan of $Q$, and the cones $\cone(G)$ over the proper faces
form its face fan, which is the normal fan of $Q^\circ$.

\subsection{The interval condition and the Cayley polytope}\label{sec:cayley}

\begin{definition}[interval condition]\label{def:intcond}
$Q$ satisfies the \emph{interval condition} if for all nonempty faces $H$ and all
nonempty proper faces $G$ of $Q$
\[
\relint N_Q(H)\cap\relint\cone(G)\neq\emptyset\iff H\subseteq G .
\tag{I}
\]
\end{definition}

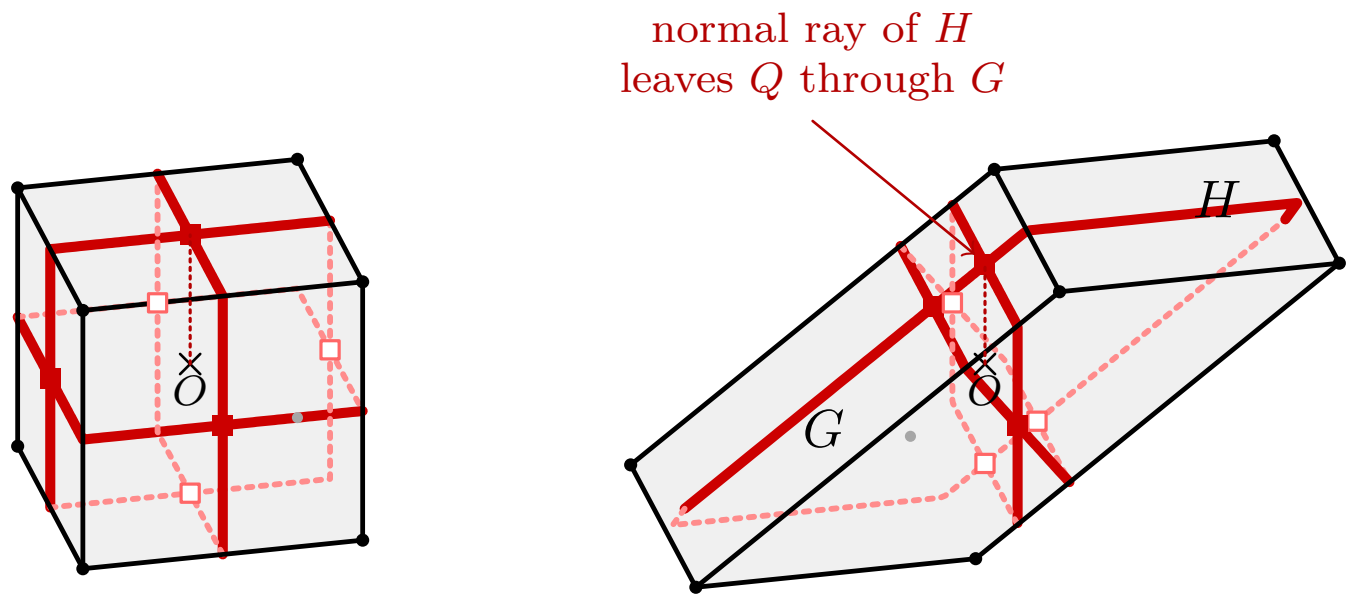
\begin{figure}[htbp]
\centering
\begin{tikzpicture}[scale=0.95,line join=round,line cap=round]
\draw[black!35,dashed] (1.201,-1.226)--(0.747,-0.373);
\draw[black!35,dashed] (0.747,-0.373)--(0.747,1.425);
\draw[black!35,dashed] (-1.201,-0.572)--(0.747,-0.373);
\fill[black!6] (-0.747,-1.425)--(-1.201,-0.572)--(-1.201,1.226)--(-0.747,0.373)--cycle;
\fill[black!6] (-0.747,-1.425)--(1.201,-1.226)--(1.201,0.572)--(-0.747,0.373)--cycle;
\fill[black!6] (-0.747,0.373)--(1.201,0.572)--(0.747,1.425)--(-1.201,1.226)--cycle;
\draw[red!80!black,line width=1.9pt] (-0.974,-0.099)--(-0.969,-0.108)--(-0.965,-0.116)--(-0.960,-0.125)--(-0.956,-0.133)--(-0.951,-0.141)--(-0.947,-0.150)--(-0.942,-0.158)--(-0.938,-0.167)--(-0.933,-0.176)--(-0.929,-0.184)--(-0.924,-0.193)--(-0.919,-0.202)--(-0.915,-0.211)--(-0.910,-0.220)--(-0.905,-0.229)--(-0.900,-0.238)--(-0.895,-0.247)--(-0.890,-0.257)--(-0.885,-0.266)--(-0.880,-0.276)--(-0.874,-0.286)--(-0.869,-0.296)--(-0.864,-0.306)--(-0.858,-0.317)--(-0.852,-0.327)--(-0.847,-0.338)--(-0.841,-0.349)--(-0.835,-0.361)--(-0.828,-0.372)--(-0.822,-0.384)--(-0.815,-0.397)--(-0.809,-0.409)--(-0.802,-0.422)--(-0.795,-0.436)--(-0.787,-0.449)--(-0.780,-0.464)--(-0.772,-0.478)--(-0.764,-0.494)--(-0.755,-0.509)--(-0.747,-0.526)--(-0.709,-0.522)--(-0.673,-0.518)--(-0.638,-0.515)--(-0.605,-0.511)--(-0.572,-0.508)--(-0.540,-0.505)--(-0.510,-0.502)--(-0.480,-0.499)--(-0.451,-0.496)--(-0.423,-0.493)--(-0.396,-0.490)--(-0.370,-0.487)--(-0.344,-0.485)--(-0.318,-0.482)--(-0.293,-0.480)--(-0.269,-0.477)--(-0.245,-0.475)--(-0.222,-0.472)--(-0.199,-0.470)--(-0.176,-0.468)--(-0.154,-0.465)--(-0.132,-0.463)--(-0.110,-0.461)--(-0.089,-0.459)--(-0.068,-0.457)--(-0.047,-0.454)--(-0.027,-0.452)--(-0.007,-0.450)--(0.014,-0.448)--(0.034,-0.446)--(0.053,-0.444)--(0.073,-0.442)--(0.093,-0.440)--(0.112,-0.438)--(0.131,-0.436)--(0.151,-0.434)--(0.170,-0.432)--(0.189,-0.430)--(0.208,-0.428)--(0.227,-0.426);
\draw[red!80!black,line width=1.9pt] (-0.974,-0.099)--(-0.978,-0.091)--(-0.983,-0.083)--(-0.987,-0.074)--(-0.992,-0.066)--(-0.996,-0.057)--(-1.001,-0.049)--(-1.005,-0.041)--(-1.010,-0.032)--(-1.014,-0.023)--(-1.019,-0.015)--(-1.024,-0.006)--(-1.028,0.003)--(-1.033,0.012)--(-1.038,0.021)--(-1.043,0.030)--(-1.048,0.039)--(-1.053,0.048)--(-1.058,0.058)--(-1.063,0.067)--(-1.068,0.077)--(-1.073,0.087)--(-1.079,0.097)--(-1.084,0.107)--(-1.090,0.118)--(-1.095,0.128)--(-1.101,0.139)--(-1.107,0.150)--(-1.113,0.162)--(-1.119,0.173)--(-1.126,0.185)--(-1.132,0.198)--(-1.139,0.210)--(-1.146,0.223)--(-1.153,0.237)--(-1.160,0.250)--(-1.168,0.265)--(-1.176,0.279)--(-1.184,0.295)--(-1.192,0.310)--(-1.201,0.327);
\draw[red!45,line width=1.1pt,dash pattern=on 1pt off 1.6pt] (-1.201,0.327)--(-1.164,0.331)--(-1.127,0.334)--(-1.093,0.338)--(-1.059,0.341)--(-1.026,0.345)--(-0.995,0.348)--(-0.964,0.351)--(-0.935,0.354)--(-0.906,0.357)--(-0.878,0.360)--(-0.851,0.363)--(-0.824,0.365)--(-0.798,0.368)--(-0.773,0.371)--(-0.748,0.373)--(-0.723,0.376)--(-0.700,0.378)--(-0.676,0.380)--(-0.653,0.383)--(-0.631,0.385)--(-0.608,0.387)--(-0.586,0.390)--(-0.565,0.392)--(-0.544,0.394)--(-0.523,0.396)--(-0.502,0.398)--(-0.481,0.400)--(-0.461,0.402)--(-0.441,0.405)--(-0.421,0.407)--(-0.401,0.409)--(-0.381,0.411)--(-0.362,0.413)--(-0.342,0.415)--(-0.323,0.417)--(-0.304,0.419)--(-0.285,0.420)--(-0.265,0.422)--(-0.246,0.424)--(-0.227,0.426);
\draw[red!80!black,line width=1.9pt] (-0.974,-0.099)--(-0.974,-0.117)--(-0.974,-0.135)--(-0.974,-0.152)--(-0.974,-0.170)--(-0.974,-0.188)--(-0.974,-0.206)--(-0.974,-0.224)--(-0.974,-0.242)--(-0.974,-0.260)--(-0.974,-0.278)--(-0.974,-0.297)--(-0.974,-0.315)--(-0.974,-0.334)--(-0.974,-0.353)--(-0.974,-0.372)--(-0.974,-0.392)--(-0.974,-0.411)--(-0.974,-0.431)--(-0.974,-0.451)--(-0.974,-0.472)--(-0.974,-0.493)--(-0.974,-0.514)--(-0.974,-0.536)--(-0.974,-0.558)--(-0.974,-0.580)--(-0.974,-0.603)--(-0.974,-0.626)--(-0.974,-0.650)--(-0.974,-0.675)--(-0.974,-0.700)--(-0.974,-0.726)--(-0.974,-0.753)--(-0.974,-0.780)--(-0.974,-0.808)--(-0.974,-0.837)--(-0.974,-0.867)--(-0.974,-0.898)--(-0.974,-0.931)--(-0.974,-0.964)--(-0.974,-0.999);
\draw[red!45,line width=1.1pt,dash pattern=on 1pt off 1.6pt] (-0.974,-0.999)--(-0.936,-0.995)--(-0.900,-0.991)--(-0.865,-0.987)--(-0.832,-0.984)--(-0.799,-0.981)--(-0.768,-0.978)--(-0.737,-0.974)--(-0.708,-0.971)--(-0.679,-0.968)--(-0.651,-0.966)--(-0.623,-0.963)--(-0.597,-0.960)--(-0.571,-0.957)--(-0.545,-0.955)--(-0.521,-0.952)--(-0.496,-0.950)--(-0.472,-0.947)--(-0.449,-0.945)--(-0.426,-0.943)--(-0.403,-0.940)--(-0.381,-0.938)--(-0.359,-0.936)--(-0.338,-0.934)--(-0.316,-0.931)--(-0.295,-0.929)--(-0.275,-0.927)--(-0.254,-0.925)--(-0.234,-0.923)--(-0.214,-0.921)--(-0.194,-0.919)--(-0.174,-0.917)--(-0.154,-0.915)--(-0.135,-0.913)--(-0.115,-0.911)--(-0.096,-0.909)--(-0.077,-0.907)--(-0.057,-0.905)--(-0.038,-0.903)--(-0.019,-0.901)--(0.000,-0.899);
\draw[red!80!black,line width=1.9pt] (-0.974,-0.099)--(-0.974,-0.082)--(-0.974,-0.064)--(-0.974,-0.046)--(-0.974,-0.029)--(-0.974,-0.011)--(-0.974,0.007)--(-0.974,0.025)--(-0.974,0.043)--(-0.974,0.061)--(-0.974,0.079)--(-0.974,0.098)--(-0.974,0.116)--(-0.974,0.135)--(-0.974,0.154)--(-0.974,0.173)--(-0.974,0.193)--(-0.974,0.212)--(-0.974,0.232)--(-0.974,0.252)--(-0.974,0.273)--(-0.974,0.294)--(-0.974,0.315)--(-0.974,0.337)--(-0.974,0.359)--(-0.974,0.381)--(-0.974,0.404)--(-0.974,0.427)--(-0.974,0.451)--(-0.974,0.476)--(-0.974,0.501)--(-0.974,0.527)--(-0.974,0.554)--(-0.974,0.581)--(-0.974,0.609)--(-0.974,0.638)--(-0.974,0.668)--(-0.974,0.699)--(-0.974,0.732)--(-0.974,0.765)--(-0.974,0.800)--(-0.936,0.803)--(-0.900,0.807)--(-0.865,0.811)--(-0.832,0.814)--(-0.799,0.817)--(-0.768,0.821)--(-0.737,0.824)--(-0.708,0.827)--(-0.679,0.830)--(-0.651,0.833)--(-0.623,0.835)--(-0.597,0.838)--(-0.571,0.841)--(-0.545,0.843)--(-0.521,0.846)--(-0.496,0.848)--(-0.472,0.851)--(-0.449,0.853)--(-0.426,0.856)--(-0.403,0.858)--(-0.381,0.860)--(-0.359,0.862)--(-0.338,0.865)--(-0.316,0.867)--(-0.295,0.869)--(-0.275,0.871)--(-0.254,0.873)--(-0.234,0.875)--(-0.214,0.877)--(-0.194,0.879)--(-0.174,0.881)--(-0.154,0.883)--(-0.135,0.885)--(-0.115,0.887)--(-0.096,0.889)--(-0.077,0.891)--(-0.057,0.893)--(-0.038,0.895)--(-0.019,0.897)--(0.000,0.899);
\draw[red!45,line width=1.1pt,dash pattern=on 1pt off 1.6pt] (0.974,0.099)--(0.978,0.091)--(0.983,0.083)--(0.987,0.074)--(0.992,0.066)--(0.996,0.057)--(1.001,0.049)--(1.005,0.041)--(1.010,0.032)--(1.014,0.023)--(1.019,0.015)--(1.024,0.006)--(1.028,-0.003)--(1.033,-0.012)--(1.038,-0.021)--(1.043,-0.030)--(1.048,-0.039)--(1.053,-0.048)--(1.058,-0.058)--(1.063,-0.067)--(1.068,-0.077)--(1.073,-0.087)--(1.079,-0.097)--(1.084,-0.107)--(1.090,-0.118)--(1.095,-0.128)--(1.101,-0.139)--(1.107,-0.150)--(1.113,-0.162)--(1.119,-0.173)--(1.126,-0.185)--(1.132,-0.198)--(1.139,-0.210)--(1.146,-0.223)--(1.153,-0.237)--(1.160,-0.250)--(1.168,-0.265)--(1.176,-0.279)--(1.184,-0.295)--(1.192,-0.310)--(1.201,-0.327);
\draw[red!80!black,line width=1.9pt] (1.201,-0.327)--(1.164,-0.331)--(1.127,-0.334)--(1.093,-0.338)--(1.059,-0.341)--(1.026,-0.345)--(0.995,-0.348)--(0.964,-0.351)--(0.935,-0.354)--(0.906,-0.357)--(0.878,-0.360)--(0.851,-0.363)--(0.824,-0.365)--(0.798,-0.368)--(0.773,-0.371)--(0.748,-0.373)--(0.723,-0.376)--(0.700,-0.378)--(0.676,-0.380)--(0.653,-0.383)--(0.631,-0.385)--(0.608,-0.387)--(0.586,-0.390)--(0.565,-0.392)--(0.544,-0.394)--(0.523,-0.396)--(0.502,-0.398)--(0.481,-0.400)--(0.461,-0.402)--(0.441,-0.405)--(0.421,-0.407)--(0.401,-0.409)--(0.381,-0.411)--(0.362,-0.413)--(0.342,-0.415)--(0.323,-0.417)--(0.304,-0.419)--(0.285,-0.420)--(0.265,-0.422)--(0.246,-0.424)--(0.227,-0.426);
\draw[red!45,line width=1.1pt,dash pattern=on 1pt off 1.6pt] (0.974,0.099)--(0.969,0.108)--(0.965,0.116)--(0.960,0.125)--(0.956,0.133)--(0.951,0.141)--(0.947,0.150)--(0.942,0.158)--(0.938,0.167)--(0.933,0.176)--(0.929,0.184)--(0.924,0.193)--(0.919,0.202)--(0.915,0.211)--(0.910,0.220)--(0.905,0.229)--(0.900,0.238)--(0.895,0.247)--(0.890,0.257)--(0.885,0.266)--(0.880,0.276)--(0.874,0.286)--(0.869,0.296)--(0.864,0.306)--(0.858,0.317)--(0.852,0.327)--(0.847,0.338)--(0.841,0.349)--(0.835,0.361)--(0.828,0.372)--(0.822,0.384)--(0.815,0.397)--(0.809,0.409)--(0.802,0.422)--(0.795,0.436)--(0.787,0.449)--(0.780,0.464)--(0.772,0.478)--(0.764,0.494)--(0.755,0.509)--(0.747,0.526)--(0.709,0.522)--(0.673,0.518)--(0.638,0.515)--(0.605,0.511)--(0.572,0.508)--(0.540,0.505)--(0.510,0.502)--(0.480,0.499)--(0.451,0.496)--(0.423,0.493)--(0.396,0.490)--(0.370,0.487)--(0.344,0.485)--(0.318,0.482)--(0.293,0.480)--(0.269,0.477)--(0.245,0.475)--(0.222,0.472)--(0.199,0.470)--(0.176,0.468)--(0.154,0.465)--(0.132,0.463)--(0.110,0.461)--(0.089,0.459)--(0.068,0.457)--(0.047,0.454)--(0.027,0.452)--(0.007,0.450)--(-0.014,0.448)--(-0.034,0.446)--(-0.053,0.444)--(-0.073,0.442)--(-0.093,0.440)--(-0.112,0.438)--(-0.131,0.436)--(-0.151,0.434)--(-0.170,0.432)--(-0.189,0.430)--(-0.208,0.428)--(-0.227,0.426);
\draw[red!45,line width=1.1pt,dash pattern=on 1pt off 1.6pt] (0.974,0.099)--(0.974,0.082)--(0.974,0.064)--(0.974,0.046)--(0.974,0.029)--(0.974,0.011)--(0.974,-0.007)--(0.974,-0.025)--(0.974,-0.043)--(0.974,-0.061)--(0.974,-0.079)--(0.974,-0.098)--(0.974,-0.116)--(0.974,-0.135)--(0.974,-0.154)--(0.974,-0.173)--(0.974,-0.193)--(0.974,-0.212)--(0.974,-0.232)--(0.974,-0.252)--(0.974,-0.273)--(0.974,-0.294)--(0.974,-0.315)--(0.974,-0.337)--(0.974,-0.359)--(0.974,-0.381)--(0.974,-0.404)--(0.974,-0.427)--(0.974,-0.451)--(0.974,-0.476)--(0.974,-0.501)--(0.974,-0.527)--(0.974,-0.554)--(0.974,-0.581)--(0.974,-0.609)--(0.974,-0.638)--(0.974,-0.668)--(0.974,-0.699)--(0.974,-0.732)--(0.974,-0.765)--(0.974,-0.800)--(0.936,-0.803)--(0.900,-0.807)--(0.865,-0.811)--(0.832,-0.814)--(0.799,-0.817)--(0.768,-0.821)--(0.737,-0.824)--(0.708,-0.827)--(0.679,-0.830)--(0.651,-0.833)--(0.623,-0.835)--(0.597,-0.838)--(0.571,-0.841)--(0.545,-0.843)--(0.521,-0.846)--(0.496,-0.848)--(0.472,-0.851)--(0.449,-0.853)--(0.426,-0.856)--(0.403,-0.858)--(0.381,-0.860)--(0.359,-0.862)--(0.338,-0.865)--(0.316,-0.867)--(0.295,-0.869)--(0.275,-0.871)--(0.254,-0.873)--(0.234,-0.875)--(0.214,-0.877)--(0.194,-0.879)--(0.174,-0.881)--(0.154,-0.883)--(0.135,-0.885)--(0.115,-0.887)--(0.096,-0.889)--(0.077,-0.891)--(0.057,-0.893)--(0.038,-0.895)--(0.019,-0.897)--(0.000,-0.899);
\draw[red!45,line width=1.1pt,dash pattern=on 1pt off 1.6pt] (0.974,0.099)--(0.974,0.117)--(0.974,0.135)--(0.974,0.152)--(0.974,0.170)--(0.974,0.188)--(0.974,0.206)--(0.974,0.224)--(0.974,0.242)--(0.974,0.260)--(0.974,0.278)--(0.974,0.297)--(0.974,0.315)--(0.974,0.334)--(0.974,0.353)--(0.974,0.372)--(0.974,0.392)--(0.974,0.411)--(0.974,0.431)--(0.974,0.451)--(0.974,0.472)--(0.974,0.493)--(0.974,0.514)--(0.974,0.536)--(0.974,0.558)--(0.974,0.580)--(0.974,0.603)--(0.974,0.626)--(0.974,0.650)--(0.974,0.675)--(0.974,0.700)--(0.974,0.726)--(0.974,0.753)--(0.974,0.780)--(0.974,0.808)--(0.974,0.837)--(0.974,0.867)--(0.974,0.898)--(0.974,0.931)--(0.974,0.964)--(0.974,0.999);
\draw[red!80!black,line width=1.9pt] (0.974,0.999)--(0.936,0.995)--(0.900,0.991)--(0.865,0.987)--(0.832,0.984)--(0.799,0.981)--(0.768,0.978)--(0.737,0.974)--(0.708,0.971)--(0.679,0.968)--(0.651,0.966)--(0.623,0.963)--(0.597,0.960)--(0.571,0.957)--(0.545,0.955)--(0.521,0.952)--(0.496,0.950)--(0.472,0.947)--(0.449,0.945)--(0.426,0.943)--(0.403,0.940)--(0.381,0.938)--(0.359,0.936)--(0.338,0.934)--(0.316,0.931)--(0.295,0.929)--(0.275,0.927)--(0.254,0.925)--(0.234,0.923)--(0.214,0.921)--(0.194,0.919)--(0.174,0.917)--(0.154,0.915)--(0.135,0.913)--(0.115,0.911)--(0.096,0.909)--(0.077,0.907)--(0.057,0.905)--(0.038,0.903)--(0.019,0.901)--(0.000,0.899);
\draw[red!80!black,line width=1.9pt] (0.227,-0.426)--(0.227,-0.444)--(0.227,-0.462)--(0.227,-0.479)--(0.227,-0.497)--(0.227,-0.515)--(0.227,-0.533)--(0.227,-0.551)--(0.227,-0.569)--(0.227,-0.587)--(0.227,-0.605)--(0.227,-0.624)--(0.227,-0.642)--(0.227,-0.661)--(0.227,-0.680)--(0.227,-0.699)--(0.227,-0.718)--(0.227,-0.738)--(0.227,-0.758)--(0.227,-0.778)--(0.227,-0.799)--(0.227,-0.820)--(0.227,-0.841)--(0.227,-0.862)--(0.227,-0.884)--(0.227,-0.907)--(0.227,-0.930)--(0.227,-0.953)--(0.227,-0.977)--(0.227,-1.002)--(0.227,-1.027)--(0.227,-1.053)--(0.227,-1.080)--(0.227,-1.107)--(0.227,-1.135)--(0.227,-1.164)--(0.227,-1.194)--(0.227,-1.225)--(0.227,-1.257)--(0.227,-1.291)--(0.227,-1.325);
\draw[red!45,line width=1.1pt,dash pattern=on 1pt off 1.6pt] (0.227,-1.325)--(0.218,-1.309)--(0.210,-1.293)--(0.202,-1.278)--(0.194,-1.263)--(0.186,-1.249)--(0.179,-1.235)--(0.172,-1.222)--(0.165,-1.209)--(0.158,-1.196)--(0.152,-1.184)--(0.145,-1.172)--(0.139,-1.160)--(0.133,-1.149)--(0.127,-1.138)--(0.121,-1.127)--(0.116,-1.116)--(0.110,-1.106)--(0.105,-1.096)--(0.099,-1.086)--(0.094,-1.076)--(0.089,-1.066)--(0.084,-1.056)--(0.079,-1.047)--(0.074,-1.038)--(0.069,-1.028)--(0.064,-1.019)--(0.059,-1.010)--(0.055,-1.001)--(0.050,-0.993)--(0.045,-0.984)--(0.041,-0.975)--(0.036,-0.967)--(0.031,-0.958)--(0.027,-0.950)--(0.022,-0.941)--(0.018,-0.933)--(0.013,-0.924)--(0.009,-0.916)--(0.004,-0.907)--(0.000,-0.899);
\draw[red!80!black,line width=1.9pt] (0.227,-0.426)--(0.227,-0.409)--(0.227,-0.391)--(0.227,-0.373)--(0.227,-0.356)--(0.227,-0.338)--(0.227,-0.320)--(0.227,-0.302)--(0.227,-0.284)--(0.227,-0.266)--(0.227,-0.247)--(0.227,-0.229)--(0.227,-0.210)--(0.227,-0.192)--(0.227,-0.173)--(0.227,-0.154)--(0.227,-0.134)--(0.227,-0.115)--(0.227,-0.095)--(0.227,-0.074)--(0.227,-0.054)--(0.227,-0.033)--(0.227,-0.012)--(0.227,0.010)--(0.227,0.032)--(0.227,0.054)--(0.227,0.077)--(0.227,0.101)--(0.227,0.125)--(0.227,0.149)--(0.227,0.174)--(0.227,0.200)--(0.227,0.227)--(0.227,0.254)--(0.227,0.282)--(0.227,0.312)--(0.227,0.342)--(0.227,0.373)--(0.227,0.405)--(0.227,0.438)--(0.227,0.473)--(0.218,0.489)--(0.210,0.505)--(0.202,0.520)--(0.194,0.535)--(0.186,0.549)--(0.179,0.563)--(0.172,0.576)--(0.165,0.589)--(0.158,0.602)--(0.152,0.614)--(0.145,0.626)--(0.139,0.638)--(0.133,0.649)--(0.127,0.660)--(0.121,0.671)--(0.116,0.682)--(0.110,0.692)--(0.105,0.703)--(0.099,0.713)--(0.094,0.722)--(0.089,0.732)--(0.084,0.742)--(0.079,0.751)--(0.074,0.761)--(0.069,0.770)--(0.064,0.779)--(0.059,0.788)--(0.055,0.797)--(0.050,0.806)--(0.045,0.814)--(0.041,0.823)--(0.036,0.832)--(0.031,0.840)--(0.027,0.849)--(0.022,0.857)--(0.018,0.866)--(0.013,0.874)--(0.009,0.882)--(0.004,0.891)--(0.000,0.899);
\draw[red!45,line width=1.1pt,dash pattern=on 1pt off 1.6pt] (-0.227,0.426)--(-0.227,0.409)--(-0.227,0.391)--(-0.227,0.373)--(-0.227,0.356)--(-0.227,0.338)--(-0.227,0.320)--(-0.227,0.302)--(-0.227,0.284)--(-0.227,0.266)--(-0.227,0.247)--(-0.227,0.229)--(-0.227,0.210)--(-0.227,0.192)--(-0.227,0.173)--(-0.227,0.154)--(-0.227,0.134)--(-0.227,0.115)--(-0.227,0.095)--(-0.227,0.074)--(-0.227,0.054)--(-0.227,0.033)--(-0.227,0.012)--(-0.227,-0.010)--(-0.227,-0.032)--(-0.227,-0.054)--(-0.227,-0.077)--(-0.227,-0.101)--(-0.227,-0.125)--(-0.227,-0.149)--(-0.227,-0.174)--(-0.227,-0.200)--(-0.227,-0.227)--(-0.227,-0.254)--(-0.227,-0.282)--(-0.227,-0.312)--(-0.227,-0.342)--(-0.227,-0.373)--(-0.227,-0.405)--(-0.227,-0.438)--(-0.227,-0.473)--(-0.218,-0.489)--(-0.210,-0.505)--(-0.202,-0.520)--(-0.194,-0.535)--(-0.186,-0.549)--(-0.179,-0.563)--(-0.172,-0.576)--(-0.165,-0.589)--(-0.158,-0.602)--(-0.152,-0.614)--(-0.145,-0.626)--(-0.139,-0.638)--(-0.133,-0.649)--(-0.127,-0.660)--(-0.121,-0.671)--(-0.116,-0.682)--(-0.110,-0.692)--(-0.105,-0.703)--(-0.099,-0.713)--(-0.094,-0.722)--(-0.089,-0.732)--(-0.084,-0.742)--(-0.079,-0.751)--(-0.074,-0.761)--(-0.069,-0.770)--(-0.064,-0.779)--(-0.059,-0.788)--(-0.055,-0.797)--(-0.050,-0.806)--(-0.045,-0.814)--(-0.041,-0.823)--(-0.036,-0.832)--(-0.031,-0.840)--(-0.027,-0.849)--(-0.022,-0.857)--(-0.018,-0.866)--(-0.013,-0.874)--(-0.009,-0.882)--(-0.004,-0.891)--(0.000,-0.899);
\draw[red!45,line width=1.1pt,dash pattern=on 1pt off 1.6pt] (-0.227,0.426)--(-0.227,0.444)--(-0.227,0.462)--(-0.227,0.479)--(-0.227,0.497)--(-0.227,0.515)--(-0.227,0.533)--(-0.227,0.551)--(-0.227,0.569)--(-0.227,0.587)--(-0.227,0.605)--(-0.227,0.624)--(-0.227,0.642)--(-0.227,0.661)--(-0.227,0.680)--(-0.227,0.699)--(-0.227,0.718)--(-0.227,0.738)--(-0.227,0.758)--(-0.227,0.778)--(-0.227,0.799)--(-0.227,0.820)--(-0.227,0.841)--(-0.227,0.862)--(-0.227,0.884)--(-0.227,0.907)--(-0.227,0.930)--(-0.227,0.953)--(-0.227,0.977)--(-0.227,1.002)--(-0.227,1.027)--(-0.227,1.053)--(-0.227,1.080)--(-0.227,1.107)--(-0.227,1.135)--(-0.227,1.164)--(-0.227,1.194)--(-0.227,1.225)--(-0.227,1.257)--(-0.227,1.291)--(-0.227,1.325);
\draw[red!80!black,line width=1.9pt] (-0.227,1.325)--(-0.218,1.309)--(-0.210,1.293)--(-0.202,1.278)--(-0.194,1.263)--(-0.186,1.249)--(-0.179,1.235)--(-0.172,1.222)--(-0.165,1.209)--(-0.158,1.196)--(-0.152,1.184)--(-0.145,1.172)--(-0.139,1.160)--(-0.133,1.149)--(-0.127,1.138)--(-0.121,1.127)--(-0.116,1.116)--(-0.110,1.106)--(-0.105,1.096)--(-0.099,1.086)--(-0.094,1.076)--(-0.089,1.066)--(-0.084,1.056)--(-0.079,1.047)--(-0.074,1.038)--(-0.069,1.028)--(-0.064,1.019)--(-0.059,1.010)--(-0.055,1.001)--(-0.050,0.993)--(-0.045,0.984)--(-0.041,0.975)--(-0.036,0.967)--(-0.031,0.958)--(-0.027,0.950)--(-0.022,0.941)--(-0.018,0.933)--(-0.013,0.924)--(-0.009,0.916)--(-0.004,0.907)--(0.000,0.899);
\draw[black,line width=0.9pt] (-0.747,-1.425)--(-1.201,-0.572);
\draw[black,line width=0.9pt] (-1.201,-0.572)--(-1.201,1.226);
\draw[black,line width=0.9pt] (-0.747,0.373)--(-1.201,1.226);
\draw[black,line width=0.9pt] (-0.747,-1.425)--(-0.747,0.373);
\draw[black,line width=0.9pt] (1.201,0.572)--(0.747,1.425);
\draw[black,line width=0.9pt] (1.201,-1.226)--(1.201,0.572);
\draw[black,line width=0.9pt] (-0.747,-1.425)--(1.201,-1.226);
\draw[black,line width=0.9pt] (-0.747,0.373)--(1.201,0.572);
\draw[black,line width=0.9pt] (-1.201,1.226)--(0.747,1.425);
\fill[black] (-0.747,-1.425) circle (1.3pt);
\fill[black] (-0.747,0.373) circle (1.3pt);
\fill[black] (-1.201,-0.572) circle (1.3pt);
\fill[black] (-1.201,1.226) circle (1.3pt);
\fill[black] (1.201,-1.226) circle (1.3pt);
\fill[black] (1.201,0.572) circle (1.3pt);
\fill[black!35] (0.747,-0.373) circle (1.1pt);
\fill[black] (0.747,1.425) circle (1.3pt);
\fill[red!80!black] (-0.974,-0.099) +(-2pt,-2pt) rectangle +(2pt,2pt);
\draw[red!60,fill=white,line width=0.6pt] (0.974,0.099) +(-1.8pt,-1.8pt) rectangle +(1.8pt,1.8pt);
\fill[red!80!black] (0.227,-0.426) +(-2pt,-2pt) rectangle +(2pt,2pt);
\draw[red!60,fill=white,line width=0.6pt] (-0.227,0.426) +(-1.8pt,-1.8pt) rectangle +(1.8pt,1.8pt);
\draw[red!60,fill=white,line width=0.6pt] (0.000,-0.899) +(-1.8pt,-1.8pt) rectangle +(1.8pt,1.8pt);
\fill[red!80!black] (0.000,0.899) +(-2pt,-2pt) rectangle +(2pt,2pt);
\draw[black,line width=0.5pt] (0.000,0.000) ++(-2pt,-2pt)--++(4pt,4pt) (0.000,0.000) ++(-2pt,2pt)--++(4pt,-4pt);
\node[font=\scriptsize,anchor=north,inner sep=2pt] at (0.000,0.000) {$O$};
\draw[red!70!black,densely dotted,line width=0.6pt] (0.000,0.000)--(0.000,0.899);

\node[font=\small] at (0.000,-2.004) {(a) (I) holds};
\end{tikzpicture}\hspace{1.4cm}
\begin{tikzpicture}[scale=0.95,line join=round,line cap=round]
\draw[black!35,dashed] (-0.065,-1.355)--(-0.519,-0.503);
\draw[black!35,dashed] (-0.519,-0.503)--(2.013,1.554);
\draw[black!35,dashed] (-2.467,-0.702)--(-0.519,-0.503);
\fill[black!6] (-2.013,-1.554)--(-2.467,-0.702)--(0.065,1.355)--(0.519,0.503)--cycle;
\fill[black!6] (-2.013,-1.554)--(-0.065,-1.355)--(2.467,0.702)--(0.519,0.503)--cycle;
\fill[black!6] (0.519,0.503)--(2.467,0.702)--(2.013,1.554)--(0.065,1.355)--cycle;
\draw[red!80!black,line width=1.9pt] (-0.362,0.398)--(-0.359,0.392)--(-0.357,0.387)--(-0.354,0.382)--(-0.351,0.377)--(-0.348,0.372)--(-0.346,0.367)--(-0.343,0.362)--(-0.340,0.356)--(-0.337,0.351)--(-0.334,0.346)--(-0.332,0.340)--(-0.329,0.335)--(-0.326,0.330)--(-0.323,0.324)--(-0.320,0.319)--(-0.317,0.313)--(-0.314,0.307)--(-0.311,0.302)--(-0.308,0.296)--(-0.305,0.290)--(-0.301,0.284)--(-0.298,0.278)--(-0.295,0.271)--(-0.291,0.265)--(-0.288,0.259)--(-0.284,0.252)--(-0.281,0.245)--(-0.277,0.238)--(-0.273,0.231)--(-0.269,0.224)--(-0.265,0.216)--(-0.261,0.209)--(-0.257,0.201)--(-0.253,0.193)--(-0.248,0.184)--(-0.244,0.176)--(-0.239,0.167)--(-0.234,0.157)--(-0.229,0.148)--(-0.223,0.138)--(-0.218,0.127)--(-0.212,0.116)--(-0.206,0.105)--(-0.200,0.093)--(-0.193,0.081)--(-0.186,0.068)--(-0.179,0.054)--(-0.171,0.040)--(-0.163,0.025)--(-0.155,0.008)--(-0.146,-0.009)--(-0.136,-0.027)--(-0.121,-0.044)--(-0.105,-0.061)--(-0.090,-0.078)--(-0.075,-0.094)--(-0.061,-0.110)--(-0.046,-0.126)--(-0.032,-0.141)--(-0.019,-0.156)--(-0.005,-0.171)--(0.008,-0.186)--(0.021,-0.200)--(0.034,-0.214)--(0.047,-0.229)--(0.060,-0.242)--(0.072,-0.256)--(0.085,-0.270)--(0.097,-0.283)--(0.109,-0.297)--(0.121,-0.310)--(0.133,-0.323)--(0.145,-0.336)--(0.157,-0.349)--(0.169,-0.362)--(0.180,-0.375)--(0.192,-0.388)--(0.204,-0.401)--(0.216,-0.414)--(0.227,-0.426);
\draw[red!80!black,line width=1.9pt] (-0.362,0.398)--(-0.365,0.403)--(-0.367,0.408)--(-0.370,0.413)--(-0.373,0.418)--(-0.376,0.423)--(-0.378,0.428)--(-0.381,0.433)--(-0.384,0.439)--(-0.387,0.444)--(-0.390,0.449)--(-0.392,0.455)--(-0.395,0.460)--(-0.398,0.465)--(-0.401,0.471)--(-0.404,0.476)--(-0.407,0.482)--(-0.410,0.488)--(-0.413,0.493)--(-0.416,0.499)--(-0.419,0.505)--(-0.423,0.511)--(-0.426,0.517)--(-0.429,0.524)--(-0.433,0.530)--(-0.436,0.536)--(-0.440,0.543)--(-0.443,0.550)--(-0.447,0.557)--(-0.451,0.564)--(-0.455,0.571)--(-0.459,0.579)--(-0.463,0.586)--(-0.467,0.594)--(-0.471,0.602)--(-0.476,0.611)--(-0.480,0.620)--(-0.485,0.629)--(-0.490,0.638)--(-0.495,0.647)--(-0.501,0.657)--(-0.506,0.668)--(-0.512,0.679)--(-0.518,0.690)--(-0.524,0.702)--(-0.531,0.714)--(-0.538,0.727)--(-0.545,0.741)--(-0.553,0.755)--(-0.561,0.770)--(-0.569,0.787)--(-0.578,0.804)--(-0.588,0.822);
\draw[red!45,line width=1.1pt,dash pattern=on 1pt off 1.6pt] (-0.588,0.822)--(-0.575,0.808)--(-0.560,0.791)--(-0.545,0.775)--(-0.530,0.759)--(-0.515,0.743)--(-0.501,0.727)--(-0.487,0.712)--(-0.473,0.696)--(-0.460,0.682)--(-0.446,0.667)--(-0.433,0.652)--(-0.420,0.638)--(-0.407,0.624)--(-0.395,0.610)--(-0.382,0.596)--(-0.370,0.583)--(-0.358,0.569)--(-0.345,0.556)--(-0.333,0.543)--(-0.321,0.530)--(-0.309,0.517)--(-0.298,0.503)--(-0.286,0.491)--(-0.274,0.478)--(-0.262,0.465)--(-0.251,0.452)--(-0.239,0.439)--(-0.227,0.426);
\draw[red!80!black,line width=1.9pt] (-0.362,0.398)--(-0.377,0.386)--(-0.391,0.374)--(-0.406,0.362)--(-0.421,0.350)--(-0.436,0.338)--(-0.451,0.325)--(-0.466,0.313)--(-0.481,0.300)--(-0.497,0.288)--(-0.513,0.275)--(-0.529,0.262)--(-0.546,0.248)--(-0.563,0.234)--(-0.581,0.220)--(-0.599,0.205)--(-0.617,0.190)--(-0.637,0.174)--(-0.657,0.158)--(-0.677,0.141)--(-0.699,0.124)--(-0.722,0.105)--(-0.745,0.086)--(-0.770,0.066)--(-0.797,0.044)--(-0.825,0.022)--(-0.854,-0.002)--(-0.886,-0.028)--(-0.920,-0.056)--(-0.956,-0.085)--(-0.996,-0.117)--(-1.039,-0.152)--(-1.086,-0.191)--(-1.138,-0.233)--(-1.195,-0.279)--(-1.259,-0.331)--(-1.331,-0.390)--(-1.413,-0.456)--(-1.507,-0.533)--(-1.616,-0.621)--(-1.746,-0.727)--(-1.901,-0.853)--(-2.091,-1.007);
\draw[red!45,line width=1.1pt,dash pattern=on 1pt off 1.6pt] (-2.091,-1.007)--(-2.174,-1.121)--(-2.005,-1.104)--(-1.856,-1.089)--(-1.724,-1.075)--(-1.605,-1.063)--(-1.498,-1.052)--(-1.401,-1.042)--(-1.312,-1.033)--(-1.230,-1.025)--(-1.154,-1.017)--(-1.084,-1.010)--(-1.019,-1.003)--(-0.957,-0.997)--(-0.900,-0.991)--(-0.845,-0.985)--(-0.793,-0.980)--(-0.744,-0.975)--(-0.697,-0.970)--(-0.653,-0.966)--(-0.610,-0.961)--(-0.568,-0.957)--(-0.528,-0.953)--(-0.490,-0.949)--(-0.453,-0.945)--(-0.416,-0.942)--(-0.381,-0.938)--(-0.346,-0.934)--(-0.313,-0.931)--(-0.276,-0.916)--(-0.236,-0.884)--(-0.199,-0.854)--(-0.165,-0.826)--(-0.133,-0.800)--(-0.104,-0.776)--(-0.075,-0.753)--(-0.049,-0.731)--(-0.024,-0.711)--(0.000,-0.692);
\draw[red!80!black,line width=1.9pt] (-0.362,0.398)--(-0.358,0.401)--(-0.354,0.404)--(-0.350,0.407)--(-0.347,0.410)--(-0.343,0.413)--(-0.339,0.416)--(-0.335,0.419)--(-0.331,0.423)--(-0.327,0.426)--(-0.323,0.429)--(-0.319,0.432)--(-0.316,0.435)--(-0.312,0.438)--(-0.308,0.442)--(-0.304,0.445)--(-0.300,0.448)--(-0.296,0.451)--(-0.292,0.454)--(-0.288,0.458)--(-0.284,0.461)--(-0.280,0.464)--(-0.276,0.467)--(-0.272,0.470)--(-0.268,0.474)--(-0.264,0.477)--(-0.260,0.480)--(-0.256,0.484)--(-0.252,0.487)--(-0.248,0.490)--(-0.244,0.493)--(-0.240,0.497)--(-0.236,0.500)--(-0.232,0.504)--(-0.227,0.507)--(-0.223,0.510)--(-0.219,0.514)--(-0.215,0.517)--(-0.211,0.521)--(-0.206,0.524)--(-0.202,0.528)--(-0.198,0.531)--(-0.193,0.535)--(-0.189,0.538)--(-0.185,0.542)--(-0.180,0.545)--(-0.176,0.549)--(-0.171,0.553)--(-0.167,0.556)--(-0.162,0.560)--(-0.158,0.564)--(-0.153,0.567)--(-0.148,0.571)--(-0.144,0.575)--(-0.139,0.579)--(-0.134,0.583)--(-0.129,0.586)--(-0.125,0.590)--(-0.120,0.594)--(-0.115,0.598)--(-0.110,0.602)--(-0.105,0.606)--(-0.100,0.610)--(-0.095,0.615)--(-0.090,0.619)--(-0.085,0.623)--(-0.079,0.627)--(-0.074,0.631)--(-0.069,0.636)--(-0.063,0.640)--(-0.058,0.645)--(-0.052,0.649)--(-0.047,0.654)--(-0.041,0.658)--(-0.036,0.663)--(-0.030,0.667)--(-0.024,0.672)--(-0.018,0.677)--(-0.012,0.682)--(-0.006,0.687)--(0.000,0.692);
\draw[red!45,line width=1.1pt,dash pattern=on 1pt off 1.6pt] (0.362,-0.398)--(0.365,-0.403)--(0.367,-0.408)--(0.370,-0.413)--(0.373,-0.418)--(0.376,-0.423)--(0.378,-0.428)--(0.381,-0.433)--(0.384,-0.439)--(0.387,-0.444)--(0.390,-0.449)--(0.392,-0.455)--(0.395,-0.460)--(0.398,-0.465)--(0.401,-0.471)--(0.404,-0.476)--(0.407,-0.482)--(0.410,-0.488)--(0.413,-0.493)--(0.416,-0.499)--(0.419,-0.505)--(0.423,-0.511)--(0.426,-0.517)--(0.429,-0.524)--(0.433,-0.530)--(0.436,-0.536)--(0.440,-0.543)--(0.443,-0.550)--(0.447,-0.557)--(0.451,-0.564)--(0.455,-0.571)--(0.459,-0.579)--(0.463,-0.586)--(0.467,-0.594)--(0.471,-0.602)--(0.476,-0.611)--(0.480,-0.620)--(0.485,-0.629)--(0.490,-0.638)--(0.495,-0.647)--(0.501,-0.657)--(0.506,-0.668)--(0.512,-0.679)--(0.518,-0.690)--(0.524,-0.702)--(0.531,-0.714)--(0.538,-0.727)--(0.545,-0.741)--(0.553,-0.755)--(0.561,-0.770)--(0.569,-0.787)--(0.578,-0.804)--(0.588,-0.822);
\draw[red!80!black,line width=1.9pt] (0.588,-0.822)--(0.575,-0.808)--(0.560,-0.791)--(0.545,-0.775)--(0.530,-0.759)--(0.515,-0.743)--(0.501,-0.727)--(0.487,-0.712)--(0.473,-0.696)--(0.460,-0.682)--(0.446,-0.667)--(0.433,-0.652)--(0.420,-0.638)--(0.407,-0.624)--(0.395,-0.610)--(0.382,-0.596)--(0.370,-0.583)--(0.358,-0.569)--(0.345,-0.556)--(0.333,-0.543)--(0.321,-0.530)--(0.309,-0.517)--(0.298,-0.503)--(0.286,-0.491)--(0.274,-0.478)--(0.262,-0.465)--(0.251,-0.452)--(0.239,-0.439)--(0.227,-0.426);
\draw[red!45,line width=1.1pt,dash pattern=on 1pt off 1.6pt] (0.362,-0.398)--(0.359,-0.392)--(0.357,-0.387)--(0.354,-0.382)--(0.351,-0.377)--(0.348,-0.372)--(0.346,-0.367)--(0.343,-0.362)--(0.340,-0.356)--(0.337,-0.351)--(0.334,-0.346)--(0.332,-0.340)--(0.329,-0.335)--(0.326,-0.330)--(0.323,-0.324)--(0.320,-0.319)--(0.317,-0.313)--(0.314,-0.307)--(0.311,-0.302)--(0.308,-0.296)--(0.305,-0.290)--(0.301,-0.284)--(0.298,-0.278)--(0.295,-0.271)--(0.291,-0.265)--(0.288,-0.259)--(0.284,-0.252)--(0.281,-0.245)--(0.277,-0.238)--(0.273,-0.231)--(0.269,-0.224)--(0.265,-0.216)--(0.261,-0.209)--(0.257,-0.201)--(0.253,-0.193)--(0.248,-0.184)--(0.244,-0.176)--(0.239,-0.167)--(0.234,-0.157)--(0.229,-0.148)--(0.223,-0.138)--(0.218,-0.127)--(0.212,-0.116)--(0.206,-0.105)--(0.200,-0.093)--(0.193,-0.081)--(0.186,-0.068)--(0.179,-0.054)--(0.171,-0.040)--(0.163,-0.025)--(0.155,-0.008)--(0.146,0.009)--(0.136,0.027)--(0.121,0.044)--(0.105,0.061)--(0.090,0.078)--(0.075,0.094)--(0.061,0.110)--(0.046,0.126)--(0.032,0.141)--(0.019,0.156)--(0.005,0.171)--(-0.008,0.186)--(-0.021,0.200)--(-0.034,0.214)--(-0.047,0.229)--(-0.060,0.242)--(-0.072,0.256)--(-0.085,0.270)--(-0.097,0.283)--(-0.109,0.297)--(-0.121,0.310)--(-0.133,0.323)--(-0.145,0.336)--(-0.157,0.349)--(-0.169,0.362)--(-0.180,0.375)--(-0.192,0.388)--(-0.204,0.401)--(-0.216,0.414)--(-0.227,0.426);
\draw[red!45,line width=1.1pt,dash pattern=on 1pt off 1.6pt] (0.362,-0.398)--(0.358,-0.401)--(0.354,-0.404)--(0.350,-0.407)--(0.347,-0.410)--(0.343,-0.413)--(0.339,-0.416)--(0.335,-0.419)--(0.331,-0.423)--(0.327,-0.426)--(0.323,-0.429)--(0.319,-0.432)--(0.316,-0.435)--(0.312,-0.438)--(0.308,-0.442)--(0.304,-0.445)--(0.300,-0.448)--(0.296,-0.451)--(0.292,-0.454)--(0.288,-0.458)--(0.284,-0.461)--(0.280,-0.464)--(0.276,-0.467)--(0.272,-0.470)--(0.268,-0.474)--(0.264,-0.477)--(0.260,-0.480)--(0.256,-0.484)--(0.252,-0.487)--(0.248,-0.490)--(0.244,-0.493)--(0.240,-0.497)--(0.236,-0.500)--(0.232,-0.504)--(0.227,-0.507)--(0.223,-0.510)--(0.219,-0.514)--(0.215,-0.517)--(0.211,-0.521)--(0.206,-0.524)--(0.202,-0.528)--(0.198,-0.531)--(0.193,-0.535)--(0.189,-0.538)--(0.185,-0.542)--(0.180,-0.545)--(0.176,-0.549)--(0.171,-0.553)--(0.167,-0.556)--(0.162,-0.560)--(0.158,-0.564)--(0.153,-0.567)--(0.148,-0.571)--(0.144,-0.575)--(0.139,-0.579)--(0.134,-0.583)--(0.129,-0.586)--(0.125,-0.590)--(0.120,-0.594)--(0.115,-0.598)--(0.110,-0.602)--(0.105,-0.606)--(0.100,-0.610)--(0.095,-0.615)--(0.090,-0.619)--(0.085,-0.623)--(0.079,-0.627)--(0.074,-0.631)--(0.069,-0.636)--(0.063,-0.640)--(0.058,-0.645)--(0.052,-0.649)--(0.047,-0.654)--(0.041,-0.658)--(0.036,-0.663)--(0.030,-0.667)--(0.024,-0.672)--(0.018,-0.677)--(0.012,-0.682)--(0.006,-0.687)--(0.000,-0.692);
\draw[red!45,line width=1.1pt,dash pattern=on 1pt off 1.6pt] (0.362,-0.398)--(0.377,-0.386)--(0.391,-0.374)--(0.406,-0.362)--(0.421,-0.350)--(0.436,-0.338)--(0.451,-0.325)--(0.466,-0.313)--(0.481,-0.300)--(0.497,-0.288)--(0.513,-0.275)--(0.529,-0.262)--(0.546,-0.248)--(0.563,-0.234)--(0.581,-0.220)--(0.599,-0.205)--(0.617,-0.190)--(0.637,-0.174)--(0.657,-0.158)--(0.677,-0.141)--(0.699,-0.124)--(0.722,-0.105)--(0.745,-0.086)--(0.770,-0.066)--(0.797,-0.044)--(0.825,-0.022)--(0.854,0.002)--(0.886,0.028)--(0.920,0.056)--(0.956,0.085)--(0.996,0.117)--(1.039,0.152)--(1.086,0.191)--(1.138,0.233)--(1.195,0.279)--(1.259,0.331)--(1.331,0.390)--(1.413,0.456)--(1.507,0.533)--(1.616,0.621)--(1.746,0.727)--(1.901,0.853)--(2.091,1.007);
\draw[red!80!black,line width=1.9pt] (2.091,1.007)--(2.174,1.121)--(2.005,1.104)--(1.856,1.089)--(1.724,1.075)--(1.605,1.063)--(1.498,1.052)--(1.401,1.042)--(1.312,1.033)--(1.230,1.025)--(1.154,1.017)--(1.084,1.010)--(1.019,1.003)--(0.957,0.997)--(0.900,0.991)--(0.845,0.985)--(0.793,0.980)--(0.744,0.975)--(0.697,0.970)--(0.653,0.966)--(0.610,0.961)--(0.568,0.957)--(0.528,0.953)--(0.490,0.949)--(0.453,0.945)--(0.416,0.942)--(0.381,0.938)--(0.346,0.934)--(0.313,0.931)--(0.276,0.916)--(0.236,0.884)--(0.199,0.854)--(0.165,0.826)--(0.133,0.800)--(0.104,0.776)--(0.075,0.753)--(0.049,0.731)--(0.024,0.711)--(0.000,0.692);
\draw[red!80!black,line width=1.9pt] (0.227,-0.426)--(0.227,-0.444)--(0.227,-0.462)--(0.227,-0.479)--(0.227,-0.497)--(0.227,-0.515)--(0.227,-0.533)--(0.227,-0.551)--(0.227,-0.569)--(0.227,-0.587)--(0.227,-0.605)--(0.227,-0.624)--(0.227,-0.642)--(0.227,-0.661)--(0.227,-0.680)--(0.227,-0.699)--(0.227,-0.718)--(0.227,-0.738)--(0.227,-0.758)--(0.227,-0.778)--(0.227,-0.799)--(0.227,-0.820)--(0.227,-0.841)--(0.227,-0.862)--(0.227,-0.884)--(0.227,-0.907)--(0.227,-0.930)--(0.227,-0.953)--(0.227,-0.977)--(0.227,-1.002)--(0.227,-1.027)--(0.227,-1.053)--(0.227,-1.080)--(0.227,-1.107);
\draw[red!45,line width=1.1pt,dash pattern=on 1pt off 1.6pt] (0.227,-1.107)--(0.222,-1.108)--(0.213,-1.091)--(0.205,-1.076)--(0.197,-1.061)--(0.189,-1.046)--(0.182,-1.033)--(0.175,-1.020)--(0.168,-1.007)--(0.162,-0.995)--(0.155,-0.983)--(0.149,-0.972)--(0.143,-0.961)--(0.138,-0.950)--(0.132,-0.940)--(0.127,-0.930)--(0.122,-0.920)--(0.117,-0.911)--(0.112,-0.902)--(0.107,-0.893)--(0.102,-0.884)--(0.098,-0.875)--(0.093,-0.867)--(0.089,-0.859)--(0.085,-0.851)--(0.081,-0.843)--(0.076,-0.835)--(0.072,-0.827)--(0.068,-0.820)--(0.064,-0.813)--(0.061,-0.805)--(0.057,-0.798)--(0.053,-0.791)--(0.049,-0.784)--(0.046,-0.777)--(0.042,-0.770)--(0.038,-0.764)--(0.035,-0.757)--(0.031,-0.750)--(0.028,-0.744)--(0.024,-0.737)--(0.021,-0.730)--(0.017,-0.724)--(0.014,-0.717)--(0.010,-0.711)--(0.007,-0.704)--(0.003,-0.698)--(0.000,-0.692);
\draw[red!80!black,line width=1.9pt] (0.227,-0.426)--(0.227,-0.409)--(0.227,-0.391)--(0.227,-0.373)--(0.227,-0.356)--(0.227,-0.338)--(0.227,-0.320)--(0.227,-0.302)--(0.227,-0.284)--(0.227,-0.266)--(0.227,-0.247)--(0.227,-0.229)--(0.227,-0.210)--(0.227,-0.192)--(0.227,-0.173)--(0.227,-0.154)--(0.227,-0.134)--(0.227,-0.115)--(0.227,-0.095)--(0.227,-0.074)--(0.227,-0.054)--(0.227,-0.033)--(0.227,-0.012)--(0.227,0.010)--(0.227,0.032)--(0.227,0.054)--(0.227,0.077)--(0.227,0.101)--(0.227,0.125)--(0.227,0.149)--(0.227,0.174)--(0.227,0.200)--(0.227,0.227)--(0.227,0.254)--(0.222,0.276)--(0.213,0.292)--(0.205,0.308)--(0.197,0.323)--(0.189,0.337)--(0.182,0.351)--(0.175,0.364)--(0.168,0.376)--(0.162,0.388)--(0.155,0.400)--(0.149,0.412)--(0.143,0.422)--(0.138,0.433)--(0.132,0.443)--(0.127,0.453)--(0.122,0.463)--(0.117,0.472)--(0.112,0.482)--(0.107,0.491)--(0.102,0.499)--(0.098,0.508)--(0.093,0.516)--(0.089,0.525)--(0.085,0.533)--(0.081,0.540)--(0.076,0.548)--(0.072,0.556)--(0.068,0.563)--(0.064,0.571)--(0.061,0.578)--(0.057,0.585)--(0.053,0.592)--(0.049,0.599)--(0.046,0.606)--(0.042,0.613)--(0.038,0.620)--(0.035,0.626)--(0.031,0.633)--(0.028,0.640)--(0.024,0.646)--(0.021,0.653)--(0.017,0.659)--(0.014,0.666)--(0.010,0.672)--(0.007,0.679)--(0.003,0.685)--(0.000,0.692);
\draw[red!45,line width=1.1pt,dash pattern=on 1pt off 1.6pt] (-0.227,0.426)--(-0.227,0.409)--(-0.227,0.391)--(-0.227,0.373)--(-0.227,0.356)--(-0.227,0.338)--(-0.227,0.320)--(-0.227,0.302)--(-0.227,0.284)--(-0.227,0.266)--(-0.227,0.247)--(-0.227,0.229)--(-0.227,0.210)--(-0.227,0.192)--(-0.227,0.173)--(-0.227,0.154)--(-0.227,0.134)--(-0.227,0.115)--(-0.227,0.095)--(-0.227,0.074)--(-0.227,0.054)--(-0.227,0.033)--(-0.227,0.012)--(-0.227,-0.010)--(-0.227,-0.032)--(-0.227,-0.054)--(-0.227,-0.077)--(-0.227,-0.101)--(-0.227,-0.125)--(-0.227,-0.149)--(-0.227,-0.174)--(-0.227,-0.200)--(-0.227,-0.227)--(-0.227,-0.254)--(-0.222,-0.276)--(-0.213,-0.292)--(-0.205,-0.308)--(-0.197,-0.323)--(-0.189,-0.337)--(-0.182,-0.351)--(-0.175,-0.364)--(-0.168,-0.376)--(-0.162,-0.388)--(-0.155,-0.400)--(-0.149,-0.412)--(-0.143,-0.422)--(-0.138,-0.433)--(-0.132,-0.443)--(-0.127,-0.453)--(-0.122,-0.463)--(-0.117,-0.472)--(-0.112,-0.482)--(-0.107,-0.491)--(-0.102,-0.499)--(-0.098,-0.508)--(-0.093,-0.516)--(-0.089,-0.525)--(-0.085,-0.533)--(-0.081,-0.540)--(-0.076,-0.548)--(-0.072,-0.556)--(-0.068,-0.563)--(-0.064,-0.571)--(-0.061,-0.578)--(-0.057,-0.585)--(-0.053,-0.592)--(-0.049,-0.599)--(-0.046,-0.606)--(-0.042,-0.613)--(-0.038,-0.620)--(-0.035,-0.626)--(-0.031,-0.633)--(-0.028,-0.640)--(-0.024,-0.646)--(-0.021,-0.653)--(-0.017,-0.659)--(-0.014,-0.666)--(-0.010,-0.672)--(-0.007,-0.679)--(-0.003,-0.685)--(0.000,-0.692);
\draw[red!45,line width=1.1pt,dash pattern=on 1pt off 1.6pt] (-0.227,0.426)--(-0.227,0.444)--(-0.227,0.462)--(-0.227,0.479)--(-0.227,0.497)--(-0.227,0.515)--(-0.227,0.533)--(-0.227,0.551)--(-0.227,0.569)--(-0.227,0.587)--(-0.227,0.605)--(-0.227,0.624)--(-0.227,0.642)--(-0.227,0.661)--(-0.227,0.680)--(-0.227,0.699)--(-0.227,0.718)--(-0.227,0.738)--(-0.227,0.758)--(-0.227,0.778)--(-0.227,0.799)--(-0.227,0.820)--(-0.227,0.841)--(-0.227,0.862)--(-0.227,0.884)--(-0.227,0.907)--(-0.227,0.930)--(-0.227,0.953)--(-0.227,0.977)--(-0.227,1.002)--(-0.227,1.027)--(-0.227,1.053)--(-0.227,1.080)--(-0.227,1.107);
\draw[red!80!black,line width=1.9pt] (-0.227,1.107)--(-0.222,1.108)--(-0.213,1.091)--(-0.205,1.076)--(-0.197,1.061)--(-0.189,1.046)--(-0.182,1.033)--(-0.175,1.020)--(-0.168,1.007)--(-0.162,0.995)--(-0.155,0.983)--(-0.149,0.972)--(-0.143,0.961)--(-0.138,0.950)--(-0.132,0.940)--(-0.127,0.930)--(-0.122,0.920)--(-0.117,0.911)--(-0.112,0.902)--(-0.107,0.893)--(-0.102,0.884)--(-0.098,0.875)--(-0.093,0.867)--(-0.089,0.859)--(-0.085,0.851)--(-0.081,0.843)--(-0.076,0.835)--(-0.072,0.827)--(-0.068,0.820)--(-0.064,0.813)--(-0.061,0.805)--(-0.057,0.798)--(-0.053,0.791)--(-0.049,0.784)--(-0.046,0.777)--(-0.042,0.770)--(-0.038,0.764)--(-0.035,0.757)--(-0.031,0.750)--(-0.028,0.744)--(-0.024,0.737)--(-0.021,0.730)--(-0.017,0.724)--(-0.014,0.717)--(-0.010,0.711)--(-0.007,0.704)--(-0.003,0.698)--(0.000,0.692);
\draw[black,line width=0.9pt] (-2.013,-1.554)--(-2.467,-0.702);
\draw[black,line width=0.9pt] (-2.467,-0.702)--(0.065,1.355);
\draw[black,line width=0.9pt] (0.519,0.503)--(0.065,1.355);
\draw[black,line width=0.9pt] (-2.013,-1.554)--(0.519,0.503);
\draw[black,line width=0.9pt] (2.467,0.702)--(2.013,1.554);
\draw[black,line width=0.9pt] (-0.065,-1.355)--(2.467,0.702);
\draw[black,line width=0.9pt] (-2.013,-1.554)--(-0.065,-1.355);
\draw[black,line width=0.9pt] (0.519,0.503)--(2.467,0.702);
\draw[black,line width=0.9pt] (0.065,1.355)--(2.013,1.554);
\fill[black] (-2.013,-1.554) circle (1.3pt);
\fill[black] (0.519,0.503) circle (1.3pt);
\fill[black] (-2.467,-0.702) circle (1.3pt);
\fill[black] (0.065,1.355) circle (1.3pt);
\fill[black] (-0.065,-1.355) circle (1.3pt);
\fill[black] (2.467,0.702) circle (1.3pt);
\fill[black!35] (-0.519,-0.503) circle (1.1pt);
\fill[black] (2.013,1.554) circle (1.3pt);
\fill[red!80!black] (-0.362,0.398) +(-2pt,-2pt) rectangle +(2pt,2pt);
\draw[red!60,fill=white,line width=0.6pt] (0.362,-0.398) +(-1.8pt,-1.8pt) rectangle +(1.8pt,1.8pt);
\fill[red!80!black] (0.227,-0.426) +(-2pt,-2pt) rectangle +(2pt,2pt);
\draw[red!60,fill=white,line width=0.6pt] (-0.227,0.426) +(-1.8pt,-1.8pt) rectangle +(1.8pt,1.8pt);
\draw[red!60,fill=white,line width=0.6pt] (0.000,-0.692) +(-1.8pt,-1.8pt) rectangle +(1.8pt,1.8pt);
\fill[red!80!black] (0.000,0.692) +(-2pt,-2pt) rectangle +(2pt,2pt);
\draw[black,line width=0.5pt] (0.000,0.000) ++(-2pt,-2pt)--++(4pt,4pt) (0.000,0.000) ++(-2pt,2pt)--++(4pt,-4pt);
\node[font=\scriptsize,anchor=north,inner sep=2pt] at (0.000,0.000) {$O$};
\draw[red!70!black,densely dotted,line width=0.6pt] (0.000,0.000)--(0.000,0.692);
\node[font=\small] at (1.616,1.148) {$H$};
\node[font=\small] at (-1.124,-0.449) {$G$};
\draw[->,red!70!black,line width=0.5pt] (-1.200,1.692) -- (-0.070,0.762);
\node[font=\scriptsize,red!70!black,align=center,anchor=south] at (-1.200,1.692) {normal ray of $H$\\ leaves $Q$ through $G$};
\node[font=\small] at (0.000,-2.004) {(b) (I) fails};
\end{tikzpicture}
\caption{The two fans of a $3$-polytope $Q$, with $O\in\inte Q$, traced on its
boundary. The face fan meets $\partial Q$ in the faces of $Q$ (black, thin lines, round
vertices). The normal fan meets $\partial Q$ in a graph (red, thick lines, square
vertices): its vertices are the points where the outer normal rays of the facets leave
$Q$, its edges are the traces of the normal cones $N_Q(e)$ of the edges, and its
regions are the traces of the normal cones $N_Q(v)$ of the vertices; hidden parts are
dashed, respectively dotted, and hidden squares are open. Condition (I) says that
every square vertex lies in the relative interior of its own facet, every thick edge
crosses only its own thin edge, and every region of the thick graph contains its own
vertex. (a) The cube $[-1,1]^3$: (I)
holds. (b) The parallelepiped obtained by the shear
$(x,y,z)\mapsto(x+\tfrac{13}{10}z,\,y,\,z)$, again centred at $O$: the normal ray of
the top facet $H$ (dotted) leaves $Q$ through the side facet $G$. Thus
$\relint N_Q(H)$ meets $\relint\cone(G)$ although $H\not\subseteq G$, and (I) fails;
so does ($\star1$).}
\label{fig:fans}
\end{figure}

The left-hand side says that some linear functional exposes $H$ on $Q$ and, at the
same time, $G^*$ on $Q^\circ$. So (I) says that the common refinement of the normal fan
of $Q$ and its face fan has its cells indexed by the pairs $H\subseteq G$;
equivalently, the faces of the Minkowski sum $Q+Q^\circ$ are the sums $H+G^*$ with
$H\subseteq G$. For $H=Q$ both sides of (I) are false, since $N_Q(Q)=\{O\}$. Figure~\ref{fig:fans} shows how
(I) can be seen on a $3$-polytope: both fans are traced on $\partial Q$, where the
face fan gives the faces of $Q$ and the normal fan gives a graph dual to the edge
graph, and (I) says that this dual graph sits in the expected position.

For a nonempty face $F$ of $Q$ let $q_F$ be the point of $\aff F$ nearest to $O$. The
polytope $Q$ is \emph{perfectly centered} if $q_F\in\relint F$ for every nonempty face $F$
\cite{FukudaWeibel,Dobbins}; Broadie \cite{Broadie} calls this the \emph{projection
condition}, and Lawrence \cite{Lawrence} calls such polytopes \emph{adequately centered}.
Condition (I) is the balance condition $\mathrm{balance}(Q,Q)$ of Dobbins
\cite[Section~3]{Dobbins}, and the equivalence of (a) and (e) below is \cite[Lemma~3.3]{Dobbins},
derived there from Broadie's theorem. We collect the equivalent conditions and their
consequence for $W(\F(Q))$.

\begin{theorem}[Cayley realization]\label{thm:cayley}\label{conj:I}
Let $O\in\inte Q$ and $C=\mathrm{Cay}(Q,Q^\circ)$, and for faces $H\subseteq G$ of $Q$ let
$C_{[H,G]}=\conv\bigl((H\times\{0\})\cup(G^*\times\{1\})\bigr)$, with $\emptyset^*=Q^\circ$ and
$Q^*=\emptyset$. The following are equivalent:
\begin{enumerate}[label=(\alph*)]
\item $Q$ satisfies the interval condition (I);
\item the nonempty proper faces of $Q+Q^\circ$ are exactly the sums $H+G^*$ of faces
      $\emptyset\neq H\subseteq G\neq Q$ of $Q$;
\item the proper faces of $C$ are exactly the sets $C_{[H,G]}$ for faces $H\subseteq G$ of
      $Q$;
\item $C$ is an antiprism with bases $Q\times\{0\}$ and $Q^\circ\times\{1\}$, i.e.\ its facets
      are exactly the sets $C_F=C_{[F,F]}$ for $F\in\F(Q)$;
\item $Q$ is perfectly centered.
\end{enumerate}
If they hold, then $P=(C-z)^\circ$ is a geometric realization of $W(\F(Q))$ for every
$z\in\inte C$.
\end{theorem}
Note that $C_{[Q,Q]}=Q\times\{0\}$ and $C_{[\emptyset,\emptyset]}=Q^\circ\times\{1\}$ are the
two bases, that $C_{[\emptyset,Q]}=\emptyset$, and that $C_{[H,G]}$ meets the two bases in
$H\times\{0\}$ and $G^*\times\{1\}$.
\begin{proof}
(a)$\Leftrightarrow$(b) is the discussion after Definition~\ref{def:intcond}: a nonempty
proper face of $Q+Q^\circ$ is exposed by a linear functional $c\neq0$ and equals $H+G^*$ with
$c\in\relint N_Q(H)$ and $c\in\relint N_{Q^\circ}(G^*)=\relint\cone(G)$, so the sums that
occur are those with $\relint N_Q(H)\cap\relint\cone(G)\neq\emptyset$, which (I) identifies
with the pairs $H\subseteq G$.

(b)$\Leftrightarrow$(c) is the Cayley trick (\S\ref{sec:notation}). The sets $C_{[H,G]}$ with
$H=\emptyset$ or $G=Q$ are the faces of the two bases, which are faces of $C$ in any case.
The remaining proper faces of $C$ are the mixed faces $C_{[H,G]}$ with $H\neq\emptyset$ and
$G\neq Q$ exposed by a common functional, and they correspond bijectively to the nonempty
proper faces $H+G^*$ of $Q+Q^\circ$.

(c)$\Rightarrow$(d): the polytope $C\subset\R^{d+1}$ has dimension $d+1$, so its facets have
dimension $d$. For faces $H\subseteq G$ of $Q$, the direction spaces of $\aff H$ and $\aff G^*$
meet only in $0$: if $y_1,y_2\in G^*$, then $\langle y_1,x\rangle=\langle y_2,x\rangle=1$ for
all $x\in G$ by the definition of the dual face $G^*$, hence $\langle y_1-y_2,x\rangle=0$, so
$y_1-y_2$ is orthogonal to the linear span of $G$, which contains the direction space of
$\aff G$ and hence that of $\aff H$. As $H\times\{0\}$ and $G^*\times\{1\}$ lie in parallel
hyperplanes, $C_{[H,G]}$ therefore has dimension $\dim H+\dim G^*+1=d+\dim H-\dim G$, using
$\dim G^*=d-1-\dim G$ \cite[Section~2]{Broadie}; this also holds if $H=\emptyset$ or $G=Q$,
with $\dim\emptyset=-1$. So among the sets in (c) the facets are exactly those with
$\dim H=\dim G$, i.e.\ with $H=G$.

(d)$\Rightarrow$(c), following the proof of \cite[Lemma~2.1]{Lawrence}. Let $\Phi$ be a
proper face of $C$. Its vertices lie in the two bases, so
$\Phi=\conv\bigl((H\times\{0\})\cup(K\times\{1\})\bigr)$, where $H\times\{0\}=\Phi\cap(Q\times\{0\})$
and $K\times\{1\}=\Phi\cap(Q^\circ\times\{1\})$ are faces of the bases; thus $H$ is a face of
$Q$ and $K=G^*$ for a face $G$ of $Q$, i.e.\ $\Phi=C_{[H,G]}$. Moreover $\Phi$ lies in a
facet $C_F$; intersecting with the two bases gives $H\subseteq F$ and $G^*\subseteq F^*$,
i.e.\ $F\subseteq G$, so $H\subseteq G$. Conversely, let $H\subseteq G$ be faces of $Q$. Then
$C_H\cap C_G$ is a face of $C$; it meets $Q\times\{0\}$ in $(H\cap G)\times\{0\}=H\times\{0\}$
and $Q^\circ\times\{1\}$ in $(H^*\cap G^*)\times\{1\}=G^*\times\{1\}$, since $G^*\subseteq H^*$.
Being the convex hull of its vertices, it equals $C_{[H,G]}$.

(d)$\Leftrightarrow$(e) is Broadie's theorem \cite[Theorem~2.1]{Broadie}.

Finally, under (c) we have $C_{[H,G]}\subseteq C_{[H',G']}$ if and only if $H\subseteq H'$ and
$G'\subseteq G$, i.e.\ $[H',G']\subseteq[H,G]$. So the face lattice of $C$ is the lattice of
intervals of $\F(Q)$ ordered by reverse inclusion, with $C$ itself added; polarity reverses
inclusion once more, so $\F(P)\cong\widetilde I(\F(Q))$, which is the face lattice of
$W(\F(Q))$ by Lemma~\ref{lem:Int} and Remark~\ref{rem:Wdual}.
\end{proof}

\begin{remark}[vertices and facets of the Cayley realization]\label{rem:cayvert}
Under this isomorphism the vertex $w_x$ of $P$ is dual to the facet $C_{[x,x]}$ of
$C$; in particular the poles $w_\emptyset$ and $w_Q$ are dual to the bases
$Q^\circ\times\{1\}$ and $Q\times\{0\}$. The facets of $P$ are dual to the vertices of
$C$: the facet dual to $a\times\{0\}=C_{[a,Q]}$, $a$ a vertex of $Q$, has the vertices
$w_x$ with $a\in x$, and the facet dual to $\Phi^*\times\{1\}=C_{[\emptyset,\Phi]}$,
$\Phi$ a facet of $Q$, has the vertices $w_x$ with $x\subseteq\Phi$.
\end{remark}

\subsection{Known cases}\label{sec:classes}

The metric meaning of perfect centering is developed in the companion paper
\cite[Section~4.4]{Thales}: for perfectly centered $Q$ the realization of
Theorem~\ref{thm:cayley} can be chosen orthogonally bipolar, as the cutout of the apex of a
pyramid over $Q$ by orthogonal projections, and every orthogonally bipolar realization of
$W(\F(Q))$ arises in this way \cite[Theorem~4.10]{Thales}.

\begin{theorem}[known cases]\label{thm:W2}
In each of the following cases $Q$ satisfies (I), so that $W(\F(Q))$ is realized by
Theorem~\ref{thm:cayley}:
\begin{enumerate}[label=(\alph*)]
\item $Q$ is a polygon and the foot of $O$ on the line of every edge lies in that edge,
      for instance if $Q$ is inscribed in a circle with centre $O$
      \cite[Proposition~2]{GrunbaumAnti2};
\item $Q$ is a regular simplex, the cube $[-1,1]^d$ or the cross-polytope
      $\conv\{\pm e_i\}$, and $O$ is its centre \cite[Examples~2.1 and~2.2]{Broadie};
\item $Q$ is a midscribed $3$-polytope and $O$ is the centre of its midsphere
      \cite[Theorem~2]{GrunbaumAnti}, \cite[Theorem~3.4]{Dobbins}.
\end{enumerate}
In particular, every combinatorial type of polytope of dimension at most three has a
realization satisfying (I).
\end{theorem}
\begin{proof}
By Theorem~\ref{thm:cayley} it suffices that $Q$ is perfectly centered about $O$. For a vertex
$F$ always $q_F=F$. In (a) the condition for the edges is the assumption. In (b) the foot
$q_F$ is the centroid of $F$ by symmetry, which lies in $\relint F$. In (c) the foot on an
edge is its point of tangency with the midsphere, and the foot on a facet is the centre of
the circle in which the midsphere meets the plane of the facet; this circle is tangent to
all edges of the facet, so its centre lies in the relative interior of the facet. Every
combinatorial type of polygon is inscribed in a circle, and every combinatorial type of
$3$-polytope has a midscribed realization \cite{Koebe,Andreev,Thurston},
\cite[Theorems~4.12 and~4.13]{Ziegler}.
\end{proof}

For a simplex the realization of $W(\F(Q))$ needs no condition at all: $\F(Q)$ is a
Boolean lattice, whose intervals form the face lattice of a cube, so $W(\F(Q))$ is a cube
with two opposite vertices as poles (Example~\ref{ex:cube}); the apex cutout of
\cite[Proposition~4.9]{Thales} realizes it as a Hjelmslev cuboid, in general not a box
\cite[Figure~3]{Thales}.

\begin{remark}[realizations satisfying (I)]\label{q:realize}
Not every combinatorial type of polytope has a realization satisfying (I). Dobbins
\cite[Theorem~3.9]{Dobbins} constructed a combinatorial $4$-polytope $Q$ without any
antiprism; for it, $W(\F(Q))$ has no geometric realization at all, and in particular no
realization of $Q$ is perfectly centered. In dimension at most three every combinatorial
type has such a realization (Theorems~\ref{thm:W2} and~\ref{thm:W3}), as shown by
Gr\"unbaum \cite[Theorem~2]{GrunbaumAnti} and Dobbins \cite[Theorem~3.4]{Dobbins} via midscribed
polytopes.
\end{remark}

\begin{example}\label{ex:cayley-bip}
For the triangular bipyramid of Example~\ref{ex:thmB} with
$D,E=(0,0,\pm1)$ and $A,B,C$ on the unit circle at the angles
$90^\circ,210^\circ,330^\circ$, one checks by computation that $P$
satisfies (I); the polytope
$(\mathrm{Cay}(P,P^\circ)-z)^\circ$ is a convex $4$-polytope with
$22$ vertices and $11$ facets realizing $W(\F(P))$
(Figure~\ref{fig:thmB}(d)).
\end{example}

\subsection{Configurations and cutouts}

Section~\ref{sec:geom} uses the following notions, also for polytopes that are not
pyramids.

\begin{definition}[configurations and cutouts]\label{def:config}
Let $P\subset\R^n$ be an $n$-polytope. A \emph{configuration} of $P$ is a family
$(x_G)$ of points, one for each nonempty face $G$ of $P$, with $x_G\in\relint G$; thus
$x_v=v$ for every vertex $v$, and $x_P$ is an interior point of $P$. For a vertex $v$
of $P$, the \emph{cutout} of the configuration at $v$ is $\conv\{x_G:v\in G\}$; as only
the points $x_G$ with $v\in G$ enter, we also speak of a \emph{configuration at $v$}
when $x_G\in\relint G$ is given only for the faces $G\ni v$. Given $B\in\inte P$, let
$x_G$ be the orthogonal projection of $B$ onto $\aff G$ (so $x_v=v$ and $x_P=B$). If
$x_G\in\relint G$ for every nonempty face $G$ (resp.\ for every face $G\ni v$), the
family $(x_G)$ is a configuration (resp.\ a configuration at $v$), called the
\emph{orthogonal configuration} of $B$ (at $v$); in general the projections need not
lie in the relative interiors.
\end{definition}

\section{Vertex cutouts in \texorpdfstring{$\R^3$}{R3}: Theorem C}\label{sec:geom}

Throughout this section $P\subset\R^3$ is a $3$-polytope, with vertex set $V(P)$ and
edge set $E(P)$; as $P$ is $3$-dimensional, its $2$-faces are its facets. For a vertex
$v$ let $e_1,\dots,e_m$ be the edges at $v$ and $F_1,\dots,F_m$ the facets at $v$ in
cyclic order, $F_i\supset e_i,e_{i+1}$ (indices modulo $m$).

\subsection{Configurations}

\begin{definition}[regular configuration]\label{def:conf}
Let $(x_G)$ be a family of points indexed by the nonempty faces $G$ of $P$, with $x_v=v$
for every vertex $v$, and put $B=x_P$. It is a \emph{regular configuration} of $P$ if
\begin{enumerate}[label=(C\arabic*)]
\item $B\in\inte P$, $x_e\in\relint e$ for every edge $e$ and $x_F\in\relint F$ for
      every facet $F$, i.e.\ $(x_G)$ is a configuration (Definition~\ref{def:config});
\item for every edge $e=F\cap F'$ the points $B,x_F,x_{F'}$ are not collinear, the
      plane $\Pi_e=\aff\{B,x_F,x_{F'}\}$ does not contain $e$, and $x_e\in\Pi_e$
      (\emph{coplanarity});
\item $x_F\in\inte\conv\{x_e:e\subset F\}$ for every facet $F$;
\item $B\notin\conv W_v$ for every vertex $v$, where $W_v=\{x_G:v\in G\ne P\}$.
\end{enumerate}
\end{definition}

For a facet $F$ and a vertex $v\in F$ with edges $e,e'$ of $F$ at $v$ let
$R_{v,F}=\conv\{v,x_e,x_F,x_{e'}\}\subset F$; for fixed $F$, the quadrilaterals $R_{w,F}$,
$w$ a vertex of $F$, subdivide $F$. Let $\operatorname{star}(v)=\bigcup_{F\ni v}R_{v,F}$ be the
(closed) \emph{star} of $v$ in this subdivision of $\partial P$, i.e.\ the union of its
cells containing $v$, and
\[
C_v=\{B+t(p-B):p\in\operatorname{star}(v),\ 0\le t\le1\}.
\]
The set $C_v$ contains the points $x_G$ with $v\in G$ and is contained in their convex
hull; so if $C_v$ is convex, it is the cutout of the configuration at $v$
(Definition~\ref{def:config}). For regular configurations this is the case by
Theorem~C.

\begin{lemma}\label{lem:quad}
Let $F$ be a facet. Then all $R_{v,F}$, $v\in F$, are convex quadrilaterals if and only
if (C3) holds for $F$; in that case they tile $F$.
\end{lemma}
\begin{proof}
Let $v_1,\dots,v_k$ be the vertices of $F$ and $x_1,\dots,x_k$ the points of the
configuration on its sides, $x_i$ on the side $v_iv_{i+1}$. The polygon
$h=\conv\{x_1,\dots,x_k\}$ is convex, and $x_{i-1}x_i$ is one of its sides. The
quadrilateral $v_i\,x_{i-1}\,x_F\,x_i$ is convex iff its diagonals $v_ix_F$ and
$x_{i-1}x_i$ meet, i.e.\ iff $x_F$ lies in the angle of $F$ at $v_i$ (automatic) and in
the half-plane beyond the line $x_{i-1}x_i$. The intersection over $i$ of these
half-planes is $\inte h$, which means $x_F\in\inte h$. Moreover, if that holds for each
$i$, then $F$ is the union of the corner triangles $v_i\,x_{i-1}\,x_i$ and the triangles
$x_F\,x_{i-1}\,x_i$ forming $h$, and gluing along $x_{i-1}x_i$ gives the tiling.
\end{proof}

\subsection{Convex cutouts: Theorem C}

\begin{theorem}[Theorem C]\label{thm:C}
Let $(x_G)$ be a regular configuration of a $3$-polytope $P$, with $B=x_P$. Then
\begin{enumerate}[label=(\alph*)]
\item each $C_v$ is a geometric bipolar polytope with poles $v$ and $B$, with
      nonpolar vertices $x_{e_i}$, $x_{F_i}$ and $2m$ faces: the quadrilaterals
      $R_{v,F_i}$ (containing $v$) and the planar convex quadrilaterals
      $\conv\{B,x_{F_{i-1}},x_{e_i},x_{F_i}\}$ (containing $B$); combinatorially it is
      the $m$-gonal trapezohedron (Example~\ref{ex:trap});
\item $P=\bigcup_v C_v$, and the $C_v$ have pairwise disjoint interiors.
\end{enumerate}
\end{theorem}
\begin{proof}
Fix $v$ and let $T=v+\cone\{u-v:u\in P\}$ be the tangent cone of $P$ at $v$; it is a
convex polyhedral cone, $P\subseteq T$, $B\in\inte T$, and $F_1,\dots,F_m\subset\partial
T$. Write $x_i=x_{e_i}$, $y_i=x_{F_i}$, $R_i=R_{v,F_i}=\conv\{v,x_i,y_i,x_{i+1}\}$.

\emph{Step 1 (star).} By Lemma~\ref{lem:quad} and (C3) the $R_i$ are convex
quadrilaterals, consecutive ones sharing the segment $[v,x_i]$. So
$\operatorname{star}(v)$ is a closed subset of $\partial T$ whose boundary is the zigzag
polygon $Z=x_1y_1x_2y_2\cdots x_my_m$.

\emph{Step 2 (chart).} By (C4) and the separation theorem there is a linear functional
$\ell$ with $\ell(w)>\ell(B)$ for every $w\in\conv W_v\supseteq\operatorname{star}(v)$.
Let $\Lambda=\{p:\ell(p)=\ell(B)+1\}$ and
$\pi(p)=B+\frac{p-B}{\ell(p)-\ell(B)}$ for $\ell(p)>\ell(B)$. This is the central
projection from $B$: it maps lines to lines, planes through $B$ to lines of $\Lambda$,
convex sets to convex sets, and $\pi(p)=\pi(q)$ iff $p,q$ lie on the same open ray from
$B$.

\emph{Step 3.} A ray from the interior point $B$ meets $\partial T$ at most once, so
$\pi$ is injective on $\operatorname{star}(v)$. The quadrilaterals $R_i$, glued
cyclically around $v$, make $\operatorname{star}(v)$ homeomorphic to a closed
$2$-dimensional ball. Hence $D=\pi(\operatorname{star}(v))$ is a closed topological disc
in the plane $\Lambda$, namely the region bounded by the simple closed polygon $\pi(Z)$,
composed of the convex quadrilaterals $\tilde R_i=\pi(R_i)$.

\emph{Step 4.} By (C2), $y_{i-1},x_i,y_i\in\Pi_{e_i}\ni B$, so
$\pi(y_{i-1}),\pi(x_i),\pi(y_i)$ lie on a line of $\Lambda$. Because $\pi(Z)$ is simple,
the path $\pi(y_{i-1})\to\pi(x_i)\to\pi(y_i)$ cannot fold back on this line:
$\pi(x_i)$ lies strictly between $\pi(y_{i-1})$ and $\pi(y_i)$. So $D$ is the polygon
$y'_1\cdots y'_m$, $y'_i=\pi(y_i)$, with the points $\pi(x_i)$ in the interiors of its
sides. Near $y'_i$ the disc $D$ coincides with $\tilde R_i$, so the interior angle of
$D$ at $y'_i$ equals the angle of the convex quadrilateral $R_i$ there, which is $<\pi$.
A simple polygon with all interior angles $<\pi$ is convex: $D$ is a convex polygon.

\emph{Step 5.} Let $K=B+\R_{\ge0}(D-B)$, the cone over $D$: a closed convex pointed cone
whose boundary consists of the planar sectors $\Sigma_i$ spanned by the rays
$By_{i-1}$, $By_i$, contained in $\Pi_{e_i}$. We claim $C_v=T\cap K$. ``$\subseteq$'':
$[B,p]\subseteq T$ for $p\in\partial T$ and $[B,p]\subseteq K$ for
$p\in\operatorname{star}(v)$. ``$\supseteq$'': if $q\in T\cap K$, $q\ne B$, the ray from
$B$ through $q$ has a direction in $D$, so it meets $\partial T$ in a point
$p\in\operatorname{star}(v)$; as $T$ is convex and $B\in\inte T$, $T\cap\mathrm{ray}=[B,p]$,
so $q\in C_v$. Hence $C_v$ is an intersection of two convex cones, and it is bounded (it
lies in $P$). So $C_v$ is a convex polytope.

\emph{Step 6 (faces).} $C_v\cap\partial T=\operatorname{star}(v)$ gives the $m$ faces
$R_i$. Each $\Pi_{e_i}$ supports $K$, and $C_v\cap\Pi_{e_i}$ is the planar polygon swept
by $B$ over the path $y_{i-1}x_iy_i$, i.e.\ the quadrilateral $B\,y_{i-1}\,x_i\,y_i$, convex
because it is a face of the convex $C_v$; $x_i$ is a genuine vertex since $\Pi_{e_i}$ is
transversal to the edge $e_i$ (C2) and the dihedral angle of $T$ at $e_i$ is $<\pi$.
Altogether $C_v$ has vertices $v,B,x_i,y_i$ and $2m$ quadrilateral faces, hence is the
$m$-gonal trapezohedron; it is bipolar by Example~\ref{ex:trap}. This proves (a).

\emph{Step 7 (tiling).} By Lemma~\ref{lem:quad} the sets $\operatorname{star}(v)$,
$v\in V$, cover $\partial P$ with disjoint interiors. Each $q\in P\setminus\{B\}$ lies on
exactly one segment $[B,p]$ with $p\in\partial P$, so on $C_v$ for the $v$ with
$p\in\operatorname{star}(v)$. This proves (b).
\end{proof}

\begin{remark}\label{rem:C3C4}
Conditions (C3) and (C4) are used essentially: without (C4), the cone $K$ may contain
whole lines and $C_v$ need not be convex; without (C3) the quadrilaterals $R_{v,F}$ are
not convex. For instance, if $F$ is a square facet, then any $x_F\in\relint F\setminus
\inte\conv\{x_e:e\subset F\}$ (allowed by (C1)) makes the quadrilaterals $R_{v,F}$
non-convex by Lemma~\ref{lem:quad}. For (C4) we do not know an explicit configuration
with $B\in\conv W_v$ and non-convex $C_v$, and leave this as a question.
\end{remark}

\subsection{Canonical polytopes}

A $3$-polytope is \emph{canonical} if all its edges touch the unit sphere $\mathbb S^2$
centred at the origin $O$, the touching points $\tau_e$ lie in the relative interiors of
the edges, and $O$ is the barycentre of the points $\tau_e$. We use the following known
facts \cite{Koebe,Andreev,Thurston}, \cite[Theorems~4.12 and~4.13]{Ziegler}.
\begin{enumerate}[label=(K\arabic*)]
\item Every combinatorial type of $3$-polytope has a canonical realization.
\item In a canonical polytope $O\in\inte P$ (it is the barycentre of points of $P$ not
      lying in a plane), and for every facet $F$ the circle $\mathbb S^2\cap\aff F$ is
      inscribed in $F$ and touches the edges of $F$ at $\tau_e$.
\end{enumerate}

\begin{corollary}\label{thm:canon}
Every canonical $3$-polytope has a regular configuration, namely the orthogonal
configuration of its centre $O$ (Definition~\ref{def:config}): $B=O$, $x_e=\tau_e$, and
$x_F$ is the orthogonal projection of $O$ onto $\aff F$. Consequently every combinatorial type
of $3$-polytope has a realization that is the union of convex bipolar polytopes $C_v$ as
in Theorem~\ref{thm:C}.
\end{corollary}
\begin{proof}
(C1): $O\in\inte P$, $\tau_e\in\relint e$; $x_F$ is the centre of the incircle of $F$.

(C2) (coplanarity). Let
$e=F\cap F'$ with direction $d$ and let $\Pi=\{p: d\cdot p=0\}$. Then $B=O\in\Pi$;
$\tau_e\in\Pi$ since $O\tau_e\perp e$; $x_F\in\Pi$ since $Ox_F$ is normal to
$\aff F\supset e$, and likewise $x_{F'}\in\Pi$. The normals of two adjacent facets are
not parallel, so $O,x_F,x_{F'}$ are not collinear and $\Pi=\Pi_e$; it is orthogonal to
$e$, hence does not contain $e$.

(C3) The points $\tau_e$, $e\subset F$, lie on the incircle of $F$ with centre $x_F$.
Two consecutive ones are the tangency points of the two tangents from the vertex $u$
between the corresponding edges, and subtend the central angle $\pi-\theta_u<\pi$, where
$\theta_u\in(0,\pi)$ is the angle of $F$ at $u$. All arcs between consecutive points are
therefore $<\pi$, which means that the centre lies in the interior of their convex hull.

(C4) Fix $v$ and put $\ell(p)=\langle p,v\rangle$; $\ell(B)=0$. Then
$\ell(v)=\|v\|^2>0$; $\ell(\tau_e)=\|\tau_e\|^2=1$ for $e\ni v$, because
$\tau_e,v\in e$ and $O\tau_e\perp e$ give $\langle\tau_e,v-\tau_e\rangle=0$; and
$\ell(x_F)=\|x_F\|^2>0$ for $F\ni v$, because $v\in\aff F$ and $Ox_F\perp\aff F$ give
$\langle x_F,v-x_F\rangle=0$ (and $x_F\ne O$ as $O\in\inte P$). So $\ell>0$ on $W_v$,
whence $B\notin\conv W_v$.

The last statement follows from (K1) and Theorem~\ref{thm:C}.
\end{proof}

\begin{figure}[htbp]
\centering
\begin{tikzpicture}[scale=1.35]
\filldraw[fill=gray!13,draw=black,line width=0.55pt,line join=round] (0.201,-1.171)--(-0.664,-0.708)--(-1.275,-0.852)--(-0.788,-1.404)--(0.124,-1.601)--cycle;
\filldraw[fill=gray!14,draw=black,line width=0.55pt,line join=round] (-1.275,0.025)--(-1.399,0.749)--(-1.476,0.319)--(-1.399,-0.671)--(-1.275,-0.852)--cycle;
\filldraw[fill=gray!27,draw=black,line width=0.55pt,line join=round] (-1.275,0.025)--(-1.275,-0.852)--(-0.664,-0.708)--cycle;
\filldraw[fill=gray!43,draw=black,line width=0.55pt,line join=round] (1.399,0.671)--(0.788,1.404)--(0.000,1.419)--(0.124,0.696)--(0.989,0.233)--cycle;
\filldraw[fill=gray!12,draw=black,line width=0.55pt,line join=round] (0.201,-1.171)--(0.124,-1.601)--(0.989,-1.186)--cycle;
\filldraw[fill=gray!31,draw=black,line width=0.55pt,line join=round] (0.989,0.233)--(0.201,-0.294)--(0.201,-1.171)--(0.989,-1.186)--(1.476,-0.319)--cycle;
\filldraw[fill=gray!32,draw=black,line width=0.55pt,line join=round] (-0.124,1.601)--(-0.989,1.186)--(-1.399,0.749)--(-0.788,0.893)--(0.000,1.419)--cycle;
\filldraw[fill=gray!47,draw=black,line width=0.55pt,line join=round] (0.000,1.419)--(-0.788,0.893)--(0.124,0.696)--cycle;
\filldraw[fill=gray!35,draw=black,line width=0.55pt,line join=round] (-0.788,0.893)--(-1.399,0.749)--(-1.275,0.025)--cycle;
\filldraw[fill=gray!45,draw=black,line width=0.55pt,line join=round] (0.124,0.696)--(-0.788,0.893)--(-1.275,0.025)--(-0.664,-0.708)--(0.201,-0.294)--cycle;
\filldraw[fill=gray!47,draw=black,line width=0.55pt,line join=round] (0.989,0.233)--(0.124,0.696)--(0.201,-0.294)--cycle;
\filldraw[fill=gray!34,draw=black,line width=0.55pt,line join=round] (0.201,-0.294)--(-0.664,-0.708)--(0.201,-1.171)--cycle;
\filldraw[fill=gray!32,draw=black,line width=0.55pt,line join=round] (1.399,0.671)--(0.989,0.233)--(1.476,-0.319)--cycle;
\filldraw[fill=gray!12,draw=black,line width=0.55pt,line join=round] (1.476,-0.319)--(0.989,-1.186)--(1.399,-0.749)--cycle;
\filldraw[fill=gray!33,draw=black,line width=0.55pt,line join=round] (0.788,1.404)--(-0.124,1.601)--(0.000,1.419)--cycle;
\filldraw[fill=gray!23,draw=black,line width=0.55pt,line join=round] (1.275,0.852)--(0.788,1.404)--(1.399,0.671)--cycle;
\node at (0,-2.0) {\small (a)};
\end{tikzpicture}\hspace{1.5cm}
\begin{tikzpicture}[scale=1.35]
\draw[gray!45,line width=0.3pt] (0.989,-1.186)--(0.124,-1.601)--(0.000,-1.419)--(0.788,-0.893)--(1.399,-0.749)--cycle;
\draw[gray!45,line width=0.3pt] (0.000,-1.419)--(-0.788,-1.404)--(-1.399,-0.671)--(-0.989,-0.233)--(-0.124,-0.696)--cycle;
\draw[gray!45,line width=0.3pt] (0.124,-1.601)--(-0.788,-1.404)--(0.000,-1.419)--cycle;
\draw[gray!45,line width=0.3pt] (-1.275,-0.852)--(-1.399,-0.671)--(-0.788,-1.404)--cycle;
\draw[gray!45,line width=0.3pt] (-1.399,-0.671)--(-1.476,0.319)--(-0.989,-0.233)--cycle;
\draw[gray!45,line width=0.3pt] (-0.989,1.186)--(-1.476,0.319)--(-1.399,0.749)--cycle;
\draw[gray!45,line width=0.3pt] (0.788,-0.893)--(-0.124,-0.696)--(-0.201,0.294)--(0.664,0.708)--(1.275,-0.025)--cycle;
\draw[gray!45,line width=0.3pt] (0.000,-1.419)--(-0.124,-0.696)--(0.788,-0.893)--cycle;
\draw[gray!45,line width=0.3pt] (-0.201,0.294)--(-0.989,-0.233)--(-1.476,0.319)--(-0.989,1.186)--(-0.201,1.171)--cycle;
\draw[gray!45,line width=0.3pt] (1.399,-0.749)--(0.788,-0.893)--(1.275,-0.025)--cycle;
\draw[gray!45,line width=0.3pt] (1.476,-0.319)--(1.399,-0.749)--(1.275,-0.025)--(1.275,0.852)--(1.399,0.671)--cycle;
\draw[gray!45,line width=0.3pt] (-0.124,-0.696)--(-0.989,-0.233)--(-0.201,0.294)--cycle;
\draw[gray!45,line width=0.3pt] (0.664,0.708)--(-0.201,0.294)--(-0.201,1.171)--cycle;
\draw[gray!45,line width=0.3pt] (-0.201,1.171)--(-0.989,1.186)--(-0.124,1.601)--cycle;
\draw[gray!45,line width=0.3pt] (1.275,0.852)--(0.664,0.708)--(-0.201,1.171)--(-0.124,1.601)--(0.788,1.404)--cycle;
\draw[gray!45,line width=0.3pt] (1.275,-0.025)--(0.664,0.708)--(1.275,0.852)--cycle;
\draw[red!75!black,line width=0.9pt] (0.000,0.000)--(0.660,-1.170);
\draw[red!75!black,line width=0.9pt] (0.000,0.000)--(-0.660,-0.885);
\draw[red!75!black,line width=0.9pt] (0.000,0.000)--(-0.221,-1.475);
\draw[red!75!black,line width=0.9pt] (0.000,0.000)--(-1.154,-0.976);
\draw[red!75!black,line width=0.9pt] (0.000,0.000)--(-1.288,-0.195);
\draw[red!75!black,line width=0.9pt] (0.000,0.000)--(-1.288,0.751);
\draw[red!75!black,line width=0.9pt] (0.000,0.000)--(0.480,-0.122);
\draw[red!75!black,line width=0.9pt] (0.000,0.000)--(0.221,-1.003);
\draw[red!75!black,line width=0.9pt] (0.000,0.000)--(-0.771,0.547);
\draw[red!75!black,line width=0.9pt] (0.000,0.000)--(1.154,-0.555);
\draw[red!75!black,line width=0.9pt] (0.000,0.000)--(1.365,0.086);
\draw[red!75!black,line width=0.9pt] (0.000,0.000)--(-0.438,-0.212);
\draw[red!75!black,line width=0.9pt] (0.000,0.000)--(0.088,0.724);
\draw[red!75!black,line width=0.9pt] (0.000,0.000)--(-0.438,1.319);
\draw[red!75!black,line width=0.9pt] (0.000,0.000)--(0.480,1.147);
\draw[red!75!black,line width=0.9pt] (0.000,0.000)--(1.071,0.512);
\fill[gray!11,opacity=0.55] (0.201,-1.171)--(-0.664,-0.708)--(-1.275,-0.852)--(-0.788,-1.404)--(0.124,-1.601)--cycle;
\fill[gray!11,opacity=0.55] (-1.275,0.025)--(-1.399,0.749)--(-1.476,0.319)--(-1.399,-0.671)--(-1.275,-0.852)--cycle;
\fill[gray!24,opacity=0.55] (-1.275,0.025)--(-1.275,-0.852)--(-0.664,-0.708)--cycle;
\fill[gray!38,opacity=0.55] (1.399,0.671)--(0.788,1.404)--(0.000,1.419)--(0.124,0.696)--(0.989,0.233)--cycle;
\fill[gray!10,opacity=0.55] (0.201,-1.171)--(0.124,-1.601)--(0.989,-1.186)--cycle;
\fill[gray!27,opacity=0.55] (0.989,0.233)--(0.201,-0.294)--(0.201,-1.171)--(0.989,-1.186)--(1.476,-0.319)--cycle;
\fill[gray!28,opacity=0.55] (-0.124,1.601)--(-0.989,1.186)--(-1.399,0.749)--(-0.788,0.893)--(0.000,1.419)--cycle;
\fill[gray!42,opacity=0.55] (0.000,1.419)--(-0.788,0.893)--(0.124,0.696)--cycle;
\fill[gray!31,opacity=0.55] (-0.788,0.893)--(-1.399,0.749)--(-1.275,0.025)--cycle;
\fill[gray!40,opacity=0.55] (0.124,0.696)--(-0.788,0.893)--(-1.275,0.025)--(-0.664,-0.708)--(0.201,-0.294)--cycle;
\fill[gray!42,opacity=0.55] (0.989,0.233)--(0.124,0.696)--(0.201,-0.294)--cycle;
\fill[gray!30,opacity=0.55] (0.201,-0.294)--(-0.664,-0.708)--(0.201,-1.171)--cycle;
\fill[gray!28,opacity=0.55] (1.399,0.671)--(0.989,0.233)--(1.476,-0.319)--cycle;
\fill[gray!10,opacity=0.55] (1.476,-0.319)--(0.989,-1.186)--(1.399,-0.749)--cycle;
\fill[gray!29,opacity=0.55] (0.788,1.404)--(-0.124,1.601)--(0.000,1.419)--cycle;
\fill[gray!20,opacity=0.55] (1.275,0.852)--(0.788,1.404)--(1.399,0.671)--cycle;
\draw[red!80!black,line width=1.15pt] (0.000,0.000)--(-0.480,-1.147);
\draw[red!80!black,line width=1.15pt] (0.000,0.000)--(-1.365,-0.086);
\draw[red!80!black,line width=1.15pt] (0.000,0.000)--(-1.071,-0.512);
\draw[red!80!black,line width=1.15pt] (0.000,0.000)--(0.660,0.885);
\draw[red!80!black,line width=1.15pt] (0.000,0.000)--(0.438,-1.319);
\draw[red!80!black,line width=1.15pt] (0.000,0.000)--(0.771,-0.547);
\draw[red!80!black,line width=1.15pt] (0.000,0.000)--(-0.660,1.170);
\draw[red!80!black,line width=1.15pt] (0.000,0.000)--(-0.221,1.003);
\draw[red!80!black,line width=1.15pt] (0.000,0.000)--(-1.154,0.555);
\draw[red!80!black,line width=1.15pt] (0.000,0.000)--(-0.480,0.122);
\draw[red!80!black,line width=1.15pt] (0.000,0.000)--(0.438,0.212);
\draw[red!80!black,line width=1.15pt] (0.000,0.000)--(-0.088,-0.724);
\draw[red!80!black,line width=1.15pt] (0.000,0.000)--(1.288,0.195);
\draw[red!80!black,line width=1.15pt] (0.000,0.000)--(1.288,-0.751);
\draw[red!80!black,line width=1.15pt] (0.000,0.000)--(0.221,1.475);
\draw[red!80!black,line width=1.15pt] (0.000,0.000)--(1.154,0.976);
\draw[blue!65!black,line width=0.4pt,dash pattern=on 1.3pt off 0.9pt] (-0.480,-1.147)--(-0.232,-0.940);
\draw[blue!65!black,line width=0.4pt,dash pattern=on 1.3pt off 0.9pt] (-0.480,-1.147)--(-0.969,-0.780);
\draw[blue!65!black,line width=0.4pt,dash pattern=on 1.3pt off 0.9pt] (-0.480,-1.147)--(-1.031,-1.128);
\draw[blue!65!black,line width=0.4pt,dash pattern=on 1.3pt off 0.9pt] (-0.480,-1.147)--(-0.332,-1.502);
\draw[blue!65!black,line width=0.4pt,dash pattern=on 1.3pt off 0.9pt] (-0.480,-1.147)--(0.162,-1.386);
\draw[black,line width=0.55pt] (0.201,-1.171)--(-0.664,-0.708)--(-1.275,-0.852)--(-0.788,-1.404)--(0.124,-1.601)--cycle;
\draw[blue!65!black,line width=0.4pt,dash pattern=on 1.3pt off 0.9pt] (-1.365,-0.086)--(-1.337,0.387);
\draw[blue!65!black,line width=0.4pt,dash pattern=on 1.3pt off 0.9pt] (-1.365,-0.086)--(-1.437,0.534);
\draw[blue!65!black,line width=0.4pt,dash pattern=on 1.3pt off 0.9pt] (-1.365,-0.086)--(-1.437,-0.176);
\draw[blue!65!black,line width=0.4pt,dash pattern=on 1.3pt off 0.9pt] (-1.365,-0.086)--(-1.337,-0.762);
\draw[blue!65!black,line width=0.4pt,dash pattern=on 1.3pt off 0.9pt] (-1.365,-0.086)--(-1.275,-0.414);
\draw[black,line width=0.55pt] (-1.275,0.025)--(-1.399,0.749)--(-1.476,0.319)--(-1.399,-0.671)--(-1.275,-0.852)--cycle;
\draw[blue!65!black,line width=0.4pt,dash pattern=on 1.3pt off 0.9pt] (-1.071,-0.512)--(-1.275,-0.414);
\draw[blue!65!black,line width=0.4pt,dash pattern=on 1.3pt off 0.9pt] (-1.071,-0.512)--(-0.969,-0.780);
\draw[blue!65!black,line width=0.4pt,dash pattern=on 1.3pt off 0.9pt] (-1.071,-0.512)--(-0.969,-0.342);
\draw[black,line width=0.55pt] (-1.275,0.025)--(-1.275,-0.852)--(-0.664,-0.708)--cycle;
\draw[blue!65!black,line width=0.4pt,dash pattern=on 1.3pt off 0.9pt] (0.660,0.885)--(1.093,1.037);
\draw[blue!65!black,line width=0.4pt,dash pattern=on 1.3pt off 0.9pt] (0.660,0.885)--(0.394,1.412);
\draw[blue!65!black,line width=0.4pt,dash pattern=on 1.3pt off 0.9pt] (0.660,0.885)--(0.062,1.058);
\draw[blue!65!black,line width=0.4pt,dash pattern=on 1.3pt off 0.9pt] (0.660,0.885)--(0.556,0.464);
\draw[blue!65!black,line width=0.4pt,dash pattern=on 1.3pt off 0.9pt] (0.660,0.885)--(1.194,0.452);
\draw[black,line width=0.55pt] (1.399,0.671)--(0.788,1.404)--(0.000,1.419)--(0.124,0.696)--(0.989,0.233)--cycle;
\draw[blue!65!black,line width=0.4pt,dash pattern=on 1.3pt off 0.9pt] (0.438,-1.319)--(0.162,-1.386);
\draw[blue!65!black,line width=0.4pt,dash pattern=on 1.3pt off 0.9pt] (0.438,-1.319)--(0.556,-1.394);
\draw[blue!65!black,line width=0.4pt,dash pattern=on 1.3pt off 0.9pt] (0.438,-1.319)--(0.595,-1.179);
\draw[black,line width=0.55pt] (0.201,-1.171)--(0.124,-1.601)--(0.989,-1.186)--cycle;
\draw[blue!65!black,line width=0.4pt,dash pattern=on 1.3pt off 0.9pt] (0.771,-0.547)--(0.595,-0.030);
\draw[blue!65!black,line width=0.4pt,dash pattern=on 1.3pt off 0.9pt] (0.771,-0.547)--(0.201,-0.732);
\draw[blue!65!black,line width=0.4pt,dash pattern=on 1.3pt off 0.9pt] (0.771,-0.547)--(0.595,-1.179);
\draw[blue!65!black,line width=0.4pt,dash pattern=on 1.3pt off 0.9pt] (0.771,-0.547)--(1.232,-0.753);
\draw[blue!65!black,line width=0.4pt,dash pattern=on 1.3pt off 0.9pt] (0.771,-0.547)--(1.232,-0.043);
\draw[black,line width=0.55pt] (0.989,0.233)--(0.201,-0.294)--(0.201,-1.171)--(0.989,-1.186)--(1.476,-0.319)--cycle;
\draw[blue!65!black,line width=0.4pt,dash pattern=on 1.3pt off 0.9pt] (-0.660,1.170)--(-0.556,1.394);
\draw[blue!65!black,line width=0.4pt,dash pattern=on 1.3pt off 0.9pt] (-0.660,1.170)--(-1.194,0.968);
\draw[blue!65!black,line width=0.4pt,dash pattern=on 1.3pt off 0.9pt] (-0.660,1.170)--(-1.093,0.821);
\draw[blue!65!black,line width=0.4pt,dash pattern=on 1.3pt off 0.9pt] (-0.660,1.170)--(-0.394,1.156);
\draw[blue!65!black,line width=0.4pt,dash pattern=on 1.3pt off 0.9pt] (-0.660,1.170)--(-0.062,1.510);
\draw[black,line width=0.55pt] (-0.124,1.601)--(-0.989,1.186)--(-1.399,0.749)--(-0.788,0.893)--(0.000,1.419)--cycle;
\draw[blue!65!black,line width=0.4pt,dash pattern=on 1.3pt off 0.9pt] (-0.221,1.003)--(-0.394,1.156);
\draw[blue!65!black,line width=0.4pt,dash pattern=on 1.3pt off 0.9pt] (-0.221,1.003)--(-0.332,0.794);
\draw[blue!65!black,line width=0.4pt,dash pattern=on 1.3pt off 0.9pt] (-0.221,1.003)--(0.062,1.058);
\draw[black,line width=0.55pt] (0.000,1.419)--(-0.788,0.893)--(0.124,0.696)--cycle;
\draw[blue!65!black,line width=0.4pt,dash pattern=on 1.3pt off 0.9pt] (-1.154,0.555)--(-1.093,0.821);
\draw[blue!65!black,line width=0.4pt,dash pattern=on 1.3pt off 0.9pt] (-1.154,0.555)--(-1.337,0.387);
\draw[blue!65!black,line width=0.4pt,dash pattern=on 1.3pt off 0.9pt] (-1.154,0.555)--(-1.031,0.459);
\draw[black,line width=0.55pt] (-0.788,0.893)--(-1.399,0.749)--(-1.275,0.025)--cycle;
\draw[blue!65!black,line width=0.4pt,dash pattern=on 1.3pt off 0.9pt] (-0.480,0.122)--(-0.332,0.794);
\draw[blue!65!black,line width=0.4pt,dash pattern=on 1.3pt off 0.9pt] (-0.480,0.122)--(-1.031,0.459);
\draw[blue!65!black,line width=0.4pt,dash pattern=on 1.3pt off 0.9pt] (-0.480,0.122)--(-0.969,-0.342);
\draw[blue!65!black,line width=0.4pt,dash pattern=on 1.3pt off 0.9pt] (-0.480,0.122)--(-0.232,-0.501);
\draw[blue!65!black,line width=0.4pt,dash pattern=on 1.3pt off 0.9pt] (-0.480,0.122)--(0.162,0.201);
\draw[black,line width=0.55pt] (0.124,0.696)--(-0.788,0.893)--(-1.275,0.025)--(-0.664,-0.708)--(0.201,-0.294)--cycle;
\draw[blue!65!black,line width=0.4pt,dash pattern=on 1.3pt off 0.9pt] (0.438,0.212)--(0.556,0.464);
\draw[blue!65!black,line width=0.4pt,dash pattern=on 1.3pt off 0.9pt] (0.438,0.212)--(0.162,0.201);
\draw[blue!65!black,line width=0.4pt,dash pattern=on 1.3pt off 0.9pt] (0.438,0.212)--(0.595,-0.030);
\draw[black,line width=0.55pt] (0.989,0.233)--(0.124,0.696)--(0.201,-0.294)--cycle;
\draw[blue!65!black,line width=0.4pt,dash pattern=on 1.3pt off 0.9pt] (-0.088,-0.724)--(-0.232,-0.501);
\draw[blue!65!black,line width=0.4pt,dash pattern=on 1.3pt off 0.9pt] (-0.088,-0.724)--(-0.232,-0.940);
\draw[blue!65!black,line width=0.4pt,dash pattern=on 1.3pt off 0.9pt] (-0.088,-0.724)--(0.201,-0.732);
\draw[black,line width=0.55pt] (0.201,-0.294)--(-0.664,-0.708)--(0.201,-1.171)--cycle;
\draw[blue!65!black,line width=0.4pt,dash pattern=on 1.3pt off 0.9pt] (1.288,0.195)--(1.194,0.452);
\draw[blue!65!black,line width=0.4pt,dash pattern=on 1.3pt off 0.9pt] (1.288,0.195)--(1.232,-0.043);
\draw[blue!65!black,line width=0.4pt,dash pattern=on 1.3pt off 0.9pt] (1.288,0.195)--(1.437,0.176);
\draw[black,line width=0.55pt] (1.399,0.671)--(0.989,0.233)--(1.476,-0.319)--cycle;
\draw[blue!65!black,line width=0.4pt,dash pattern=on 1.3pt off 0.9pt] (1.288,-0.751)--(1.232,-0.753);
\draw[blue!65!black,line width=0.4pt,dash pattern=on 1.3pt off 0.9pt] (1.288,-0.751)--(1.194,-0.968);
\draw[blue!65!black,line width=0.4pt,dash pattern=on 1.3pt off 0.9pt] (1.288,-0.751)--(1.437,-0.534);
\draw[black,line width=0.55pt] (1.476,-0.319)--(0.989,-1.186)--(1.399,-0.749)--cycle;
\draw[blue!65!black,line width=0.4pt,dash pattern=on 1.3pt off 0.9pt] (0.221,1.475)--(0.332,1.502);
\draw[blue!65!black,line width=0.4pt,dash pattern=on 1.3pt off 0.9pt] (0.221,1.475)--(-0.062,1.510);
\draw[blue!65!black,line width=0.4pt,dash pattern=on 1.3pt off 0.9pt] (0.221,1.475)--(0.394,1.412);
\draw[black,line width=0.55pt] (0.788,1.404)--(-0.124,1.601)--(0.000,1.419)--cycle;
\draw[blue!65!black,line width=0.4pt,dash pattern=on 1.3pt off 0.9pt] (1.154,0.976)--(1.031,1.128);
\draw[blue!65!black,line width=0.4pt,dash pattern=on 1.3pt off 0.9pt] (1.154,0.976)--(1.093,1.037);
\draw[blue!65!black,line width=0.4pt,dash pattern=on 1.3pt off 0.9pt] (1.154,0.976)--(1.337,0.762);
\draw[black,line width=0.55pt] (1.275,0.852)--(0.788,1.404)--(1.399,0.671)--cycle;
\filldraw[fill=white,draw=red!80!black,line width=0.6pt] (0.000,0.000) circle (1.3pt);
\node[red!70!black,font=\small,below right=1pt] at (0.000,0.000) {$O$};
\node at (0,-2.0) {\small (b)};
\end{tikzpicture}\\[0.4cm]
\begin{tikzpicture}[scale=1.35]
\filldraw[fill=blue!33,draw=black!75,line width=0.3pt,line join=round] (0.000,-0.000)--(0.480,-0.122)--(0.232,0.501)--(0.088,0.724)--cycle;
\filldraw[fill=blue!15,draw=black!75,line width=0.3pt,line join=round] (-0.438,-0.212)--(-0.162,-0.201)--(0.480,-0.122)--(0.000,-0.000)--cycle;
\filldraw[fill=blue!47,draw=black!75,line width=0.3pt,line join=round] (-0.771,0.547)--(0.000,-0.000)--(0.088,0.724)--(-0.201,0.732)--cycle;
\filldraw[fill=blue!17,draw=black!75,line width=0.3pt,line join=round] (-0.595,0.030)--(-0.438,-0.212)--(0.000,-0.000)--(-0.771,0.547)--cycle;
\filldraw[fill=orange!28,draw=black!75,line width=0.3pt,line join=round] (-0.438,-0.212)--(-0.556,-0.464)--(-0.660,-0.885)--(0.000,0.000)--cycle;
\filldraw[fill=orange!27,draw=black!75,line width=0.3pt,line join=round] (0.000,0.000)--(-0.660,-0.885)--(-0.062,-1.058)--(0.221,-1.003)--cycle;
\filldraw[fill=orange!25,draw=black!75,line width=0.3pt,line join=round] (0.480,-0.122)--(0.000,0.000)--(0.221,-1.003)--(0.332,-0.794)--cycle;
\filldraw[fill=green!70!black!15,draw=black!75,line width=0.3pt,line join=round] (0.480,-0.122)--(0.969,0.342)--(1.071,0.512)--(0.000,0.000)--cycle;
\filldraw[fill=green!70!black!28,draw=black!75,line width=0.3pt,line join=round] (0.232,0.940)--(0.088,0.724)--(0.000,0.000)--(0.480,1.147)--cycle;
\filldraw[fill=green!70!black!46,draw=black!75,line width=0.3pt,line join=round] (0.480,1.147)--(0.000,0.000)--(1.071,0.512)--(0.969,0.780)--cycle;
\filldraw[fill=green!70!black!48,draw=black!75,line width=0.3pt,line join=round] (-0.771,0.547)--(-1.232,0.043)--(-1.288,-0.195)--(0.000,-0.000)--cycle;
\filldraw[fill=green!70!black!36,draw=black!75,line width=0.3pt,line join=round] (0.000,-0.000)--(-1.288,-0.195)--(-1.194,-0.452)--(-0.660,-0.885)--cycle;
\filldraw[fill=orange!54,draw=black!75,line width=0.3pt,line join=round] (0.480,1.147)--(-0.162,1.386)--(-0.438,1.319)--(0.000,-0.000)--cycle;
\filldraw[fill=orange!36,draw=black!75,line width=0.3pt,line join=round] (0.000,-0.000)--(-0.438,1.319)--(-0.595,1.179)--(-0.771,0.547)--cycle;
\filldraw[fill=green!70!black!46,draw=black!75,line width=0.3pt,line join=round] (-0.000,0.000)--(0.660,-1.170)--(1.093,-0.821)--(1.154,-0.555)--cycle;
\filldraw[fill=green!70!black!17,draw=black!75,line width=0.3pt,line join=round] (0.221,-1.003)--(0.394,-1.156)--(0.660,-1.170)--(-0.000,0.000)--cycle;
\filldraw[fill=green!70!black!33,draw=black!75,line width=0.3pt,line join=round] (0.480,-0.122)--(-0.000,0.000)--(1.154,-0.555)--(1.031,-0.459)--cycle;
\filldraw[fill=blue!47,draw=black!75,line width=0.3pt,line join=round] (1.365,0.086)--(1.275,0.414)--(1.071,0.512)--(-0.000,-0.000)--cycle;
\filldraw[fill=blue!27,draw=black!75,line width=0.3pt,line join=round] (1.365,0.086)--(-0.000,-0.000)--(1.154,-0.555)--(1.337,-0.387)--cycle;
\filldraw[fill=blue!36,draw=black!75,line width=0.3pt,line join=round] (-0.394,-1.412)--(-0.221,-1.475)--(0.000,0.000)--(-0.660,-0.885)--cycle;
\filldraw[fill=blue!47,draw=black!75,line width=0.3pt,line join=round] (-0.221,-1.475)--(0.062,-1.510)--(0.660,-1.170)--(0.000,0.000)--cycle;
\filldraw[fill=blue!15,draw=black!75,line width=0.3pt,line join=round] (0.000,0.000)--(-1.365,-0.086)--(-1.437,-0.176)--(-1.288,-0.195)--cycle;
\filldraw[fill=blue!17,draw=black!75,line width=0.3pt,line join=round] (-1.437,-0.176)--(-1.365,-0.086)--(-1.437,0.534)--(-1.476,0.319)--cycle;
\filldraw[fill=blue!54,draw=black!75,line width=0.3pt,line join=round] (-1.288,0.751)--(-1.437,0.534)--(-1.365,-0.086)--(0.000,0.000)--cycle;
\filldraw[fill=blue!33,draw=black!75,line width=0.3pt,line join=round] (-0.771,0.547)--(-1.232,0.753)--(-1.288,0.751)--(0.000,0.000)--cycle;
\filldraw[fill=orange!48,draw=black!75,line width=0.3pt,line join=round] (1.154,0.976)--(1.031,1.128)--(0.480,1.147)--(-0.000,-0.000)--cycle;
\filldraw[fill=orange!47,draw=black!75,line width=0.3pt,line join=round] (1.365,0.086)--(1.337,0.762)--(1.154,0.976)--(-0.000,-0.000)--cycle;
\filldraw[fill=orange!27,draw=black!75,line width=0.3pt,line join=round] (1.031,1.128)--(1.154,0.976)--(1.337,0.762)--(1.275,0.852)--cycle;
\filldraw[fill=green!70!black!36,draw=black!75,line width=0.3pt,line join=round] (-0.660,1.170)--(-1.194,0.968)--(-1.288,0.751)--(0.000,-0.000)--cycle;
\filldraw[fill=green!70!black!36,draw=black!75,line width=0.3pt,line join=round] (-0.989,1.186)--(-1.194,0.968)--(-0.660,1.170)--(-0.556,1.394)--cycle;
\filldraw[fill=green!70!black!33,draw=black!75,line width=0.3pt,line join=round] (-0.660,1.170)--(0.000,-0.000)--(-0.438,1.319)--(-0.556,1.394)--cycle;
\filldraw[fill=orange!17,draw=black!75,line width=0.3pt,line join=round] (-1.437,-0.176)--(-1.399,-0.671)--(-1.337,-0.762)--(-1.365,-0.086)--cycle;
\filldraw[fill=orange!47,draw=black!75,line width=0.3pt,line join=round] (-1.365,-0.086)--(-1.337,-0.762)--(-1.154,-0.976)--(0.000,-0.000)--cycle;
\filldraw[fill=orange!15,draw=black!75,line width=0.3pt,line join=round] (-0.660,-0.885)--(0.000,-0.000)--(-1.154,-0.976)--(-1.093,-1.037)--cycle;
\filldraw[fill=orange!54,draw=black!75,line width=0.3pt,line join=round] (-0.000,0.000)--(1.288,-0.751)--(1.437,-0.534)--(1.365,0.086)--cycle;
\filldraw[fill=orange!36,draw=black!75,line width=0.3pt,line join=round] (0.660,-1.170)--(1.194,-0.968)--(1.288,-0.751)--(-0.000,0.000)--cycle;
\filldraw[fill=orange!15,draw=black!75,line width=0.3pt,line join=round] (1.437,-0.534)--(1.288,-0.751)--(1.194,-0.968)--(1.399,-0.749)--cycle;
\filldraw[fill=blue!36,draw=black!75,line width=0.3pt,line join=round] (-0.556,1.394)--(-0.660,1.170)--(-0.062,1.510)--(-0.124,1.601)--cycle;
\filldraw[fill=blue!47,draw=black!75,line width=0.3pt,line join=round] (0.221,1.475)--(-0.062,1.510)--(-0.660,1.170)--(0.000,-0.000)--cycle;
\filldraw[fill=blue!25,draw=black!75,line width=0.3pt,line join=round] (0.480,1.147)--(0.332,1.502)--(0.221,1.475)--(0.000,-0.000)--cycle;
\filldraw[fill=blue!37,draw=black!75,line width=0.3pt,line join=round] (-0.062,1.510)--(0.221,1.475)--(0.332,1.502)--(-0.124,1.601)--cycle;
\filldraw[fill=green!70!black!16,draw=black!75,line width=0.3pt,line join=round] (-0.480,-1.147)--(-1.031,-1.128)--(-0.788,-1.404)--(-0.332,-1.502)--cycle;
\filldraw[fill=green!70!black!48,draw=black!75,line width=0.3pt,line join=round] (0.000,0.000)--(-1.154,-0.976)--(-1.031,-1.128)--(-0.480,-1.147)--cycle;
\filldraw[fill=green!70!black!25,draw=black!75,line width=0.3pt,line join=round] (-0.480,-1.147)--(-0.332,-1.502)--(-0.221,-1.475)--(0.000,0.000)--cycle;
\filldraw[fill=white,draw=black,line width=0.5pt] (0.000,0.000) circle (1.2pt);
\node[font=\small,above right=1pt] at (0.000,0.000) {$O$};
\node at (0,-2.11) {\small (c)};
\end{tikzpicture}\hspace{1.5cm}
\begin{tikzpicture}[scale=0.85]
\filldraw[fill=blue!33,draw=black!75,line width=0.3pt,line join=round] (-0.160,0.235)--(0.320,0.113)--(0.071,0.736)--(-0.073,0.959)--cycle;
\filldraw[fill=blue!15,draw=black!75,line width=0.3pt,line join=round] (-0.598,0.023)--(-0.323,0.034)--(0.320,0.113)--(-0.160,0.235)--cycle;
\filldraw[fill=blue!47,draw=black!75,line width=0.3pt,line join=round] (-0.931,0.782)--(-0.160,0.235)--(-0.073,0.959)--(-0.361,0.967)--cycle;
\filldraw[fill=blue!17,draw=black!75,line width=0.3pt,line join=round] (-0.755,0.265)--(-0.598,0.023)--(-0.160,0.235)--(-0.931,0.782)--cycle;
\filldraw[fill=orange!28,draw=black!75,line width=0.3pt,line join=round] (-0.537,-0.768)--(-0.655,-1.021)--(-0.759,-1.441)--(-0.099,-0.557)--cycle;
\filldraw[fill=orange!27,draw=black!75,line width=0.3pt,line join=round] (-0.099,-0.557)--(-0.759,-1.441)--(-0.161,-1.614)--(0.122,-1.559)--cycle;
\filldraw[fill=orange!25,draw=black!75,line width=0.3pt,line join=round] (0.381,-0.679)--(-0.099,-0.557)--(0.122,-1.559)--(0.233,-1.351)--cycle;
\filldraw[fill=green!70!black!15,draw=black!75,line width=0.3pt,line join=round] (1.012,0.444)--(1.501,0.908)--(1.602,1.078)--(0.531,0.567)--cycle;
\filldraw[fill=green!70!black!28,draw=black!75,line width=0.3pt,line join=round] (0.763,1.506)--(0.619,1.291)--(0.531,0.567)--(1.012,1.714)--cycle;
\filldraw[fill=green!70!black!46,draw=black!75,line width=0.3pt,line join=round] (1.012,1.714)--(0.531,0.567)--(1.602,1.078)--(1.501,1.347)--cycle;
\filldraw[fill=green!70!black!48,draw=black!75,line width=0.3pt,line join=round] (-1.562,0.361)--(-2.023,-0.144)--(-2.078,-0.381)--(-0.791,-0.186)--cycle;
\filldraw[fill=green!70!black!36,draw=black!75,line width=0.3pt,line join=round] (-0.791,-0.186)--(-2.078,-0.381)--(-1.985,-0.638)--(-1.451,-1.071)--cycle;
\filldraw[fill=orange!54,draw=black!75,line width=0.3pt,line join=round] (0.320,2.084)--(-0.323,2.323)--(-0.598,2.256)--(-0.160,0.937)--cycle;
\filldraw[fill=orange!36,draw=black!75,line width=0.3pt,line join=round] (-0.160,0.937)--(-0.598,2.256)--(-0.755,2.115)--(-0.931,1.484)--cycle;
\filldraw[fill=green!70!black!46,draw=black!75,line width=0.3pt,line join=round] (0.630,-0.714)--(1.290,-1.884)--(1.724,-1.535)--(1.784,-1.270)--cycle;
\filldraw[fill=green!70!black!17,draw=black!75,line width=0.3pt,line join=round] (0.852,-1.717)--(1.024,-1.870)--(1.290,-1.884)--(0.630,-0.714)--cycle;
\filldraw[fill=green!70!black!33,draw=black!75,line width=0.3pt,line join=round] (1.111,-0.836)--(0.630,-0.714)--(1.784,-1.270)--(1.662,-1.173)--cycle;
\filldraw[fill=blue!47,draw=black!75,line width=0.3pt,line join=round] (2.385,0.066)--(2.295,0.394)--(2.091,0.492)--(1.020,-0.020)--cycle;
\filldraw[fill=blue!27,draw=black!75,line width=0.3pt,line join=round] (2.385,0.066)--(1.020,-0.020)--(2.174,-0.575)--(2.357,-0.407)--cycle;
\filldraw[fill=blue!36,draw=black!75,line width=0.3pt,line join=round] (-0.394,-2.547)--(-0.221,-2.610)--(0.000,-1.136)--(-0.660,-2.020)--cycle;
\filldraw[fill=blue!47,draw=black!75,line width=0.3pt,line join=round] (-0.221,-2.610)--(0.062,-2.646)--(0.660,-2.305)--(0.000,-1.136)--cycle;
\filldraw[fill=blue!15,draw=black!75,line width=0.3pt,line join=round] (-1.180,0.255)--(-2.545,0.169)--(-2.618,0.079)--(-2.468,0.060)--cycle;
\filldraw[fill=blue!17,draw=black!75,line width=0.3pt,line join=round] (-2.618,0.079)--(-2.545,0.169)--(-2.618,0.789)--(-2.656,0.574)--cycle;
\filldraw[fill=blue!54,draw=black!75,line width=0.3pt,line join=round] (-2.468,1.006)--(-2.618,0.789)--(-2.545,0.169)--(-1.180,0.255)--cycle;
\filldraw[fill=blue!33,draw=black!75,line width=0.3pt,line join=round] (-1.951,0.802)--(-2.412,1.008)--(-2.468,1.006)--(-1.180,0.255)--cycle;
\filldraw[fill=orange!48,draw=black!75,line width=0.3pt,line join=round] (2.174,1.657)--(2.051,1.810)--(1.500,1.829)--(1.020,0.682)--cycle;
\filldraw[fill=orange!47,draw=black!75,line width=0.3pt,line join=round] (2.385,0.768)--(2.357,1.443)--(2.174,1.657)--(1.020,0.682)--cycle;
\filldraw[fill=orange!27,draw=black!75,line width=0.3pt,line join=round] (2.051,1.810)--(2.174,1.657)--(2.357,1.443)--(2.295,1.534)--cycle;
\filldraw[fill=green!70!black!36,draw=black!75,line width=0.3pt,line join=round] (-1.451,2.119)--(-1.985,1.917)--(-2.078,1.700)--(-0.791,0.949)--cycle;
\filldraw[fill=green!70!black!36,draw=black!75,line width=0.3pt,line join=round] (-1.779,2.136)--(-1.985,1.917)--(-1.451,2.119)--(-1.347,2.343)--cycle;
\filldraw[fill=green!70!black!33,draw=black!75,line width=0.3pt,line join=round] (-1.451,2.119)--(-0.791,0.949)--(-1.229,2.269)--(-1.347,2.343)--cycle;
\filldraw[fill=orange!17,draw=black!75,line width=0.3pt,line join=round] (-2.556,-0.713)--(-2.518,-1.207)--(-2.456,-1.298)--(-2.484,-0.623)--cycle;
\filldraw[fill=orange!47,draw=black!75,line width=0.3pt,line join=round] (-2.484,-0.623)--(-2.456,-1.298)--(-2.273,-1.512)--(-1.119,-0.537)--cycle;
\filldraw[fill=orange!15,draw=black!75,line width=0.3pt,line join=round] (-1.779,-1.421)--(-1.119,-0.537)--(-2.273,-1.512)--(-2.213,-1.574)--cycle;
\filldraw[fill=orange!54,draw=black!75,line width=0.3pt,line join=round] (1.119,-0.599)--(2.407,-1.350)--(2.556,-1.133)--(2.484,-0.513)--cycle;
\filldraw[fill=orange!36,draw=black!75,line width=0.3pt,line join=round] (1.779,-1.769)--(2.313,-1.567)--(2.407,-1.350)--(1.119,-0.599)--cycle;
\filldraw[fill=orange!15,draw=black!75,line width=0.3pt,line join=round] (2.556,-1.133)--(2.407,-1.350)--(2.313,-1.567)--(2.518,-1.348)--cycle;
\filldraw[fill=blue!36,draw=black!75,line width=0.3pt,line join=round] (-0.655,2.674)--(-0.759,2.450)--(-0.161,2.791)--(-0.223,2.882)--cycle;
\filldraw[fill=blue!47,draw=black!75,line width=0.3pt,line join=round] (0.122,2.756)--(-0.161,2.791)--(-0.759,2.450)--(-0.099,1.281)--cycle;
\filldraw[fill=blue!25,draw=black!75,line width=0.3pt,line join=round] (0.381,2.428)--(0.233,2.783)--(0.122,2.756)--(-0.099,1.281)--cycle;
\filldraw[fill=blue!37,draw=black!75,line width=0.3pt,line join=round] (-0.161,2.791)--(0.122,2.756)--(0.233,2.783)--(-0.223,2.882)--cycle;
\filldraw[fill=green!70!black!16,draw=black!75,line width=0.3pt,line join=round] (-1.111,-2.270)--(-1.662,-2.251)--(-1.418,-2.527)--(-0.962,-2.626)--cycle;
\filldraw[fill=green!70!black!48,draw=black!75,line width=0.3pt,line join=round] (-0.630,-1.123)--(-1.784,-2.099)--(-1.662,-2.251)--(-1.111,-2.270)--cycle;
\filldraw[fill=green!70!black!25,draw=black!75,line width=0.3pt,line join=round] (-1.111,-2.270)--(-0.962,-2.626)--(-0.852,-2.598)--(-0.630,-1.123)--cycle;
\filldraw[fill=green!70!black!15,draw=black!75,line width=0.3pt,line join=round] (1.290,2.008)--(0.630,1.123)--(1.784,2.099)--(1.724,2.161)--cycle;
\filldraw[fill=green!70!black!36,draw=black!75,line width=0.3pt,line join=round] (1.290,2.008)--(1.024,2.535)--(0.852,2.598)--(0.630,1.123)--cycle;
\filldraw[fill=green!70!black!47,draw=black!75,line width=0.3pt,line join=round] (1.418,2.527)--(1.024,2.535)--(1.290,2.008)--(1.724,2.161)--cycle;
\filldraw[fill=green!70!black!27,draw=black!75,line width=0.3pt,line join=round] (1.724,2.161)--(1.784,2.099)--(1.662,2.251)--(1.418,2.527)--cycle;
\filldraw[fill=green!70!black!37,draw=black!75,line width=0.3pt,line join=round] (1.418,2.527)--(0.962,2.626)--(0.852,2.598)--(1.024,2.535)--cycle;
\filldraw[fill=orange!16,draw=black!75,line width=0.3pt,line join=round] (-0.381,-2.428)--(-0.233,-2.783)--(0.223,-2.882)--(0.261,-2.667)--cycle;
\filldraw[fill=orange!33,draw=black!75,line width=0.3pt,line join=round] (0.759,-2.450)--(0.099,-1.281)--(0.537,-2.600)--(0.655,-2.674)--cycle;
\filldraw[fill=orange!54,draw=black!75,line width=0.3pt,line join=round] (-0.381,-2.428)--(0.261,-2.667)--(0.537,-2.600)--(0.099,-1.281)--cycle;
\filldraw[fill=orange!15,draw=black!75,line width=0.3pt,line join=round] (0.655,-2.674)--(0.537,-2.600)--(0.261,-2.667)--(0.223,-2.882)--cycle;
\filldraw[fill=orange!17,draw=black!75,line width=0.3pt,line join=round] (-2.484,0.513)--(-2.456,0.986)--(-2.518,1.348)--(-2.556,1.133)--cycle;
\filldraw[fill=orange!36,draw=black!75,line width=0.3pt,line join=round] (-2.518,1.348)--(-2.213,1.420)--(-1.779,1.769)--(-2.313,1.567)--cycle;
\filldraw[fill=orange!27,draw=black!75,line width=0.3pt,line join=round] (-1.119,0.599)--(-2.273,1.154)--(-2.456,0.986)--(-2.484,0.513)--cycle;
\filldraw[fill=orange!46,draw=black!75,line width=0.3pt,line join=round] (-2.213,1.420)--(-2.273,1.154)--(-1.119,0.599)--(-1.779,1.769)--cycle;
\filldraw[fill=orange!39,draw=black!75,line width=0.3pt,line join=round] (-2.518,1.348)--(-2.456,0.986)--(-2.273,1.154)--(-2.213,1.420)--cycle;
\filldraw[fill=violet!15,draw=black!75,line width=0.3pt,line join=round] (2.484,0.623)--(2.556,0.713)--(2.407,0.732)--(1.119,0.537)--cycle;
\filldraw[fill=violet!36,draw=black!75,line width=0.3pt,line join=round] (1.119,0.537)--(2.407,0.732)--(2.313,0.988)--(1.779,1.421)--cycle;
\filldraw[fill=violet!36,draw=black!75,line width=0.3pt,line join=round] (2.518,1.207)--(2.313,0.988)--(2.407,0.732)--(2.556,0.713)--cycle;
\filldraw[fill=violet!47,draw=black!75,line width=0.3pt,line join=round] (2.313,0.988)--(2.518,1.207)--(2.213,1.574)--(1.779,1.421)--cycle;
\filldraw[fill=violet!27,draw=black!75,line width=0.3pt,line join=round] (2.518,1.207)--(2.456,1.298)--(2.273,1.512)--(2.213,1.574)--cycle;
\filldraw[fill=green!70!black!35,draw=black!75,line width=0.3pt,line join=round] (1.562,-1.497)--(1.385,-2.128)--(1.779,-2.136)--(2.023,-1.702)--cycle;
\filldraw[fill=green!70!black!33,draw=black!75,line width=0.3pt,line join=round] (1.562,-1.497)--(2.023,-1.702)--(2.078,-1.700)--(0.791,-0.949)--cycle;
\filldraw[fill=green!70!black!36,draw=black!75,line width=0.3pt,line join=round] (0.791,-0.949)--(1.229,-2.269)--(1.385,-2.128)--(1.562,-1.497)--cycle;
\filldraw[fill=green!70!black!15,draw=black!75,line width=0.3pt,line join=round] (1.779,-2.136)--(1.985,-1.917)--(2.078,-1.700)--(2.023,-1.702)--cycle;
\filldraw[fill=green!70!black!15,draw=black!75,line width=0.3pt,line join=round] (1.779,-2.136)--(1.385,-2.128)--(1.229,-2.269)--(1.347,-2.343)--cycle;
\filldraw[fill=blue!16,draw=black!75,line width=0.3pt,line join=round] (-2.295,-1.534)--(-2.051,-1.810)--(-1.500,-1.829)--(-1.989,-1.462)--cycle;
\filldraw[fill=blue!46,draw=black!75,line width=0.3pt,line join=round] (-1.020,-0.682)--(-2.091,-1.194)--(-1.989,-1.462)--(-1.500,-1.829)--cycle;
\filldraw[fill=blue!47,draw=black!75,line width=0.3pt,line join=round] (-2.385,-0.768)--(-2.295,-1.095)--(-2.091,-1.194)--(-1.020,-0.682)--cycle;
\filldraw[fill=blue!31,draw=black!75,line width=0.3pt,line join=round] (-1.989,-1.462)--(-2.091,-1.194)--(-2.295,-1.095)--(-2.295,-1.534)--cycle;
\filldraw[fill=blue!17,draw=black!75,line width=0.3pt,line join=round] (-2.385,-0.768)--(-2.357,-1.443)--(-2.295,-1.534)--(-2.295,-1.095)--cycle;
\filldraw[fill=teal!15,draw=black!75,line width=0.3pt,line join=round] (2.618,-0.789)--(2.656,-0.574)--(2.412,-1.008)--(2.468,-1.006)--cycle;
\filldraw[fill=teal!35,draw=black!75,line width=0.3pt,line join=round] (2.412,-1.008)--(2.656,-0.574)--(2.412,-0.298)--(1.951,-0.802)--cycle;
\filldraw[fill=teal!36,draw=black!75,line width=0.3pt,line join=round] (2.618,-0.079)--(2.468,-0.060)--(2.412,-0.298)--(2.656,-0.574)--cycle;
\filldraw[fill=teal!48,draw=black!75,line width=0.3pt,line join=round] (1.951,-0.802)--(2.412,-0.298)--(2.468,-0.060)--(1.180,-0.255)--cycle;
\filldraw[fill=orange!17,draw=black!75,line width=0.3pt,line join=round] (-0.660,2.305)--(0.000,1.136)--(-0.221,2.138)--(-0.394,2.292)--cycle;
\filldraw[fill=orange!27,draw=black!75,line width=0.3pt,line join=round] (0.000,1.136)--(0.660,2.020)--(0.062,2.193)--(-0.221,2.138)--cycle;
\filldraw[fill=orange!52,draw=black!75,line width=0.3pt,line join=round] (-0.221,2.138)--(0.062,2.193)--(0.000,2.555)--(-0.394,2.292)--cycle;
\filldraw[fill=orange!47,draw=black!75,line width=0.3pt,line join=round] (0.062,2.193)--(0.660,2.020)--(0.394,2.547)--(0.000,2.555)--cycle;
\filldraw[fill=orange!36,draw=black!75,line width=0.3pt,line join=round] (-0.394,2.292)--(0.000,2.555)--(-0.062,2.646)--(-0.660,2.305)--cycle;
\filldraw[fill=orange!37,draw=black!75,line width=0.3pt,line join=round] (0.394,2.547)--(0.221,2.610)--(-0.062,2.646)--(0.000,2.555)--cycle;
\filldraw[fill=green!70!black!50,draw=black!75,line width=0.3pt,line join=round] (-1.989,-0.322)--(-1.500,0.142)--(-2.051,0.479)--(-2.295,0.045)--cycle;
\filldraw[fill=green!70!black!15,draw=black!75,line width=0.3pt,line join=round] (-1.500,0.142)--(-1.989,-0.322)--(-2.091,-0.492)--(-1.020,0.020)--cycle;
\filldraw[fill=green!70!black!33,draw=black!75,line width=0.3pt,line join=round] (-1.020,0.020)--(-2.174,0.575)--(-2.051,0.479)--(-1.500,0.142)--cycle;
\filldraw[fill=green!70!black!31,draw=black!75,line width=0.3pt,line join=round] (-2.295,0.045)--(-2.295,-0.394)--(-2.091,-0.492)--(-1.989,-0.322)--cycle;
\filldraw[fill=green!70!black!17,draw=black!75,line width=0.3pt,line join=round] (-2.295,-0.394)--(-2.295,0.045)--(-2.357,0.407)--(-2.385,-0.066)--cycle;
\filldraw[fill=green!70!black!39,draw=black!75,line width=0.3pt,line join=round] (-2.051,0.479)--(-2.174,0.575)--(-2.357,0.407)--(-2.295,0.045)--cycle;
\filldraw[fill=blue!50,draw=black!75,line width=0.3pt,line join=round] (-0.962,1.508)--(-1.418,1.607)--(-1.662,1.173)--(-1.111,0.836)--cycle;
\filldraw[fill=blue!36,draw=black!75,line width=0.3pt,line join=round] (-1.290,1.884)--(-1.724,1.535)--(-1.418,1.607)--(-1.024,1.870)--cycle;
\filldraw[fill=blue!39,draw=black!75,line width=0.3pt,line join=round] (-1.724,1.535)--(-1.784,1.270)--(-1.662,1.173)--(-1.418,1.607)--cycle;
\filldraw[fill=blue!25,draw=black!75,line width=0.3pt,line join=round] (-0.630,0.714)--(-0.852,1.717)--(-0.962,1.508)--(-1.111,0.836)--cycle;
\filldraw[fill=blue!52,draw=black!75,line width=0.3pt,line join=round] (-0.962,1.508)--(-0.852,1.717)--(-1.024,1.870)--(-1.418,1.607)--cycle;
\filldraw[fill=blue!28,draw=black!75,line width=0.3pt,line join=round] (0.160,-0.937)--(-0.320,-2.084)--(-0.071,-1.876)--(0.073,-1.661)--cycle;
\filldraw[fill=blue!47,draw=black!75,line width=0.3pt,line join=round] (0.931,-1.484)--(0.160,-0.937)--(0.073,-1.661)--(0.361,-1.669)--cycle;
\filldraw[fill=blue!16,draw=black!75,line width=0.3pt,line join=round] (0.361,-2.108)--(-0.071,-1.876)--(-0.320,-2.084)--(0.323,-2.323)--cycle;
\filldraw[fill=blue!38,draw=black!75,line width=0.3pt,line join=round] (-0.071,-1.876)--(0.361,-2.108)--(0.361,-1.669)--(0.073,-1.661)--cycle;
\filldraw[fill=blue!35,draw=black!75,line width=0.3pt,line join=round] (0.361,-1.669)--(0.361,-2.108)--(0.755,-2.115)--(0.931,-1.484)--cycle;
\filldraw[fill=blue!15,draw=black!75,line width=0.3pt,line join=round] (0.361,-2.108)--(0.323,-2.323)--(0.598,-2.256)--(0.755,-2.115)--cycle;
\filldraw[fill=blue!35,draw=black!75,line width=0.3pt,line join=round] (2.023,0.144)--(1.779,0.419)--(1.385,0.156)--(1.562,-0.361)--cycle;
\filldraw[fill=blue!17,draw=black!75,line width=0.3pt,line join=round] (1.385,0.156)--(1.229,0.398)--(0.791,0.186)--(1.562,-0.361)--cycle;
\filldraw[fill=blue!36,draw=black!75,line width=0.3pt,line join=round] (2.023,0.144)--(2.078,0.381)--(1.985,0.638)--(1.779,0.419)--cycle;
\filldraw[fill=blue!28,draw=black!75,line width=0.3pt,line join=round] (0.791,0.186)--(1.229,0.398)--(1.347,0.651)--(1.451,1.071)--cycle;
\filldraw[fill=blue!47,draw=black!75,line width=0.3pt,line join=round] (1.451,1.071)--(1.347,0.651)--(1.779,0.419)--(1.985,0.638)--cycle;
\filldraw[fill=blue!52,draw=black!75,line width=0.3pt,line join=round] (1.779,0.419)--(1.347,0.651)--(1.229,0.398)--(1.385,0.156)--cycle;
\filldraw[fill=orange!16,draw=black!75,line width=0.3pt,line join=round] (-1.195,-1.275)--(-1.501,-1.347)--(-1.012,-1.714)--(-0.763,-1.506)--cycle;
\filldraw[fill=orange!31,draw=black!75,line width=0.3pt,line join=round] (-1.501,-1.347)--(-1.195,-1.275)--(-1.501,-0.908)--(-1.602,-1.078)--cycle;
\filldraw[fill=orange!33,draw=black!75,line width=0.3pt,line join=round] (-1.012,-0.444)--(-0.763,-1.068)--(-0.619,-1.291)--(-0.531,-0.567)--cycle;
\filldraw[fill=orange!38,draw=black!75,line width=0.3pt,line join=round] (-0.763,-1.506)--(-0.619,-1.291)--(-0.763,-1.068)--(-1.195,-1.275)--cycle;
\filldraw[fill=orange!50,draw=black!75,line width=0.3pt,line join=round] (-1.501,-0.908)--(-1.195,-1.275)--(-0.763,-1.068)--(-1.012,-0.444)--cycle;
\filldraw[fill=green!70!black!50,draw=black!75,line width=0.3pt,line join=round] (0.261,0.758)--(0.223,1.252)--(-0.233,1.351)--(-0.381,0.679)--cycle;
\filldraw[fill=green!70!black!15,draw=black!75,line width=0.3pt,line join=round] (0.099,0.557)--(0.537,0.768)--(0.261,0.758)--(-0.381,0.679)--cycle;
\filldraw[fill=green!70!black!47,draw=black!75,line width=0.3pt,line join=round] (0.655,1.021)--(0.759,1.441)--(0.161,1.614)--(0.223,1.252)--cycle;
\filldraw[fill=green!70!black!52,draw=black!75,line width=0.3pt,line join=round] (0.223,1.252)--(0.261,0.758)--(0.537,0.768)--(0.655,1.021)--cycle;
\filldraw[fill=green!70!black!52,draw=black!75,line width=0.3pt,line join=round] (0.161,1.614)--(-0.122,1.559)--(-0.233,1.351)--(0.223,1.252)--cycle;
\filldraw[fill=violet!50,draw=black!75,line width=0.3pt,line join=round] (0.323,-0.034)--(-0.320,-0.113)--(-0.071,-0.736)--(0.361,-0.529)--cycle;
\filldraw[fill=violet!38,draw=black!75,line width=0.3pt,line join=round] (0.073,-0.959)--(0.361,-0.967)--(0.361,-0.529)--(-0.071,-0.736)--cycle;
\filldraw[fill=violet!35,draw=black!75,line width=0.3pt,line join=round] (0.931,-0.782)--(0.755,-0.265)--(0.361,-0.529)--(0.361,-0.967)--cycle;
\filldraw[fill=violet!52,draw=black!75,line width=0.3pt,line join=round] (0.755,-0.265)--(0.598,-0.023)--(0.323,-0.034)--(0.361,-0.529)--cycle;
\node at (0,-3.35) {\small (d)};
\end{tikzpicture}
\caption{The tiling of Theorem~\ref{thm:C} for the icosidodecahedron, a canonical
$3$-polytope with $30$ vertices, with the configuration of
Corollary~\ref{thm:canon} ($B=O$). (a) The polytope. (b) The same polytope drawn
translucent: each facet is divided by the segments from its incentre $x_F$ to the
points of tangency $x_e$ of its edges (blue, dashed), and the inner edges of the tiling
are the segments $Ox_F$ (red, thick). (c) The cells assembled, with the cells of the near half removed:
the cells meet along the quadrilaterals $x_e\,x_F\,O\,x_{F'}$, all of which contain $O$;
neighbouring cells have different colours and shades. (d) Exploded view of all $30$
cells, coloured as in (c): each cell $C_v$ is translated by $0.8\,v$. Every vertex has degree $4$, so
every cell is a bipolar polytope with poles $v$ and $O$ whose vertex lattice is the face
lattice of a quadrilateral; it has $10$ vertices ($v$, four points $x_e$, four points
$x_F$, and $O$) and $8$ facets, four containing $v$ and four containing $O$. The cells are
orthogonally bipolar and are the cells of the power diagram of the Thales spheres
\cite[Remark~6.4]{Thales}.}
\label{fig:tiling}
\end{figure}
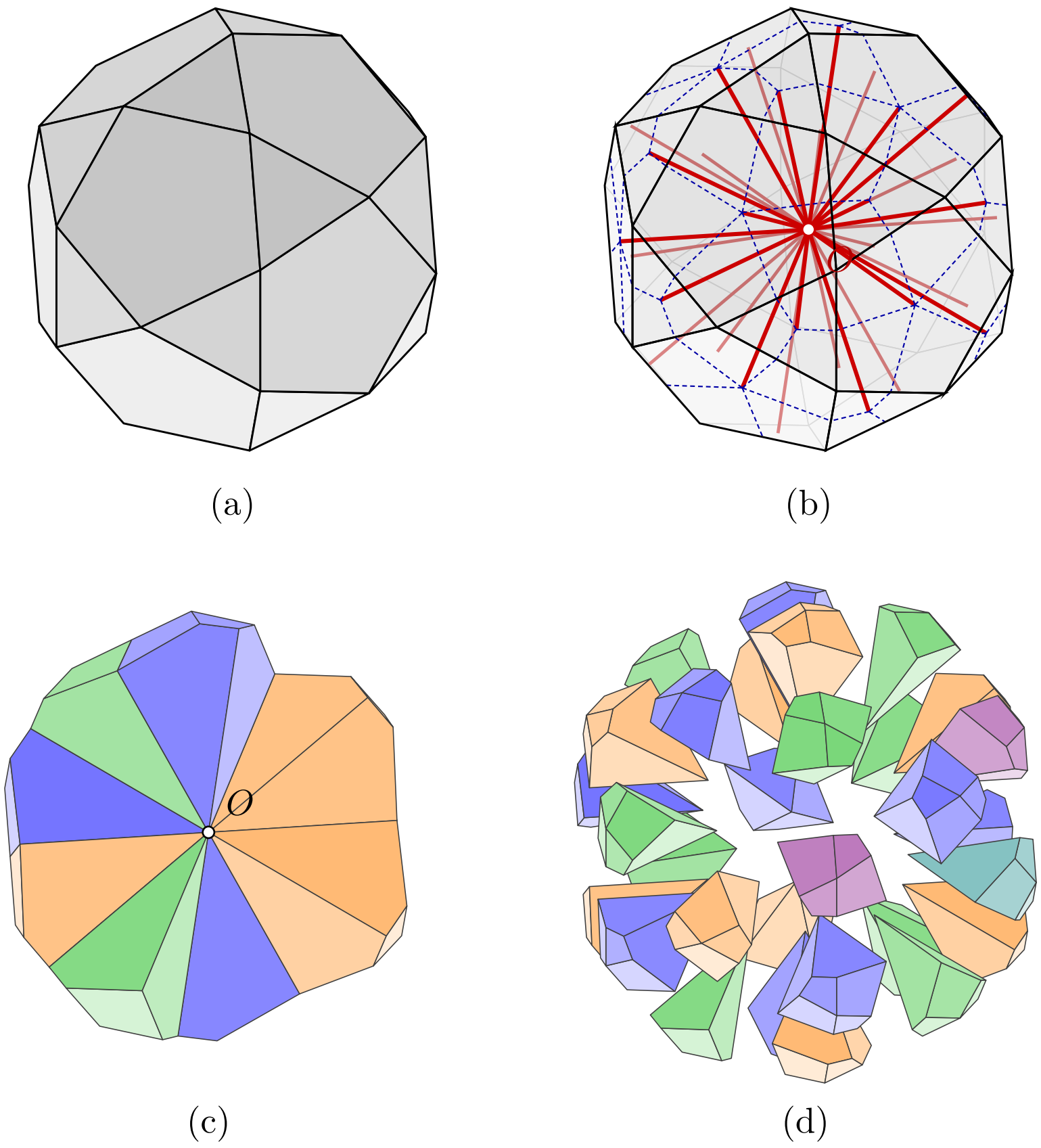

Figure~\ref{fig:tiling} shows the resulting tiling for the icosidodecahedron, whose
vertices all have degree~$4$.

\begin{remark}\label{rem:canon-orth}
In the canonical configuration, $\tau_e$ is the orthogonal projection of
$O$ onto the line of $e$ (as $O\tau_e\perp e$), and $x_F$ is that onto
$\aff F$. Hence the face $[v,\tau_e]\subseteq e$ of $C_v$ is perpendicular to
the face $[\tau_e,O]$, which lies in the plane $\Pi_e$ through $O$ orthogonal
to $e$; so each $C_v$ is an orthogonally bipolar polytope in the sense of
\cite{Thales}, inscribed in the sphere with diameter $[v,O]$. Thus a canonical $3$-polytope is tiled by orthogonally bipolar
polytopes \cite[Theorem~6.5]{Thales}. The proof of
Corollary~\ref{thm:canon} uses only that the edges touch a sphere whose
centre lies in $\inte P$.
\end{remark}

\begin{theorem}\label{thm:W3}
Every canonical $3$-polytope $Q$ (with centre $O$) satisfies the interval
condition (I). Consequently, for every $3$-polytope $Q$ the bipolar
polytope $W(\F(Q))$ is realized by a convex $4$-polytope, namely
$(\mathrm{Cay}(Q',Q'^\circ)-z)^\circ$ for a canonical realization $Q'$ of
$Q$ (Theorem~\ref{thm:cayley}). Together with
Theorem~\ref{thm:W2}, $W(\F(Q))$ has a geometric realization for every
polytope $Q$ of dimension at most three.
\end{theorem}
\begin{proof}
Use the canonical configuration $B=O$, $x_e=\tau_e$, $x_F$ (the incentre
of $F$) of Corollary~\ref{thm:canon}; it is regular. For a facet $F$ the
dual vertex is $F^*=x_F/\|x_F\|^2$, so $\relint N_Q(F)$ is the open ray
through $x_F\in\relint F$. For an edge $e=F\cap F'$, $N_Q(e)=\cone(e^*)$
with $e^*=[F^*,F'^*]$; this wedge lies in the plane $\Pi_e$ through $O$,
$x_F$, $x_{F'}$, which contains $\tau_e$, and by Step~4 of the proof of
Theorem~\ref{thm:C} the ray through $\tau_e$ lies strictly between the
rays through $x_F$ and $x_{F'}$. So $\relint N_Q(e)$ consists of the rays
through $(x_F,\tau_e)\cup\{\tau_e\}\cup(\tau_e,x_{F'})$, which lie in
$\relint F$, $\relint e$, $\relint F'$, and meets each of these three.
For a vertex $v$, $N_Q(v)=\cone(v^*)$, where $v^*$ is the polygon with the
vertices $F^*$, $F\ni v$, and the edges $e^*$, $e\ni v$. By the edge case,
the boundary rays of $\cone(v^*)$ are the rays through the zigzag
$Z=\tau_{e_1}x_{F_1}\tau_{e_2}x_{F_2}\cdots$ bounding $\operatorname{star}(v)$
(Steps~1--4 of the proof of Theorem~\ref{thm:C}); as both $\cone(v^*)$
and the cone over $\operatorname{star}(v)$ are convex (Step~4) and have
the same boundary, they coincide. The interior of $\operatorname{star}(v)$
is the union of $v$, the open segments $(v,\tau_e)$, $e\ni v$, and the
interiors of the quadrilaterals $R_{v,F}$, $F\ni v$; these lie in
$\relint\{v\}$, $\relint e$, $\relint F$, and each is met. Hence for every
face $H$, $\relint N_Q(H)$ meets $\relint\cone(G)$ exactly for the faces
$G\supseteq H$, which is (I). The last statements follow from (K1) and
Theorem~\ref{thm:cayley}.
\end{proof}

\subsection{A conjecture for general \texorpdfstring{$3$}{3}-polytopes}

By Corollary~\ref{thm:canon}, every combinatorial type of $3$-polytope has a
realization with a regular configuration. We expect that every realization has one.

\begin{conjecture}\label{conj:regular}
Every $3$-polytope has a regular configuration.
\end{conjecture}

\noindent\textbf{Computational evidence.} For $20$ random $3$-polytopes (simplicial ones
obtained as convex hulls of random points, and polytopes with $3$- to $6$-gonal faces:
perturbed cubes, truncated pentagonal prisms, pentagonal pyramids, intersections of
random half-spaces) we searched for a regular configuration by maximizing the smallest
of the margins in (C1), (C3), (C4). A configuration with positive margin was found in
every case (margins from $0.003$ for pentagonal pyramids with an almost degenerate base
to $0.19$). For each, the volumes of the cells $C_v$ and of their convex hulls agreed to
machine precision, and the cells tiled $P$, as predicted by Theorem~\ref{thm:C}.

The conjecture is invariant under projective transformations, which extends it from
the perfectly centered polytopes, whose orthogonal configurations are regular, to their
projective images.

\begin{proposition}[projective invariance]\label{prop:projinv}
Let $P$ be a $3$-polytope and $\varphi$ an admissible projective transformation for $P$.
\begin{enumerate}[label=(\alph*)]
\item If $(x_G)$ is a regular configuration of $P$, then $(\varphi(x_G))$ is a regular
      configuration of $\varphi(P)$, indexed by the corresponding faces.
\item If $P$ is projectively equivalent to a polytope that is perfectly centered about some
      point, then $P$ has a regular configuration, and it is tiled by convex bipolar
      polytopes $C_v$, $v\in V(P)$, with a common pole.
\end{enumerate}
\end{proposition}
\begin{proof}
(a) Each of the conditions (C1)--(C4) is expressed by relative interiors of faces,
collinearity, coplanarity and incidences of points of $P$, and convex hulls of subsets of
$P$, which are preserved by $\varphi$ (\S\ref{sec:notation}). (b) Let $P'=\varphi(P)$ be
perfectly centered about $B'$. Then $B'$ lies in the orthogonal kernel of $P'$, so the
orthogonal projections of $B'$ onto the faces of $P'$ form a regular configuration of $P'$
\cite[Remark~6.6]{Thales}, and by (a) its image under $\varphi^{-1}$, which is admissible for
$P'$, is a regular configuration of $P$. The tiling follows from Theorem~\ref{thm:C}.
\end{proof}

Perfect centering itself is not invariant under projective, or even affine, maps: a
parallelogram with sides $a,b$ and acute angle $\alpha$ is perfectly centered about some
point if and only if $a\cos\alpha<b$ and $b\cos\alpha<a$, since the perpendiculars to two
opposite sides must overlap; a sufficiently skew parallelogram
\cite[Figure~5]{GrunbaumAnti2} therefore fails, although it is an affine image of a
square. In numerical
experiments, random simplicial $3$-polytopes with $7$ to $11$ vertices were in no case
perfectly centered, but most of them had a perfectly centered projective image.

\section{Open problems}\label{sec:open}

The realization question for $W(\F(Q))$, Lindstr\"om's question whether every polytope
has an interval polytope, has been answered by Dobbins: there is a combinatorial
$4$-polytope $Q$ for which $W(\F(Q))$ has no geometric realization
(Remark~\ref{q:realize}). Every combinatorial type with a perfectly centered realization
has one (Theorem~\ref{thm:cayley}). The questions left open here
concern decompositions.

\paragraph{Bipolar tilings.} Does every polytope $P\subset\R^n$ have a tiling by convex
bipolar polytopes $C_v$, $v\in V(P)$, with poles $v$ and a common pole $B\in\inte P$? For
$n=3$ this holds for every canonical polytope (Corollary~\ref{thm:canon}), hence for a
suitable realization of every combinatorial type. In every dimension it holds for
every polytope with nonempty orthogonal kernel, i.e.\ with a perfectly centered position
of the origin \cite[Theorem~6.5]{Thales}. For $n\ge4$ it is open for the
combinatorial types without a perfectly centered realization, which exist by
\cite[Theorem~3.9]{Dobbins}.

\paragraph{Derived questions.}
\begin{enumerate}[label=\arabic*.,leftmargin=*]
\item \emph{Regular configurations (Conjecture~\ref{conj:regular}).} Does every
      $3$-polytope have a regular configuration? By Theorem~\ref{thm:C} this would
      answer the question above for $n=3$ in every realization. It is known for
      canonical $3$-polytopes (Corollary~\ref{thm:canon}) and, more generally, for
      all $3$-polytopes projectively equivalent to perfectly centered ones
      (Proposition~\ref{prop:projinv}). Is every $3$-polytope projectively equivalent
      to a perfectly centered one?
\item \emph{The role of (C4) (Remark~\ref{rem:C3C4}).} Is there a configuration
      satisfying (C1)--(C3) with $B\in\conv W_v$ for which the cell $C_v$ is not
      convex?
\item \emph{Configurations in higher dimensions.} For $n\ge4$, find conditions on a
      configuration of an $n$-polytope under which the cells $C_v$ are convex bipolar
      polytopes that tile $P$, i.e.\ an analogue of Theorem~\ref{thm:C}.
\end{enumerate}

\section*{Acknowledgements}
The author thanks Martin Winter for valuable advice, in particular for pointing out the
connection with interval polytopes and antiprisms and the work of Lindstr\"om, Broadie,
Fukuda and Weibel, Dobbins, and Gleason and Hubard. The author acknowledges the support by the Specific Research Project no.\
FSI-S-26-9581 of Brno University of Technology. The funding source had no involvement
in the study, in the writing of the article, or in the decision to submit it for
publication.

\section*{Declaration of competing interest}
The author declares that he has no known competing financial interests or personal
relationships that could have appeared to influence the work reported in this paper.

\section*{Data availability}
No data were generated or analysed in this study. The numerical experiments mentioned
in Section~\ref{sec:geom} only illustrate the results; the
scripts are available from the author on request.

\section*{Declaration of generative AI and AI-assisted technologies in the manuscript
preparation process}
During the preparation of this paper the author used Claude (Anthropic) as an
assistant: for checking and discussing proofs, searching the literature, suggesting
reformulations and restructuring, carrying out numerical experiments, and preparing the
TikZ code of the figures from exact geometric data. All results, proofs, references and
figures were verified by the author, who takes full responsibility for the content.

\end{document}